\documentclass[oneside,reqno]{amsart}

\usepackage{preamble}

\definecolor{phaseA-tri}{HTML}{305753}
\definecolor{phaseA-tri-ext}{HTML}{A0BBB8}
\definecolor{phaseA-face}{HTML}{F1FFFD}
\definecolor{phaseB-tri}{HTML}{D68E8A}
\definecolor{phaseB-tri-ext}{HTML}{EEC2BF}
\definecolor{phaseB-face}{HTML}{F9E3E2}
\colorlet{seam-edge}{black}
\colorlet{seam-d}{black}
\colorlet{seam-u}{black}
\colorlet{seam-gap}{white}
\colorlet{plane-line}{black}

\definecolor{r-line}{HTML}{BF3B3F}
\tikzset{
  oldline/.style={line width=.3pt, black},
  newline/.style={line width=.3pt, black},
  rline/.style={line width=1.3pt, r-line, dash pattern=on 4pt off 2pt},
  face/.style={fill=gray!20},
  rface/.style={pattern={Lines[angle=-45, distance=5pt, line width=1.8pt]}, pattern color=phaseA-tri-ext},
  starpt/.style={fill=black},
  lbl/.style={font=\scriptsize, inner sep=1.4pt},
  rimlbl/.style={font=\small, inner sep=2pt},
  starlbl/.style={font=\small, inner sep=3pt},
}

\tikzset{
  g-b-int-s3/.style={fill=phaseA-tri},
  g-b-ext-s3/.style={fill=phaseA-tri-ext},
  g-w/.style={fill=phaseA-face},
  g-s4-w/.style={fill=phaseA-face},
  g-int-s3-w/.style={fill=phaseB-tri},
  g-ext-s3-w/.style={fill=phaseB-tri-ext},
  g-b-int/.style={fill=phaseB-face},
  g-b-ext/.style={fill=phaseB-face},
  g-b-int-s4/.style={fill=phaseB-face},
  g-b-ext-s4/.style={fill=phaseB-face},
}

\tikzset{
  zoom-mark/.style={draw=black, line width=0.5pt},
  zoom-frame/.style={draw=black, line width=0.5pt},
}

\graphicspath{{figures/}}

\title[Non-simple Arrangements with Many Triangles]{Non-simple Arrangements of Lines\\ and Pseudolines with Many Triangles}
\author{Roman Parpalak}
\author{Denis Utkin}
\thanks{Email: \texttt{roman@parpalak.com} (R.\,Parpalak),
  \texttt{ud1@mail.ru} (D.\,Utkin).}
\date{November 2, 2026}

\begin{document}

\begin{abstract}
We study affine arrangements of an even number $n$ of lines and pseudolines with many
triangles. For simple ones, Blanc's bound limits the number of triangular faces
to $\lfloor n(n-5/2)/3\rfloor$.
Adapting the Bartholdi--Blanc--Loisel doubling to the even case, we construct
an infinite series of arrangements of $n=6\cdot2^t+2$ ($8,14,26,50,\ldots$) straight lines with two
triple points and $(n(n-5/2)+1)/3$ triangles, one more than the simple bound
allows.

We obtain exhaustive enumerations of triangle-maximal pseudoline
arrangements, simple or not, for even $n\le12$. For $n=6,8,10,12,14$ the maximum
number of triangles, for both lines and pseudolines, is $7,15,25,38,54$; for
$n=8,12,14$ every optimal arrangement has multiple points. All optimal arrangements
for even $n\le12$ are realizable by straight lines, except two arrangements of 10
pseudolines, which are not stretchable. For the odd value $n=11$, no
arrangement with $33$ triangles exists even with multiple points.

We also obtain even arrangements by adding a line to odd ones. For
$n\equiv2\pmod 6$ this again gives one triangle more than the simple bound, in
pseudolines up to $n=44$ and in straight lines for $n=20$ and $n=32$.
Doubling these pseudoline arrangements gives members of the series
$18\cdot2^t+2$ and $30\cdot2^t+2$ with the same excess, up to $n=146$ and
$n=242$, respectively; those with $n=38$ and $n=74$ are realizable by straight lines.

These results lead us to conjecture that an arrangement of an even number $n$
of pseudolines has at most $\lfloor n(n-2)/3\rfloor-\lfloor n/6\rfloor$
triangles.
\end{abstract}

\maketitle

\section{Introduction}\label{sec:intro}

\subsection{Background}

Line and pseudoline arrangements with many triangular faces have been studied by several authors. In the affine setting this is the Kobon triangle problem~\cite{fujimura78,a006066}, which asks for the maximum number of triangular faces determined by $n$ lines in the plane. A brief survey of the literature and the basic notions is given in our previous paper~\cite{pu26b}, which deals with enumeration and classification algorithms for triangle-maximal affine arrangements of an odd number of pseudolines. There we also proved that the enumeration algorithms miss no arrangements, obtained the exhaustive enumeration for $n=27$ pseudolines, and found examples up to $n=93$.

The even case turned out to be harder: its arrangements have a more complicated structure, and, unlike in the odd case, no known upper bound for arbitrary arrangements matches the known examples for any $n\ge6$. Consider the bounds first. Roudneff~\cite{roudneff96} proved that the number $p_3$ of triangles in a projective arrangement $\mathcal{A}$ of $N$ pseudolines satisfies
\begin{equation}\label{eq:p3-bound}
p_3(\mathcal{A})\le\left\lfloor\frac{N(N-1)}{3}\right\rfloor,\qquad N\ge 9.
\end{equation}

Felsner and Kriegel claim that the proof can be adapted to the affine case~\cite[Theorem~2(3)]{fk99}. This gives a bound on the number $a_3$ of triangles in an affine arrangement $\mathcal{A}$ of $n$ pseudolines:
\begin{equation}\label{eq:a3-bound}
a_3(\mathcal{A})\le\left\lfloor\frac{n(n-2)}{3}\right\rfloor.
\end{equation}

Both bounds are stated for arbitrary pseudoline arrangements. For \emph{simple} arrangements, which have no point where three or more pseudolines cross (\emph{multiple points}), \eqref{eq:a3-bound} follows from an elementary segment count~\cite[Subsection~2.3]{pu26b}. The other $n-1$ pseudolines cut each pseudoline into $n-2$ bounded segments, and each segment is a side of at most one triangle. Let $u(\mathcal A)$ be the number of \emph{unused segments}, which are sides of no triangle. Counting all segments gives
\begin{equation}\label{eq:simple-segm}
3a_3(\mathcal{A})+u(\mathcal{A})=n(n-2),
\end{equation}
which implies \eqref{eq:a3-bound}.

For odd $n$, some simple arrangements have no unused segments and attain \eqref{eq:a3-bound} without rounding. Such arrangements are called \emph{perfect}.
By divisibility, they can exist only for $n\equiv 3,5\pmod 6$.
Their triangles alternate along every pseudoline, and the alternation is consistent across the whole arrangement: in a checkerboard coloring, all triangles are black and all other faces are white~\cite[Lemma~2.2]{pu26b}.

In a simple arrangement with even $n$, the alternation cannot propagate: for every pseudoline $\ell$, some other pseudoline has an unused segment with an endpoint on $\ell$~\cite[Proposition~2.3]{blanc11}. Each such segment serves at most two pseudolines, so $u\ge n/2$. This gives Blanc's bound~\cite[Corollary~2.4]{blanc11}:
\begin{equation}\label{eq:blanc-even}
a_{3}(\mathcal{A}_\text{simple})\le\left\lfloor\dfrac{n\left(n - \frac52\right)}{3}\right\rfloor=\begin{cases}
\dfrac{n\left(n - \frac52\right)}{3}, & n\equiv0,4\pmod 6 \\
\dfrac{n\left(n - \frac52\right) - 2}{3}. &  n\equiv2\pmod 6
\end{cases}
\end{equation}
Compared with \eqref{eq:a3-bound}, the bound \eqref{eq:blanc-even} loses $n/6$ triangles, up to rounding, to the $n/2$ unavoidable unused segments.

In the even case, however, one can no longer restrict attention to simple arrangements.
Two examples were already known in which multiple points yield an extra triangle.
The arrangement of $n=8$ lines in \cref{fig:seed} has 15 triangles, while \eqref{eq:blanc-even} allows only 14.
For $n=12$, an arrangement found by Kabanovich~\cite{khodulev99} attains \eqref{eq:blanc-even} with 38 triangles, although no simple arrangement has 38 triangles (see \cref{tab:triple-points}).
We found a third example, $n=14$ with 54 triangles, and it was the starting point of this work. Maiorana~\cite{maiorana26} found it independently while this paper was in preparation.
Still, none of these examples attains the general bound \eqref{eq:a3-bound}, and the literature gives no intermediate bound for non-simple arrangements.

Thus in the even case triangles cannot be packed with the unbroken alternation of perfect arrangements: parity or multiple points break it.
The breaks are reflected in the gap between \eqref{eq:a3-bound} and the known examples, which grows with $n$.
They also affect the computer search: in the perfect case a branch can be discarded at the first break of the alternation, whereas in the even case the search has to work with a budget on the number of triangles (cf.~\cite[Section~3]{pu26b}).

A further difficulty of the even case is that the known infinite series of even line arrangements are obtained from odd ones~\cite{for,bbl} by adding a line far from all intersection points~\cite{bbl,blanc11,savchuk25}, and they do not exceed \eqref{eq:blanc-even}. No construction was previously known that takes an even arrangement as input. The known odd series are $2\cdot2^t+1$, $6\cdot2^t+1$, $14\cdot2^t+1$~\cite{bbl} and $18\cdot2^t+1$, which we found earlier~\cite{pu26}. Iterative doubling~\cite{bbl}, on which these series are based, takes a suitable arrangement of an odd number $n$ of lines and produces an arrangement of $2n-1$ lines.

\subsection{Overview of results}

This paper is concerned with maximizing the number $a_{3}(\mathcal{A})$ of triangular faces in the even case, over simple and non-simple affine arrangements of pseudolines and of straight lines.

We adapt iterative doubling to an even number of lines and apply it to an arrangement of $n=8$ lines. The adaptation temporarily removes one line, applies doubling as many times as needed, and puts the line back. This procedure gives an infinite series of arrangements $\mathcal{C}_t$ of $n=6\cdot2^t+2$ lines with one triangle more than \eqref{eq:blanc-even} allows:
\[
a_3(\mathcal{C}_t)=\frac{n\left(n-\frac52\right)+1}{3}
=\frac{n\left(n-\frac52\right)-2}{3}+1.
\]
The series contains arrangements of $n=8$, $14$, $26$, $50,\ldots$ lines with $a_3=15$, $54$, $204$, $792,\ldots$ triangles. This series (\cref{thm:series}, proved in \cref{sec:series}) constitutes the main part of our results.

\Cref{sec:enum} describes a computer search for pseudoline arrangements, possibly with multiple points. Its results include exhaustive enumerations for $n=6, 8, 10, 12$ and a proof that no arrangement of $n=14$ pseudolines has $a_3>54$. This verifies that the maximum numbers of triangles for these $n$ are 7, 15, 25, 38, 54. The result extends to straight lines: these maxima are attained by the line arrangements mentioned above and by other known ones. We also checked that no arrangement of $n=11$ pseudolines has $a_3=33$, simple or not.

\Cref{sec:augment} describes augmentation, which produces arrangements of an even number of pseudolines from odd ones.
Augmenting perfect arrangements of $n-1$ pseudolines gives arrangements that attain \eqref{eq:blanc-even} for $n\equiv0,4\pmod 6$.
For $n\equiv2\pmod 6$, augmentation can exceed \eqref{eq:blanc-even} by one, like the series above. We found such augmentations for $n=8,14,20,26,32,38,44$.

We found realizations by straight lines of all arrangements, up to isomorphism, given by the exhaustive enumerations of \cref{sec:enum} for even $n$, except two arrangements of $n=10$ pseudolines. These two are not stretchable (verified by the method of biquadratic final polynomials~\cite{br90}). Line realizations of augmented arrangements (\cref{sec:augment}) of $n=20$ and $n=32$ pseudolines with 117 and 315 triangles were found as well. All these line arrangements are collected in \cref{app:figures}.

We also applied the doubling of \cref{sec:series} to pseudolines. Starting from the augmented arrangements with $n=20$ and $n=32$ (\cref{fig:n20-collapse,fig:n32}), we obtained the members with $t\le3$ of the pseudoline series $18\cdot2^t+2$ and $30\cdot2^t+2$, up to $n=146$ and $n=242$. Each of them exceeds \eqref{eq:blanc-even} by one. The members with $n=38$ and $n=74$, with 450 and 1764 triangles, are also realizable by straight lines (\cref{fig:n38,tab:n74-lines}).

These results lead us to conjecture that every arrangement of an even number $n$ of pseudolines satisfies
\begin{equation}\label{eq:conj-bound}
a_{3}\le\left\lfloor\dfrac{n(n - 2)}{3}\right\rfloor-\left\lfloor\dfrac{n}{6}\right\rfloor=\begin{cases}
\dfrac{n(n - \frac52)}{3}, & n\equiv0,4\pmod 6 \\
\dfrac{n(n - \frac52) + 1}{3}. &  n\equiv2\pmod 6
\end{cases}
\end{equation}
\Cref{sec:discussion} discusses the results and the conjecture.

\section{An infinite series of arrangements of \texorpdfstring{$6\cdot2^t+2$}{6·2ᵗ + 2} straight lines}\label{sec:series}

We now construct an infinite series of arrangements of an even number of lines.
The starting arrangement of the series has $N_0=8$ lines and two triple points.
The other arrangements are obtained from it by temporarily removing one line, applying the Bartholdi--Blanc--Loisel doubling~\cite{bbl} the required number of times, and restoring the removed line.

In the notation of~\cite{bbl}, $n$ has a special meaning: it is the number of added lines.
To avoid confusion, in this section, unlike the others, the number of lines in an arrangement is always denoted by $N$.
An arrangement $\mathcal A$ of $N$ straight lines in the affine plane is a set of $N$ distinct lines, no two of which are parallel.
The \emph{multiplicity} of a point is the number of lines of the arrangement passing through it.
A point of multiplicity 2 is called \emph{simple}, and a point of higher multiplicity is called \emph{multiple}.
An arrangement is \emph{simple} if all its intersection points are simple.

The lines divide the plane into \emph{faces}.
A face is \emph{triangular} if it is bounded and its sides lie on exactly three lines of the arrangement.
We denote by $a_3(\mathcal A)$ the number of triangular faces of the arrangement $\mathcal A$.

\begin{theorem}\label{thm:series}
For every integer $t\ge0$ there exists an arrangement $\mathcal{C}_t$ of
\begin{equation}\label{eq:N-t}
N_t=6\cdot2^t+2
\end{equation}
lines with two triple points and exactly
\begin{equation}\label{eq:a3-series}
a_3(\mathcal{C}_t)=\frac{N_t(N_t-\frac52)+1}{3}
\end{equation}
triangular faces.
\end{theorem}

\subsection{The starting arrangement}

Fix an arbitrary integer $t\ge0$ for which the theorem is to be proved. For the intermediate values $s$, $0\le s\le t$, set
\[
n_s=6\cdot2^s.
\]
Consider the eight lines given by the equations
\begin{equation}\label{eq:base-lines}
\begin{array}{ll}
R:y=x-\delta,
&Y_0:y=0,\\[4pt]
L_1:y=-\dfrac12\left(x+\sqrt3\right),
&L_2:y=-\eta\Bigl(x+\dfrac1{\sqrt3}\Bigr),\\[6pt]
L_3:y=\eta\Bigl(x-\dfrac1{\sqrt3}\Bigr),
&L_4:y=\dfrac25\left(x-\sqrt3\right),\\[8pt]
L_5:y=u(\varepsilon,\eta)(x+\varepsilon),
&L_6:y=v(\varepsilon)(x-\varepsilon).
\end{array}
\end{equation}
These equations have three free parameters. We tie them to the target value $t$:
\begin{equation}\label{eq:params}
\varepsilon=\frac1{4n_t},
\qquad
\eta=3\sqrt3\,n_t^6,
\qquad
\delta=2n_t^6.
\end{equation}
The remaining two dependent parameters $u(\varepsilon,\eta)$ and $v(\varepsilon)$ are defined as follows:
\begin{equation}\label{eq:star-slopes}
u(\varepsilon,\eta)=-\frac{4\eta}{2\eta-3+\varepsilon\sqrt3(2\eta+1)},
\qquad
v(\varepsilon)=\frac{4\sqrt3}{\sqrt3+9\varepsilon}.
\end{equation}

\begin{figure}[t]
\centering
\resizebox{\ifdim\width>\linewidth\linewidth\else\width\fi}{!}{\input{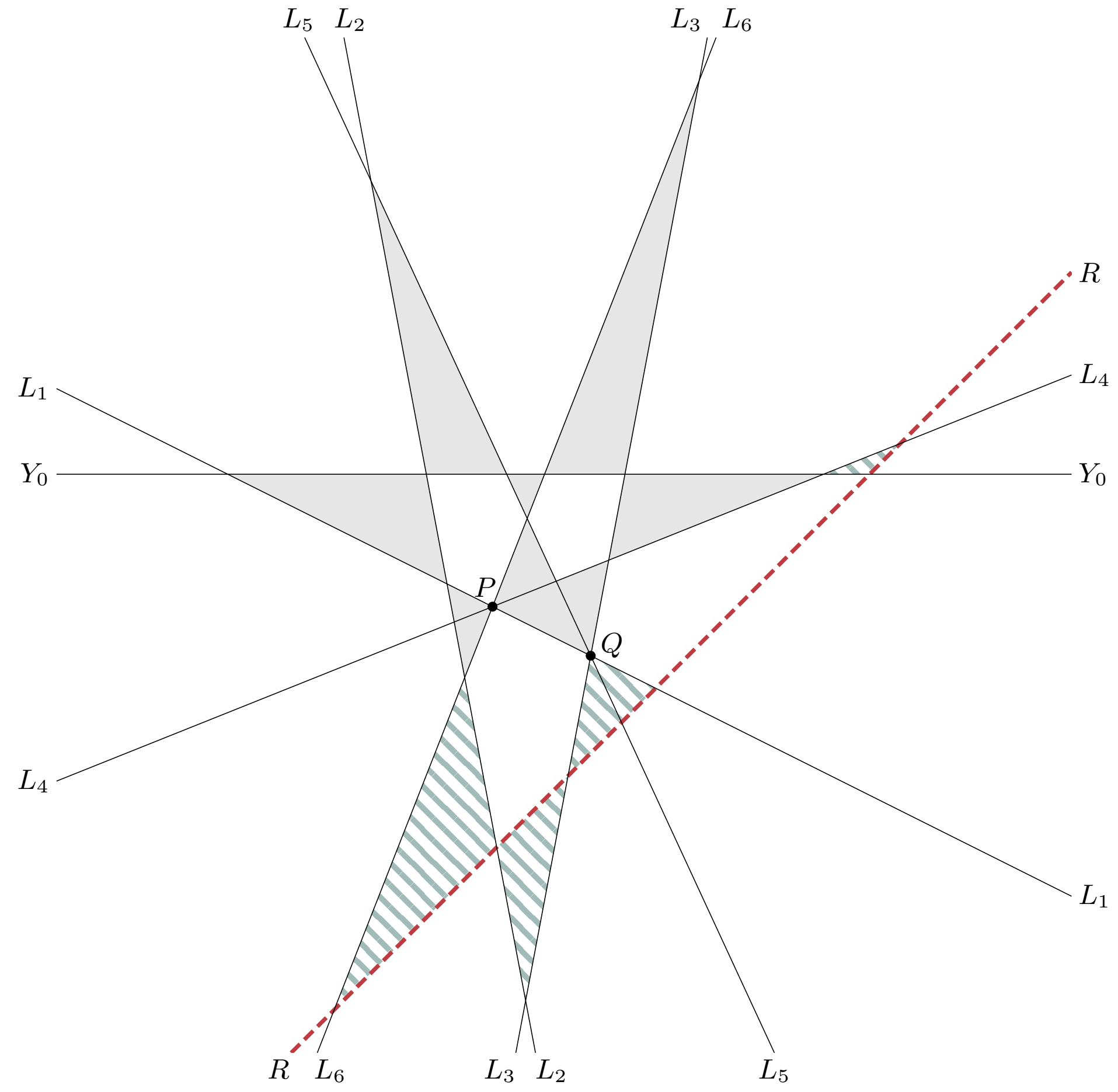}}
\caption{The starting arrangement \eqref{eq:base-lines}, with 15 triangular faces (5 of them $R$-triangles, hatched) and two triple points $P$ and $Q$. The values of $\varepsilon,\eta,\delta$ differ from \eqref{eq:params} for readability, but the order of intersection points on each line is the same.}
\label{fig:seed}
\end{figure}

With these equations, the eight lines form an arrangement with 15 triangular faces (\cref{fig:seed}).
The line $Y_0$ is horizontal, and the lines $L_1,\ldots,L_6$ intersect it at the points $x=\pm\varepsilon,\tan(\pm\pi/6),\tan(\pm2\pi/6)$,
as the doubling requires~\cite[Prop.~3.1]{bbl}. The parameter $\varepsilon$ satisfies the condition under which the doubling can be iterated~\cite[Remark~3.2]{bbl}.
The slopes of $L_1,L_4$ are fixed and equal to $-1/2, 2/5$, respectively. The slopes of $L_2,L_3$ are set by \eqref{eq:params} and depend on $t$.
The slopes \eqref{eq:star-slopes} of $L_5,L_6$ are chosen so that three lines meet at each of two points:
\[
P=L_1\cap L_4\cap L_6,
\qquad
Q=L_1\cap L_3\cap L_5.
\]
The parameters $\eta,\delta$ ensure that after $t$ iterations of the doubling all new intersection points lie on one side of the line $R$.
This line plays a special role:
removing it before the doubling makes the doubling possible, and afterwards the line is restored.
We call a triangular face with a side on $R$ an \emph{$R$-triangle}, and write $r(\cdot)$ for the number of $R$-triangles of an arrangement.

All lines in \eqref{eq:base-lines} have finite slope, that is, none of them is vertical.
Therefore the order of the intersection points along any line of the arrangement is uniquely determined by the order of their $x$-coordinates, and below the words ``to the right'' and ``rightmost'' refer to comparing $x$-coordinates.

\begin{lemma}\label{lem:seed}
The starting arrangement has exactly 15 triangular faces, and 5 of them are $R$-triangles.
After the line $R$ is removed, the arrangement
\[
\mathcal A_0=\{Y_0,L_1,\ldots,L_6\}
\]
has exactly two triple points $P$ and $Q$ (and no other multiple points) and 10 triangular faces.
\end{lemma}

\begin{proof}
\Cref{fig:seed} illustrates the statement, which is verified in \cref{app:seed-faces}.
\end{proof}

We construct the arrangements $\mathcal A_1,\ldots,\mathcal A_t$ recursively. The arrangement $\mathcal A_0$ has been defined above.
The arrangement $\mathcal A_{s+1}$ is obtained from $\mathcal A_s$ by adding the $n_s$ lines $M_j$ given by the conditions of~\cite[Prop.~3.1]{bbl} with $n=n_s$.
By construction, $\mathcal A_s$ consists of $n_s+1$ lines.
The line $R$ is not used in the doublings; after $t$ doublings, we add it to $\mathcal A_t$. This gives the required arrangement $\mathcal C_t$:
\[
\mathcal C_t=\mathcal A_t\cup\{R\},
\]
containing $N_t=n_t+2$ lines \eqref{eq:N-t}, for which we prove the statement of the theorem.

\subsection{The doubling step}

\begin{figure}[t]
\centering
\resizebox{\ifdim\width>\linewidth\linewidth\else\width\fi}{!}{\input{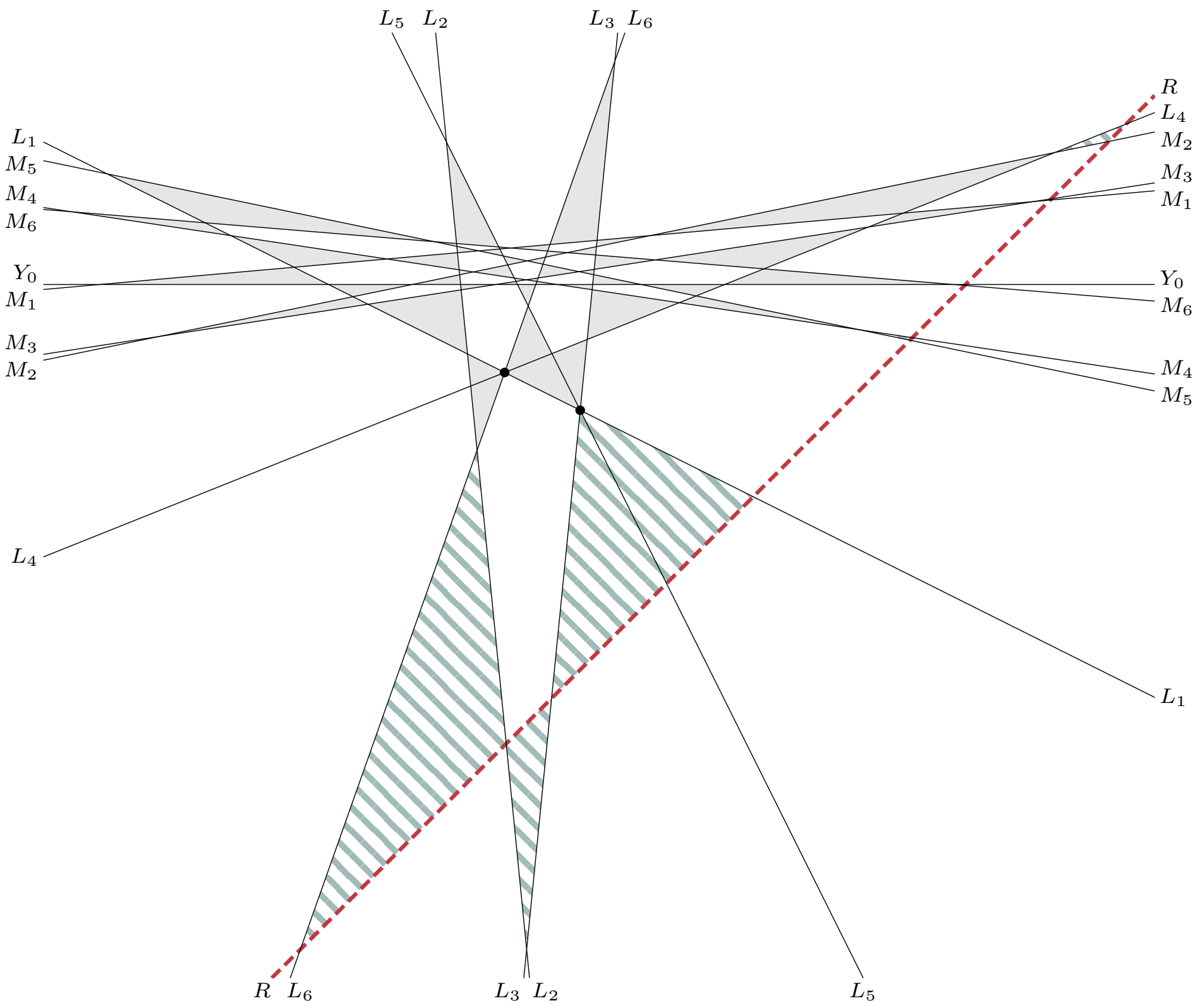}}
\caption{The starting arrangement after the doubling step \eqref{eq:new-lines}, with 54 triangular faces (8 of them $R$-triangles); the triple points are preserved.
The slopes $\mu_j$ are larger in absolute value than in \eqref{eq:new-lines}, since with the actual values the strip \eqref{eq:strip} would be indistinguishable at this scale, but the order of intersection points on each line is the same.}
\label{fig:double}
\end{figure}

Let $n$ be even, and let the arrangement
\[
\mathcal A=\{Y_0,L_1,\ldots,L_n\}
\]
be obtained at one of the steps of the construction. Choose coordinates in which
\[
Y_0:y=0,
\qquad
L_i:y=m_i(x-a_i).
\]
The added lines are those of~\cite[Prop.~3.1]{bbl}:
\begin{equation}\label{eq:new-lines}
\begin{gathered}
M_j:y=\mu_j\left(x-b_j\right),
\qquad
b_j=\tan\beta_j,
\qquad
\beta_j=-\frac\pi2+\left(j-\frac12\right)\frac\pi n,\\
\mu_j=-\frac{m_{\min}}{n^{10}}
\left(\sin(2\beta_j)+\frac1{n^6b_j}\right),
\qquad
m_{\min}=\min_{1\le i\le n}|m_i|,
\qquad
j=1,\ldots,n.
\end{gathered}
\end{equation}
The factor $-1$ in $\mu_j$ is the sign $\sigma$ of~\cite[Prop.~3.1]{bbl}, which equals $-1$ when $L_{n-1}$ and $L_n$ meet below $Y_0$, as $L_5$ and $L_6$ in \eqref{eq:base-lines} do.
The result of the first doubling step is shown in \cref{fig:double}.
From now on, the indices $i$ and $j$ always run over $1,\ldots,n$. The subscript of $L_i$, $a_i$, $m_i$ refers to an old line, and the subscript of $M_j$, $b_j$, $\beta_j$, $\mu_j$ refers to a new one.
Below we compare $m_{\min}$ for different arrangements, so where confusion is possible, we add the arrangement as an argument: $m_{\min}(\mathcal A)$.

The lines $L_1,\ldots,L_{n-2}$ are numbered in increasing order of their $x$-intercepts, so that $a_1<\cdots<a_{n-2}$. The lines $L_{n-1},L_n$ form the \emph{special pair}.
The doubling step imposes the following conditions on $\mathcal A$:
\begin{enumerate}[label=\textup{(D\arabic*)}]
\item the $x$-intercepts form the grid
\begin{equation}\label{eq:x-grid}
\begin{gathered}
\{a_1,\ldots,a_{n-2}\}
=\Bigl\{\tan\frac{k\pi}{n}:1\le |k|\le\frac n2-1\Bigr\},
\\
-\frac1n<a_{n-1}<0<a_n<\frac1n,
\end{gathered}
\end{equation}
in particular, $a_{n-2}=\cot(\pi/n)$;
\item the slopes of the special pair satisfy $m_{n-1}<0<m_n$;
\item each of the $n-1$ bounded segments of $Y_0$ is a side of exactly one triangular face of $\mathcal A$, and the faces on consecutive segments lie on opposite sides of $Y_0$.
\end{enumerate}
Under (D1) and (D2), the point $L_{n-1}\cap L_n$ lies below $Y_0$, and its $x$-coordinate lies in the interval $(a_{n-1},a_n)$.

\begin{lemma}[doubling step]\label{lem:doubling}
Let $n\ge2$ be even, and let $\mathcal A$ be a not necessarily simple arrangement satisfying \textup{(D1)--(D3)}. Then
\[
\mathcal B=\mathcal A\cup\{M_1,\ldots,M_n\}
\]
preserves all multiple points of $\mathcal A$ together with their multiplicities, has no new multiple points, and satisfies
\[
a_3(\mathcal B)=a_3(\mathcal A)+n^2.
\]
The arrangement $\mathcal B$ satisfies \textup{(D3)} with $2n$ in place of $n$. If, in addition,
\begin{equation}\label{eq:special-slopes}
|a_{n-1}|<\frac1{2n},\quad |a_n|<\frac1{2n},
\end{equation}
then $\mathcal B$ (after relabeling its lines) satisfies \textup{(D1)--(D3)} with $2n$ in place of $n$.
Finally,
\begin{equation}\label{eq:m-min}
m_{\min}(\mathcal B)<\frac{2\,m_{\min}(\mathcal A)}{n^{10}}<m_{\min}(\mathcal A).
\end{equation}
\end{lemma}

\begin{proof}
The computations in items~1--4 of the proof of~\cite[Prop.~3.1]{bbl} rely only on the equations of the lines: on condition \eqref{eq:x-grid}, on the signs of the slopes of the special pair,\footnote{However, item~4 of~\cite[Prop.~3.1]{bbl} needs a correction. There the signs of the slopes of the special pair are derived from the position of the point $L_{n-1}\cap L_n$ relative to $Y_0$. But the position of this point does not determine them: the two lines can also meet on that same side of $Y_0$ when their slopes have the same sign (then the $x$-coordinate of their intersection point lies outside $(a_{n-1},a_n)$). Therefore the signs $m_{n-1}<0<m_n$ are included in the hypotheses of the lemma as (D2), and everything the proof of~\cite[Prop.~3.1]{bbl} derives from the signs, namely $x(M_j\cap L_{n-1})<0<x(M_j\cap L_n)$, follows from the inequalities \eqref{eq:special-x} proved below.}
and on formulas \eqref{eq:new-lines}.
These computations do not use the simplicity of the arrangement. Therefore the estimates~\cite[(3)--(9)]{bbl}, the case analysis 3a--3c, and the resulting strict order of the intersection points on each line $M_j$ remain valid for our $\mathcal A$.
We restate the results of these computations as the following inequalities.
First, by~\cite[(6)]{bbl},
\begin{equation}\label{eq:L-M-coord}
|y(L_i\cap M_j)|<\frac{3m_{\min}}{n^9},
\qquad
|x(L_i\cap M_j)-a_i|<\frac3{n^9}.
\end{equation}
Second, by~\cite[(5)]{bbl}, the $y$-coordinates of the old intersection points off $Y_0$ are bounded away from zero:
\begin{equation}\label{eq:L-L-coord}
|y(L_i\cap L_j)|>\frac{m_{\min}}n,
\qquad \{i,j\}\neq\{n-1,n\},
\end{equation}
where the special pair $L_{n-1},L_n$ is excluded from the estimate because the $y$-coordinate of $L_{n-1}\cap L_n$ can be arbitrarily close to zero. Third, the case analysis 3a--3c and the estimates~\cite[(7)--(9)]{bbl} give the following:
\begin{equation}\label{eq:M-M-coord}
x(M_i\cap M_j)=
\begin{cases}
0, & \beta_i+\beta_j=0 \quad\text{(case 3a)},\\
-n^6\sin(2\beta_i), & \beta_i+\beta_j=\pm\pi/2 \quad\text{(case 3b)},\\
\tan(\beta_i+\beta_j)+\theta_{ij}, & \text{otherwise} \quad\text{(case 3c)},
\end{cases}
\end{equation}
where in case 3c
\begin{equation}\label{eq:theta}
\lvert\tan(\beta_i+\beta_j)\rvert<\frac n\pi,
\quad \frac3{n^8}<|\theta_{ij}|<\frac1{n^2},
\quad \operatorname{sgn}\theta_{ij}=\operatorname{sgn}\bigl((\beta_i+\beta_j)\beta_i\beta_j\bigr).
\end{equation}
Hence
\begin{equation}\label{eq:x-far}
|x(M_i\cap M_j)|\le n^6.
\end{equation}
Fourth,
\begin{equation}\label{eq:special-x}
x(M_j\cap L_{n-1})<a_{n-1}<0<a_n<x(M_j\cap L_n).
\end{equation}
Indeed,
\[
x(M_j\cap L_n)-a_n=\frac{\mu_j\,(a_n-b_j)}{m_n-\mu_j}.
\]
The denominator is positive: $|\mu_j|<2m_{\min}/n^{10}<m_{\min}\le m_n$ by~\cite[(3)]{bbl}.
Both factors $\mu_j$ and $a_n-b_j$ in the numerator have the sign of $-b_j$. For $\mu_j$ this is~\cite[(3)]{bbl} (for $\sigma=1$ the values $\mu_j$ and $b_j$ have the same sign, and here $\sigma=-1$). For $a_n-b_j$ it follows from $|b_j|\ge\tan\frac{\pi}{2n}>\frac1n>|a_n|$.
The computation for $L_{n-1}$ is similar.
Together with the signs of the slopes $m_{n-1}<0<m_n$, this implies that the new lines intersect the special pair strictly above $Y_0$.

Moreover, by~\cite[(3)]{bbl}, $|b_j|<2n/\pi$ and $|\mu_j|<2m_{\min}/n^{10}$.
The point $M_i\cap M_j$ lies on $M_j$, and its $y$-coordinate equals $\mu_j(x-b_j)$, where $x=x(M_i\cap M_j)$. For $n\ge2$ we have $|x-b_j|\le|x|+|b_j|<n^6+2n/\pi<2n^6$, and therefore
\[
|y(M_i\cap M_j)|=|\mu_j|\,|x-b_j|
<\frac{2m_{\min}}{n^{10}}\cdot2n^6=\frac{4m_{\min}}{n^4}.
\]
Since $4/n^4<1/n$ and $3/n^9<1/n$ for $n\ge2$, all intersection points involving at least one line $M_j$ lie in the strip
\begin{equation}\label{eq:strip}
|y|<\frac{m_{\min}}n.
\end{equation}
The threshold $m_{\min}/n$ is the same as in \eqref{eq:L-L-coord}. It separates the region where the old intersection points remain from the region where the new points appear.

Extending~\cite{bbl} to the non-simple case requires three more checks.

\emph{Old multiple points are preserved.}
Condition \eqref{eq:x-grid} gives $n$ distinct points $L_i\cap Y_0$, so all intersection points on $Y_0$ are simple, and every multiple point of $\mathcal A$ has the form $L_i\cap L_j$.
By \eqref{eq:L-L-coord}, it lies outside the strip \eqref{eq:strip}.
This holds for $L_{n-1}\cap L_n$ as well: this point lies below $Y_0$, so if it is multiple, it equals $L_k\cap L_{n-1}$ for some $k\notin\{n-1,n\}$, and \eqref{eq:L-L-coord} applies.
All intersection points on the new lines lie in the strip, so no new line passes through a multiple point of $\mathcal A$, and these points are preserved together with their multiplicities.

\emph{No new multiple points appear.}
A multiple point of $\mathcal B$ on a new line $M_j$ would mean that two intersection points on $M_j$ coincide.
But item~4 of the proof of~\cite[Prop.~3.1]{bbl} (with \eqref{eq:special-x} in place of the signs of the special pair) places the intersection points of $M_j$ with all other lines of $\mathcal B$ at pairwise distinct $x$-coordinates.
Also, $L_{n-1}\cap L_n$ is the only old intersection point off $Y_0$ that may lie in the strip \eqref{eq:strip}. Its $x$-coordinate lies in the interval $(a_{n-1},a_n)$, so by \eqref{eq:special-x} no new line passes through it.

\emph{Multiple points do not affect the triangle count.}
By \eqref{eq:L-L-coord}, the strip \eqref{eq:strip} contains no old intersection points except the $n$ intersection points on $Y_0$ and possibly $L_{n-1}\cap L_n$.
The new lines enter only the faces of $\mathcal A$ that have a side on $Y_0$.\footnote{In~\cite{bbl} this is stated with reference to estimates (5) and (6), which do not take the point $L_{n-1}\cap L_n$ into account. For $L\in\{L_{n-1},L_n\}$ this point lies below $Y_0$, while by \eqref{eq:special-x} $M_j$ meets $L$ above $Y_0$; for any other line $L$ through this point, the point is multiple and lies outside the strip (see above). Hence inside the strip every line of $\mathcal A$ meets $M_j$ on an edge ending on $Y_0$, so every face that $M_j$ enters has a side on $Y_0$.}
By \eqref{eq:L-M-coord}, the new lines intersect the old ones in the $3/n^{9}$-neighborhoods of the points $a_i=x(L_i\cap Y_0)$.
By (D3), the triangles with a side on $Y_0$ alternate between the two sides of $Y_0$, as they automatically do in a simple arrangement.
Therefore inside the strip the arrangement looks the same as in the simple case, and the count from the last paragraph of the same proof carries over unchanged: the new lines cut only the $n-1$ triangles with a side on $Y_0$, and each of them is replaced by $n+1$ triangles; another $n$ triangles appear at the two far ends of the strip; there are no other new triangles, and the remaining old ones are preserved.\footnote{
Strictly speaking, only $n-2$ of the triangles with a side on $Y_0$ disappear.
The triangle $Y_0L_{n-1}L_n$ bounded by the special pair is not replaced but remains a face of $\mathcal B$. It lies below $Y_0$, while by \eqref{eq:special-x} the new lines intersect $L_{n-1}$ and $L_n$ above $Y_0$ and intersect $Y_0$ outside the interval $[a_{n-1},a_n]$, since $|b_j|>1/n$.
Hence no new line crosses its boundary, and therefore no new line enters it.
However, $n$ new triangles are formed between the lines of the special pair, so this remark does not affect the triangle count.
Also, in the intuitive picture (\cref{fig:double}), where $Y_0$ ``thickens'' into the set $\{Y_0, M_1, \ldots, M_n\}$, this preserved triangle does not differ from the disappearing ones.}
Consequently,
\[
a_3(\mathcal B)
=a_3(\mathcal A)-(n-1)+(n-1)(n+1)+n
=a_3(\mathcal A)+n^2.
\]
The argument that each bounded segment of $Y_0$ is a side of at least one triangular face of $\mathcal B$ also carries over.
A segment cannot be a side of two such faces, since both its endpoints are intersection points of multiplicity 2 (see \cref{lem:double-segment}).
Suppose that the triangles on two consecutive segments of $Y_0$ lie on the same side of $Y_0$. Both have a side on the line $L$ through the common endpoint of the segments, and this side is the edge of $L$ that starts at $Y_0$ on that side. Hence the third vertex of both triangles is the other endpoint of this edge, and the three lines through the three consecutive intersection points on $Y_0$ meet there.
In $\mathcal B$ this is impossible: the points $b_j$ alternate on $Y_0$ with the points $a_1,\ldots,a_{n-2}$, and the pair $a_{n-1},a_n$ lies between $b_j=\pm\tan\frac\pi{2n}$, since $\tan\frac\pi{2n}>\frac1n$. So among any three consecutive points at least one lies on a new line, and the new lines have no multiple points.

Finally,
\[
\{a_1,\ldots,a_{n-2}\}\cup\{b_1,\ldots,b_n\}
=\left\{\tan\frac{k\pi}{2n}:1\le|k|\le n-1\right\}.
\]
If \eqref{eq:special-slopes} holds, then, after relabeling, $a_{n-1}$ and $a_n$ remain the $x$-intercepts of the special pair, and the others, numbered in increasing order, satisfy \eqref{eq:x-grid} with $2n$ in place of $n$.
Since the same two lines remain the special pair, condition (D2) and the sign $\sigma=-1$ are inherited.
The estimate \eqref{eq:m-min} follows from $|\mu_j|<2m_{\min}/n^{10}$, whose right-hand side is less than $m_{\min}$ for $n\ge2$. So the minimum $m_{\min}(\mathcal B)$ is attained at one of the new lines.
\end{proof}

\Cref{lem:doubling} describes the doubling without the line $R$, as in~\cite{bbl}.
When $R$ is restored, the new lines also change the $R$-triangles, and we now count these changes.
Since $R$ lies far to the right, it cuts off triangles from pairs of lines that meet at their last intersection point.

\begin{figure}[H]
\centering
\resizebox{\ifdim\width>\linewidth\linewidth\else\width\fi}{!}{\input{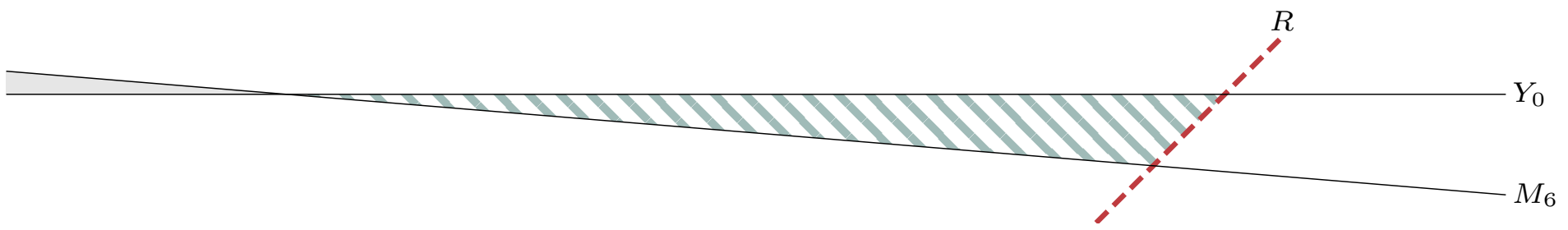}}
\caption{The $R$-triangle at the point $R\cap Y_0$ in \cref{fig:double} ($n=6$), enlarged. The line $R$ cuts it off from the right wedge $(Y_0,M_6)$. This is the triangle $\Delta_{12}$ of \textup{(W2)} at the next step, with $M_6$ in the role of $L_{10}$.}
\label{fig:double-zoom}
\end{figure}

We call a pair of lines of an arrangement a \emph{right wedge} (\cref{fig:double-zoom}) if their common point (the \emph{apex}) is
the rightmost intersection point on each of them. Beyond this point, neither line
meets any other line of the arrangement.
By \eqref{eq:x-grid}, the rightmost intersection point of the old lines $L_i$ with $Y_0$ has $x$-coordinate $a_{n-2}=\cot(\pi/n)$.
For the count of $R$-triangles, the doubling step imposes two more conditions on $\mathcal A$:
\begin{enumerate}[label=\textup{(W\arabic*)}]
\item $(Y_0,L_{n-2})$ is a right wedge of $\mathcal A$;
\item no $R$-triangle of $\mathcal A\cup\{R\}$, except the triangle $\Delta_n$ bounded by $R$, $Y_0$, and $L_{n-2}$, meets the strip \eqref{eq:strip}.
\end{enumerate}

The following lemma lists the right wedges of $\mathcal B$ that involve a new line, for $n\equiv0\pmod4$.
In our series, this case covers all iterations with $s\ge1$: $n_s=6\cdot2^s\equiv0\pmod{12}$.
The only remaining case $s=0$ (that is, $n=6$) is treated in \cref{lem:first-step}.

\begin{lemma}\label{lem:right-wedges}
Let $n\ge8$, with $n\equiv0\pmod4$, let the arrangement $\mathcal A$ satisfy \textup{(D1)--(D3)} and \textup{(W1)}, let the slope of the line $L_{n-2}$ be negative, and let $\mathcal B=\mathcal A\cup\{M_1,\ldots,M_n\}$.
Then the right wedges of $\mathcal B$ that involve at least one new line are exactly the pairs
\begin{gather}
(M_i, M_{n/2+1-i}),
\qquad 1\le i\le\frac n4,\label{eq:wedges-low}\\
(M_i,M_{3n/2-i}),
\qquad \frac n2<i<\frac{3n}{4},\label{eq:wedges-high}\\
(L_{n-2},M_{3n/4}),
\qquad
(Y_0,M_n).\label{eq:wedges-old}
\end{gather}
There are $n/2+1$ of them, and each new line belongs to exactly one.
Moreover, $\mathcal B$ satisfies \textup{(W1)} with $2n$ in place of $n$ and $M_n$ in the role of $L_{2n-2}$.
\end{lemma}

\begin{proof}
Among all pairs of new lines $(M_i,M_j)$, $i\neq j$, consider those with $\beta_i+\beta_j=\pm\pi/2$.
By \eqref{eq:new-lines}, $\beta_i+\beta_j=-\pi+(i+j-1)\pi/n$, so
\begin{gather}
\beta_i+\beta_j=-\pi/2 \iff i+j=n/2+1,\label{eq:Par1}\\
\beta_i+\beta_j=\pi/2 \iff i+j=3n/2+1.\label{eq:Par2}
\end{gather}
Each of these equalities implies $b_ib_j=1$ and $\sin(2\beta_i)=\sin(2\beta_j)$.
For such values the slopes $\mu_i$ and $\mu_j$ differ only in the term with $n^{-6}$. Hence the angle between the lines $M_i,M_j$ is small, and they intersect far from the other intersection points of the arrangement.
We call such a pair \emph{almost parallel}. It corresponds to case 3b of~\cite[Prop.~3.1]{bbl}, and by \eqref{eq:M-M-coord} $x(M_i\cap M_j)=-n^6\sin(2\beta_i)$.

For a pair of lines satisfying \eqref{eq:Par1}, both $\beta_i$ and $\beta_j$ are negative, so $\sin(2\beta_i)<0$ and $x(M_i\cap M_j)>0$. Together with the inequality $\lvert\sin(2\beta_i)\rvert>\pi/2n$ from~\cite[(3)]{bbl}, this gives $x(M_i\cap M_j)>\pi n^5/2>n^5$.
For a pair of lines satisfying \eqref{eq:Par2}, both $\beta_i$ and $\beta_j$ are positive, and a similar estimate gives $x(M_i\cap M_j)<-n^5$, so these intersection points lie far to the left and are not the rightmost points on their lines.

Let $\beta_i<0$, that is, $i\le n/2$.
Then $M_i$ has an almost parallel partner $M_j$ with $\beta_j$ determined by \eqref{eq:Par1}: $j=n/2+1-i$, and $j\neq i$, since $n/2+1$ is odd for $n\equiv0\pmod4$.
All other intersection points of $M_i$ lie in the region $|x|<n$: $x(M_i\cap Y_0)=b_i$, $x(M_i\cap L_k)$
is bounded by \eqref{eq:L-M-coord}, and the bounds on the intersection points of $M_i$ with the other new lines are given by cases 3a and 3c of~\cite[Prop.~3.1]{bbl}.
Therefore the intersection point $M_i\cap M_j$ with $x$-coordinate greater than $n^5$ is the rightmost on $M_i$.
The angle $\beta_j$ is also negative, and the same argument applies to $M_j$: its partner from \eqref{eq:Par1} is $M_i$.
Hence the pairs \eqref{eq:wedges-low} form right wedges.
With \eqref{eq:Par1}, the condition $i<j$, which orders the lines in a pair and lists each pair once, reduces to $i\le n/4$.

Now let $\beta_i>0$, that is, $i>n/2$.
As shown above, for the almost parallel line $M_j$ from \eqref{eq:Par2} we have $x(M_i\cap M_j)<-n^5$, so their intersection point lies far to the left.
We look for the rightmost of the remaining intersection points on $M_i$. In case 3c of~\cite[Prop.~3.1]{bbl}, by \eqref{eq:M-M-coord} and \eqref{eq:theta}, the intersection with a new line $M_j$ has $x=\tan(\beta_i+\beta_j)+\theta_j$ with $3/n^8<|\theta_j|<1/n^2$. In case 3a, it has exactly $x=0=\tan(\beta_i+\beta_j)$. The intersection with an old line $L_k$ lies in the $3/n^{9}$-neighborhood of the point $a_k$ by \eqref{eq:L-M-coord}. The intersection with $Y_0$ is exactly the point $b_i$.

Since $\beta_i>0$, the sums $z=\beta_i+\beta_j$ lie in the interval $(-\pi/2,\pi)$.
For $z\in(\pi/2,\pi)$ the tangent is negative, and such points lie to the left of zero and hence to the left of $a_{n-2}$.
On the interval $(-\pi/2,\pi/2)$ the tangent increases, so the rightmost point corresponds to the largest sum in this interval.
On the grid \eqref{eq:new-lines} the sums are spaced $\pi/n$ apart, and the value $z=\pi/2$ corresponds to \eqref{eq:Par2}, that is, to the almost parallel partner, whose intersection point lies far to the left.
Therefore the largest available sum is $\pi/2-\pi/n$.
It corresponds to the partner $M_{3n/2-i}$ (if $3n/2-i\neq i$) and gives
\[
x(M_i\cap M_{3n/2-i})=\tan\left(\frac\pi2-\frac\pi n\right)+\theta
=a_{n-2}+\theta,
\qquad \frac3{n^8}<|\theta|<\frac1{n^2} .
\]
Here we used $\tan(\pi/2-\pi/n)=\cot(\pi/n)=a_{n-2}$. By \eqref{eq:theta}, the sign of $\theta$ is the sign of the product $(\beta_i+\beta_j)\beta_i\beta_j$. For $i\le n-1$ all three factors are positive, since we consider $\beta_i>0$, and $\beta_j=\pi/2-\pi/n-\beta_i>0$; for $i=n$ we have $\beta_j<0$.
Among the intersections of $M_i$ with the new lines, the rightmost neighborhood is that of the point $a_{n-2}$. By \eqref{eq:L-M-coord}, the intersection of $M_i$ with the rightmost old line $L_{n-2}$ falls into the neighborhood of the same point (for the other old lines $L_k$ we have $a_k<a_{n-2}$).
Within this neighborhood, for $i\le n-1$ the intersection with $M_{3n/2-i}$ lies farther to the right, since its offset $\theta$ is greater than $3/n^{8}$, while the offset of the point $M_i\cap L_{n-2}$ is less than $3/n^{9}$ in absolute value.
It remains to take the line $Y_0$ into account. Its intersection point with $M_i$ has $x$-coordinate $b_i$ and lies in none of the neighborhoods listed, so it has to be compared with $a_{n-2}$ separately.
The values $b_i$ increase with $i$.
For $i<n$ we have
\[
a_{n-2}-b_i\ge a_{n-2}-b_{n-1}=\cot\frac\pi n-\cot\frac{3\pi}{2n}=\frac{\sin(\pi/2n)}{\sin(\pi/n)\sin(3\pi/2n)}\ge\sin\frac\pi{2n}\ge\frac1n,
\]
so the point $M_i\cap Y_0$ lies to the left of the neighborhood of $a_{n-2}$ and is not the rightmost point on $M_i$.
For $i=n$, by the identity $\cot(x/2)=\cot x+1/\sin x$,
\begin{equation}\label{eq:b-n}
b_n=\cot\frac\pi{2n}=a_{n-2}+\frac1{\sin(\pi/n)}>a_{n-2}+\frac n\pi,
\end{equation}
so on the line $M_n$ the point $M_n\cap Y_0$ is the rightmost.
It is the rightmost on $Y_0$ as well: the other intersection points on $Y_0$ have $x$-coordinates $b_j<b_n$ and $a_k\le a_{n-2}<b_n$.
Hence $(Y_0,M_n)$ is a right wedge, the second pair in \eqref{eq:wedges-old}.

The argument with the partner $M_{3n/2-i}$ has one exception, $i=3n/4$, when $3n/2-i=i$ and there is no partner with the sum $\pi/2-\pi/n$.
The next available sum $\pi/2-2\pi/n$ gives only $x\le\cot(2\pi/n)+1/n^{2}$.
This is to the left of the neighborhood of $a_{n-2}$, which has radius less than $1/n^{2}$, since
\[
a_{n-2}-\cot\frac{2\pi}n=\frac1{\sin(2\pi/n)}\ge1 .
\]
Therefore on the line $M_{3n/4}$ the rightmost intersection point is the one with the old line $L_{n-2}$.

Thus the rightmost intersection point on $M_i$ is the one with the line $M_{3n/2-i}$ for $i\notin\{3n/4,n\}$, with $L_{n-2}$ for $i=3n/4$, and with $Y_0$ for $i=n$.
Accordingly, the right wedges are the pairs \eqref{eq:wedges-high} and the two pairs \eqref{eq:wedges-old}.
For the pairs \eqref{eq:wedges-high}, the argument above suffices: both lines have indices greater than $n/2$, and the common point of the pair is the rightmost on each of them.
The pair $(Y_0,M_n)$ has also been checked on both lines.

We check that the common point of the pair $(L_{n-2},M_{3n/4})$ is the rightmost on $L_{n-2}$ as well.
By the hypothesis of the lemma,
\[
L_{n-2}:y=\ell(x-a_{n-2}),
\qquad \ell<0.
\]
The new lines are almost horizontal compared with $L_{n-2}$: $|\mu_j|<2m_{\min}/n^{10}<|\ell|$ by~\cite[(3)]{bbl}, since the minimum in $m_{\min}(\mathcal A)$ is also taken over $\ell$.
Hence $L_{n-2}$ lies above each line $M_j$ to the left of their common point and below it to the right.
Let $U=L_{n-2}\cap M_{3n/4}$. The line $L_{n-2}$ meets $M_k$ to the right of $U$ only if $M_k$ passes below $M_{3n/4}$ at $x=x(U)$.
So it suffices to show that at $x=x(U)$ the line $M_{3n/4}$ is the lowest of the new lines.

By \eqref{eq:new-lines}, the line $M_j$ passes at height
\begin{equation}\label{eq:M-height}
\mu_j(a_{n-2}-b_j)=-\frac{m_{\min}}{n^{10}}F_j,
\quad\!\!
F_j=\Bigl(\sin(2\beta_j)+\frac1{n^6\tan\beta_j}\Bigr)\left(a_{n-2}-\tan\beta_j\right),
\end{equation}
at $x=a_{n-2}$.
By \eqref{eq:L-M-coord}, $|x(U)-a_{n-2}|<3/n^9$, so between $a_{n-2}$ and $x(U)$ the height of each new line changes by less than $|\mu_j|\cdot3/n^9<6m_{\min}/n^{19}$.
Therefore $M_{3n/4}$ is the lowest at $x=x(U)$ if at $x=a_{n-2}$ it lies below every other new line by more than $12m_{\min}/n^{19}$, that is, by \eqref{eq:M-height}, if
\begin{equation}\label{eq:F-gap}
F_{3n/4}-F_k>\frac{12}{n^9}
\qquad\text{for all }k\neq\frac{3n}4.
\end{equation}

We split $F_j=G_j+H_j$ according to the two terms in parentheses:
\[
G_j=\sin(2\beta_j)\left(a_{n-2}-\tan\beta_j\right),
\qquad
H_j=\frac{a_{n-2}-\tan\beta_j}{n^{6}\tan\beta_j}.
\]
For $G_j$ we use the identity
\begin{equation}\label{eq:tan-id}
\sin(2\beta)(a_{n-2}-\tan\beta)
=\frac{\sin(2\beta+\pi/n)}{\sin(\pi/n)}-1 .
\end{equation}
By \eqref{eq:new-lines}, $2\beta_j+\pi/n=2\pi j/n-\pi$, so
\begin{equation}\label{eq:G-j}
G_j=-\frac{\sin(2\pi j/n)}{\sin(\pi/n)}-1 ,
\end{equation}
and the maximum over $j$ corresponds to the minimum of $\sin(2\pi j/n)$, that is, to $2\pi j/n=3\pi/2$:
\[
j=\frac{3n}{4},
\qquad
\beta_j=\frac\pi4-\frac\pi{2n},
\]
and the maximum equals $1/\sin(\pi/n)-1$.
For $j=3n/4\pm1$ we have $\sin(2\pi j/n)=\sin(3\pi/2\pm2\pi/n)=-\cos(2\pi/n)$, so
the maximum exceeds the values at these indices by
\[
\frac1{\sin(\pi/n)}-1-\frac{\cos(2\pi/n)}{\sin(\pi/n)}+1=
\frac{1-\cos(2\pi/n)}{\sin(\pi/n)}
=\frac{2\sin^2(\pi/n)}{\sin(\pi/n)}
=2\sin\frac\pi n\ \ge\ \frac4n .
\]
Since $|a_{n-2}-\tan\beta_j|<3n/\pi<n$~\cite[Prop.~3.1, item~2]{bbl} and
$\lvert\tan\beta_j\rvert=|b_j|>\pi/2n$ by~\cite[(3)]{bbl},
\[
|H_j|=\frac{|a_{n-2}-\tan\beta_j|}{n^6\,|\tan\beta_j|}<\frac n{n^6}\cdot\frac{2n}\pi<\frac1{n^4}.
\]
Therefore, for $k\ne3n/4$,
\[
F_{3n/4}-F_k>\frac4n-\frac2{n^4}>\frac{12}{n^9},
\]
which proves \eqref{eq:F-gap}. Hence $U$ is the rightmost point on $L_{n-2}$ among the intersection points with the new lines.

The old intersection points on $L_{n-2}$ lie to the left of $U$.
Indeed, $F_{3n/4}>0$, since both factors in \eqref{eq:M-height} are positive for $\beta_{3n/4}=\pi/4-\pi/2n$. So $M_{3n/4}$ passes below the point $(a_{n-2},0)$ of $L_{n-2}$, and $U$ lies to the right of $a_{n-2}$. By (W1), there are no old intersection points of $\mathcal A$ on $L_{n-2}$ to the right of $a_{n-2}$.
Hence $U$ is the rightmost intersection point of $\mathcal B$ on $L_{n-2}$, and the pair $(L_{n-2},M_{3n/4})$ is indeed a right wedge.

Thus the rightmost intersection point is found for each new line, and the indices in \eqref{eq:wedges-low}--\eqref{eq:wedges-old} cover $1,\ldots,n$ exactly once: \eqref{eq:wedges-low} gives $1,\ldots,n/2$, \eqref{eq:wedges-high} gives $n/2+1,\ldots,3n/4-1$ and $3n/4+1,\ldots,n-1$, and \eqref{eq:wedges-old} gives the remaining $3n/4$ and $n$.
By \cref{lem:doubling}, the new lines have no multiple points. Hence each of them has a single rightmost intersection point and belongs to at most one wedge.
Therefore \eqref{eq:wedges-low}--\eqref{eq:wedges-old} are all the right wedges that involve at least one new line, and there are exactly $n/2+1$ of them.

By \eqref{eq:b-n}, in \textup{(W1)} the role of $L_{n-2}$ with $a_{n-2}=\cot(\pi/n)$ passes to $M_n$ with $b_n=\cot(\pi/2n)$ (\cref{fig:double-zoom}).
\end{proof}

\begin{lemma}\label{lem:r-triangle}
Let $n\ge6$, with either $n\equiv0\pmod4$ or $n=6$, let the arrangement $\mathcal A$ satisfy \textup{(D1)--(D3)}, \textup{(W1)}, \textup{(W2)}, and let $\mathcal B=\mathcal A\cup\{M_1,\ldots,M_n\}$.
Suppose also that the slope of the line $L_{n-2}$ is at most $1/2$ in absolute value, negative for $n\equiv0\pmod4$ and positive for $n=6$, and that the line $R$ has the form $y=x-\delta$ with $\delta\ge2n^6$.
Then $\Delta_n$ is a face of $\mathcal A\cup\{R\}$, and in $\mathcal B\cup\{R\}$ it is replaced by exactly $n/2+1$ new $R$-triangles, while all other $R$-triangles of $\mathcal A\cup\{R\}$ remain faces.
Consequently,
\[
r(\mathcal B\cup\{R\})
=r(\mathcal A\cup\{R\})+\frac n2.
\]
Moreover, $\mathcal B$ satisfies \textup{(W2)} with $2n$ in place of $n$ and $M_n$ in the role of $L_{2n-2}$.
\end{lemma}

\begin{proof}
By \cref{lem:right-wedges} for $n\equiv0\pmod4$ and by \cref{lem:first-step} for $n=6$, the right wedges of $\mathcal B$ that involve at least one new line are the $n/2+1$ pairs \eqref{eq:wedges-low}--\eqref{eq:wedges-old} and \eqref{eq:wedges-first-step}, respectively, and each new line belongs to exactly one of them.
Other right wedges, formed only by old lines $L_i$, give no new $R$-triangles, and we do not consider them.

\emph{Step 1. The new $R$-triangles are exactly the triangles that $R$ cuts off from these wedges.} For each line $D$ that belongs to a right wedge,
$x(R\cap D)>n^6$.
Indeed, by hypothesis $R:y=x-\delta$ with $\delta\ge 2n^6$, and the slopes of all lines in the wedges are at most $1/2$ in absolute value: $Y_0$ has slope zero, the slope of $L_{n-2}$ is bounded by the hypothesis of the lemma, and the new lines have $|\mu_j|<2m_{\min}(\mathcal A)/n^{10}\le 1/n^{10}$ by~\cite[(3)]{bbl}, since $m_{\min}(\mathcal A)\le|m_{n-2}|\le1/2$.
The $x$-intercepts $a$ in the equations of these lines are less than $n$ in absolute value (they are $b_j$ or $a_{n-2}$), so
\[
x(R\cap D)=\frac{\delta-ma}{1-m}\ge\frac{\delta-n/2}{3/2}\ge\frac43n^6-\frac n3>n^6 .
\]

The apexes of the wedges in $\mathcal{B}$ have $x$-coordinates at most $n^6$. For a pair of new lines, the apex is a point $M_i\cap M_j$, which satisfies \eqref{eq:x-far}. For the two pairs with an old line, the apex has $x$-coordinate less than $n$.
For the pair with the line $L_{n-2}$, the apex is a point $L_{n-2}\cap M_j$, and by \eqref{eq:L-M-coord} and~\cite[(4)]{bbl}
\[
x(L_{n-2}\cap M_j)<a_{n-2}+\frac3{n^9}=\cot\frac\pi n+\frac3{n^9}<\frac n\pi+1 .
\]
For the pair $(Y_0,M_n)$, the apex is the point $Y_0\cap M_n$, and by~\cite[(3)]{bbl}
\[
x(Y_0\cap M_n)=b_n=\cot\frac\pi{2n}<\frac{2n}\pi .
\]
Both expressions are less than $n$ for $n\ge6$. Therefore $R$ intersects the lines of the wedge $(D,D')$ on the rays to the right of the apex, and the three lines bound a nondegenerate triangle $T$.

This triangle is a face of the arrangement $\mathcal{B}\cup \{R\}$.
Indeed, the rays forming the wedge extend to the right from its apex $Z=D\cap D'$ without crossing the other lines.
Therefore a line $D''$ of the arrangement can cross the sides of $T$ on $D$ and $D'$ only at the apex $Z$ itself.
But only two lines pass through $Z$, since each of these wedges involves at least one new line $M_j$, and by \cref{lem:doubling} the new lines have no multiple points.
Hence $D''$ crosses the boundary of $T$ at most once (possibly on the line $R$) and cannot cut off any part of $T$.

Thus $R$ closes each of the $n/2+1$ right wedges into a triangular face of $\mathcal B\cup\{R\}$.
The same argument applies to the wedge $(Y_0,L_{n-2})$ of $\mathcal A$: its apex is simple, since all intersection points on $Y_0$ are distinct by \eqref{eq:x-grid}, and the slopes $0$ and $m_{n-2}$ are at most $1/2$ in absolute value.
Hence $\Delta_n$ is also a face of $\mathcal A\cup\{R\}$, which proves the first claim of the lemma.

There are no other new $R$-triangles. Indeed, a face of $\mathcal B\cup\{R\}$ bounded by three lines of $\mathcal A\cup\{R\}$ is also a face of $\mathcal A\cup\{R\}$, since adding lines only subdivides faces. Hence every new $R$-triangle has a side on some new line $M_j$.
This side goes from the intersection point with some line $D$ to the point $M_j\cap R$, which lies to the right of all other intersection points on $M_j$. The side has no intersection points in its interior only if $M_j\cap D$ is the rightmost intersection point on $M_j$.
This point is simple by \cref{lem:doubling} and is the apex of the wedge containing $M_j$, so $(D,M_j)$ is one of these wedges.
Hence there are exactly $n/2+1$ new $R$-triangles. They replace the single face $\Delta_n$ (its side on $Y_0$ is crossed by $M_n$ at the point $b_n>a_{n-2}$), which gives an increase of $n/2$.

\emph{Step 2. The old $R$-triangles other than $\Delta_n$ remain faces.} Let $F$ be an $R$-triangle of $\mathcal A\cup\{R\}$ other than $\Delta_n$, and suppose that a new line $M_j$ enters the interior of $F$.
The face $F$ is convex and bounded, so $M_j$ crosses its boundary at exactly two points.
The side of $F$ on $R$ contains at most one of them, since $M_j$ and $R$ intersect once.
Hence at least one of these two points lies on a side belonging to some line $L$ of $\mathcal A$. This point is the intersection point $M_j\cap L$, and by \eqref{eq:strip} it lies in the strip $|y|<m_{\min}(\mathcal A)/n$.
But this point belongs to $F$, and by \textup{(W2)} $F$ does not meet this strip, a contradiction.
Thus no new line enters the interior of $F$, and $F$ remains a face of $\mathcal B\cup\{R\}$.

\emph{Step 3. The new $R$-triangles other than $\Delta_{2n}$ avoid the next strip, so \textup{(W2)} is inherited.} Let $m'=m_{\min}(\mathcal B)$, and let $S$ be the strip $|y|<m'/2n$, that is, the strip \eqref{eq:strip} for $\mathcal B$ at the next step.
For the old $R$-triangles the required condition already holds: by \eqref{eq:m-min}, $m'<m_{\min}(\mathcal A)$, so $S$ is contained in the strip $|y|<m_{\min}(\mathcal A)/n$, which they do not meet by \textup{(W2)}.

Now let $T$ be a new $R$-triangle corresponding to a wedge $(D,D')$ other than $(Y_0,M_n)$. Its vertices are $Z=D\cap D'$, $R\cap D$, and $R\cap D'$.
Neither $D$ nor $D'$ coincides with $Y_0$, so the point $D\cap Y_0=(a_D,0)$ is one of the intersection points on $D$, and $Z$ is the rightmost of them. The vertex $Z$ is simple, so $Z\neq D\cap Y_0$ and $x(Z)>a_D$.
Similarly, $x(Z)>a_{D'}$.
Therefore both factors in parentheses in the equality
\[
y(Z)=m_D\,\bigl(x(Z)-a_D\bigr)=m_{D'}\bigl(x(Z)-a_{D'}\bigr)
\]
are positive, so the slopes $m_D$ and $m_{D'}$ have the same sign.
Denote this sign by $\tau$, so that $\tau\,y(Z)>0$.
Further, $x(R\cap D)>n^6>n>|a_D|$ by the estimate of Step 1, so $\tau\,y(R\cap D)>0$, and the same holds for $D'$.
Thus all three vertices of $T$ lie on the same side of $Y_0$.

We estimate their distance from $Y_0$.
The slopes of both lines of the wedge are at least $m'$ in absolute value (both lines differ from $Y_0$) and at most $1/2$ by the estimate of Step 1, so $1-m_D\in[1/2,3/2]$ and
\[
\bigl|y(R\cap D)\bigr|=\frac{|m_D|\,(\delta-a_D)}{1-m_D}
\ \ge\ \frac{2\,m'\,(2n^6-n)}3\ >\ m'n^6 ,
\]
and the same holds for $D'$.
For the vertex $Z$, we take the line of the wedge that is new: each of the $n/2+1$ wedges except $(Y_0,M_n)$ involves a line $M_j$ with $j\le n-1$.
The point $M_j\cap L_{n-2}$ lies on $M_j$, and $Z$ is the rightmost intersection point on $M_j$, so by \eqref{eq:L-M-coord}
\[
x(Z)\ \ge\ x(M_j\cap L_{n-2})\ >\ a_{n-2}-\frac3{n^9} .
\]
The values $b_j=\tan\beta_j$ increase with $j$, so $b_j\le b_{n-1}=\cot\frac{3\pi}{2n}$ and
\[
x(Z)-b_j\ >\ \cot\frac\pi n-\cot\frac{3\pi}{2n}-\frac3{n^9}
=\frac{\sin\frac\pi{2n}}{\sin\frac\pi n\,\sin\frac{3\pi}{2n}}-\frac3{n^9}
\ >\ \frac{2n}{3\pi^2}-\frac3{n^9}\ >\ \frac n{16} ,
\]
where we used $\sin x>2x/\pi$ on $(0,\pi/2)$ for the numerator and $\sin x<x$ for the denominator.
Hence
\[
|y(Z)|=|\mu_j|\,\bigl(x(Z)-b_j\bigr)>\frac{m'n}{16} .
\]
Both bounds, $m'n^6$ and $m'n/16$, exceed $m'/2n$ for $n\ge6$, so all three vertices of $T$ lie in the half-plane $\tau y\ge m'/2n$. This half-plane is convex, so the whole triangle $T$ lies in it and does not meet the strip $S$.

The triangle $\Delta_{2n}=RY_0M_n$, by contrast, meets the strip $S$: two of its vertices, $R\cap Y_0=(\delta,0)$ and $Y_0\cap M_n=(b_n,0)$, lie on the line $Y_0$ itself.
\end{proof}

We need to count the triangular faces of the arrangement $\mathcal A_s\cup\{R\}$.
\Cref{lem:r-avoids} below shows that the line $R$ cuts no triangular face of $\mathcal A_s$.
Then the numbers of triangular faces of $\mathcal A_s$ and of $R$-triangles can simply be added.

\begin{lemma}\label{lem:r-avoids}
For each $0\le s\le t$, the line $R$ passes through no intersection point of $\mathcal A_s$ and has no common points with the interior of any triangular face of $\mathcal A_s$.
Consequently, adding $R$ preserves the multiplicities of all points, all triangular faces of $\mathcal A_s$ remain faces of $\mathcal A_s\cup\{R\}$, and
\begin{equation}\label{eq:a3-with-R}
a_3(\mathcal A_s\cup\{R\})
=a_3(\mathcal A_s)+r(\mathcal A_s\cup\{R\}).
\end{equation}
\end{lemma}

\begin{proof}
Denote by $V=L_2\cap L_3$ the intersection point of the starting arrangement on the two lines with large slopes $\mp\eta$. The point $V$ is simple in $\mathcal A_0$ by \cref{lem:seed}, and hence in all $\mathcal A_s$ by \cref{lem:doubling}. Set
\[
\Pi=\left\{(x,y):|x|\le\frac{n_t^6}{64},\ |y|\le\frac{n_t^6}{64}\right\}.
\]

\emph{All intersection points of $\mathcal A_s$ except $V$ lie in $\Pi$.}
For the intersection points of the starting arrangement, this is \cref{lem:seed-box}: all points except $V$ have both coordinates less than 14 in absolute value, and $n_t^6/64\ge6^6/64=729$.
An intersection point added at the step with $n=n_{s'}$, $s'<s$, lies on at least one new line $M_j$. Therefore by \eqref{eq:strip} its $y$-coordinate is less than $m_{\min}(\mathcal A_{s'})/n_{s'}<1$ in absolute value, and its $x$-coordinate does not exceed $n_{s'}^6\le(n_t/2)^6=n_t^6/64$ in absolute value.
Indeed, for points of the form $M_i\cap M_j$, the $x$-coordinate is bounded by \eqref{eq:x-far}. The point $M_j\cap Y_0$ has $x$-coordinate $b_j$ with $|b_j|<2n_{s'}/\pi$ by~\cite[(3)]{bbl}, and the $x$-coordinate of the point $M_j\cap L_i$ differs from $a_i$ by less than $3/n_{s'}^9$ by \eqref{eq:L-M-coord}, where $|a_i|\le\cot(\pi/n_{s'})<n_{s'}/\pi$ by \eqref{eq:x-grid}. So in the last two cases the $x$-coordinate is less than $n_{s'}$ in absolute value.

\emph{$R$ does not enter $\Pi$ and contains no intersection points of $\mathcal A_s$.}
In $\Pi$ we have
\[
x-y\le \frac{n_t^6}{32}<2n_t^6=\delta,
\]
whereas $x-y=\delta$ everywhere on $R$. For the point $V$, the difference $x-y$ equals $3n_t^6\neq\delta$. Hence $R$ contains no intersection points of $\mathcal A_s$: adding $R$ preserves the multiplicities, and all points on $R$ in $\mathcal A_s\cup\{R\}$ are simple.

\emph{Faces without the vertex $V$.} A bounded face of a line arrangement is a convex polygon, hence the convex hull of its vertices, and each of its vertices is an intersection point of lines of the arrangement.
The set $\Pi$ is also convex. Hence a face all of whose vertices lie in $\Pi$ lies entirely in $\Pi$ and does not meet $R$.

\emph{No triangular face has the vertex $V$, for any $s$.}
Here the lines $L_1,\ldots,L_6$ keep their labels from the starting arrangement.
On the line $L_2$ the slope is $-\eta<0$, so $x=-1/\sqrt3-y/\eta$, and the order of the intersection points by $x$-coordinate is opposite to their order by $y$-coordinate.
By \cref{cor:rightmost-L2}, the point $V$ is the rightmost intersection point of $\mathcal A_0$ on $L_2$, and its neighbor is $L_2\cap L_6$. Hence $y(L_2\cap L_6)$ is the smallest $y$-coordinate among the intersection points of $\mathcal A_0$ on $L_2$ other than $V$, and $y(L_2\cap L_6)<-1$ by \cref{lem:seed-box}.
Every new intersection point on $L_2$, at whatever step it appears, lies in the strip \eqref{eq:strip} of its step.
These strips shrink: at step $s'$ the strip is $|y|<m_{\min}(\mathcal A_{s'})/n_{s'}$, where the denominator grows with $s'$ and the numerator decreases by \eqref{eq:m-min}.
The widest strip corresponds to the step $s'=0$, where $m_{\min}(\mathcal A_0)=2/5$ by \cref{lem:uv} and $n_0=6$, and gives $|y|<1/15$.
Therefore the new intersection points lie above $L_2\cap L_6$, that is, to the left of it, and the neighbor of $V$ on $L_2$ is the same point $L_2\cap L_6$ in all $\mathcal A_s$.
Similarly, on $L_3$ the slope is $\eta>0$, and the order by $x$-coordinate coincides with the order by $y$-coordinate. Here $V$ is the leftmost intersection point, and its neighbor is the triple point $Q=L_1\cap L_3\cap L_5$ with $y(Q)<-1$. The new intersection points have $|y|<1/15$, so they lie above $Q$, that is, to the right of it.
A triangular face with the vertex $V$ would be bounded by $L_2$, $L_3$, and a line $D$ whose intersection points with both of them are neighbors of $V$. But then $D=L_6$ and at the same time $D\in\{L_1,L_5\}$, a contradiction.

\emph{Conclusion.} No triangular face of $\mathcal A_s$ meets $R$, so all of them remain faces of $\mathcal A_s\cup\{R\}$.
Conversely, let $F$ be a triangular face of $\mathcal A_s\cup\{R\}$ without a side on $R$.
It is contained in a unique face $G$ of $\mathcal A_s$. If $R$ cut $G$, then $F$ would be one of the pieces and would have a side on $R$. Hence $R$ does not cut $G$, and $F=G$.
Thus the triangular faces of $\mathcal A_s\cup\{R\}$ are exactly the triangular faces of $\mathcal A_s$ and the $R$-triangles, which gives \eqref{eq:a3-with-R}.
\end{proof}

\subsection{\texorpdfstring{Proof of Theorem~\ref{thm:series}}{Proof of the theorem}}

We check that \cref{lem:doubling,lem:r-triangle} apply to each $\mathcal A_s$, $0\le s<t$.
For $\mathcal A_0$, where $n=6$, we check (D1) on the lines \eqref{eq:base-lines}. The $x$-intercepts of $L_1,\ldots,L_4$ are $-\sqrt3,-1/\sqrt3,1/\sqrt3,\sqrt3$, that is, $\tan(k\pi/6)$ with $1\le|k|\le2$, and for the special pair they are $-\varepsilon<0<\varepsilon$, where $\varepsilon\le1/24<1/6$ by \eqref{eq:A-num}.
Condition (D2) holds: $L_{n-1}=L_5$ (with $x$-intercept $-\varepsilon<0$) and $L_n=L_6$ have slopes $u<0<v$ by \cref{lem:uv}.
Finally, we check (D3). By \cref{lem:seed-orders}, the lines cross $Y_0$ in the order $L_1,L_2,L_5,L_6,L_3,L_4$, and by \cref{app:seed-faces} the triangular faces of $\mathcal A_0$ with a side on $Y_0$ are
\[
Y_0L_1L_2,\quad Y_0L_2L_5,\quad Y_0L_5L_6,\quad Y_0L_3L_6,\quad Y_0L_3L_4,
\]
exactly one on each bounded segment.
They alternate between the two sides of $Y_0$: by the argument at the end of the proof of \cref{lem:doubling}, two consecutive triangles on the same side would require three consecutive lines among $L_1,L_2,L_5,L_6,L_3,L_4$ to pass through one point, while the only multiple points $P=L_1\cap L_4\cap L_6$ and $Q=L_1\cap L_3\cap L_5$ are not of this form.
The transition from $\mathcal A_s$ to $\mathcal A_{s+1}$ is provided by \cref{lem:doubling}: under condition \eqref{eq:special-slopes}, the arrangement $\mathcal B$ satisfies (D1)--(D3) with $2n$ in place of $n$. At step $s$, condition \eqref{eq:special-slopes} means $\varepsilon<1/2n_s$, which holds since $\varepsilon=1/4n_t\le1/4n_s<1/2n_s$.
Therefore the original pair $L_5,L_6$ remains the special pair in all $\mathcal A_s$, so the signs $u<0<v$ and $\sigma=-1$ carry over as well.

The other hypotheses of \cref{lem:r-triangle} also hold.
First, $\delta=2n_t^6\ge2n_s^6$, since $n_s\le n_t$.
Second, the slope of the line $L_{n_s-2}$ is at most $1/2$ in absolute value and has the required sign.
For $s=0$, this line is $L_4$ with positive slope $2/5$, as required for $n=6$.
For $s\ge1$, it is the line $M_j$ of step $s-1$ with $j=n_{s-1}$ and slope $|\mu_j|<2m_{\min}(\mathcal A_{s-1})/n_{s-1}^{10}<1/2$, since $m_{\min}(\mathcal A_{s-1})\le m_{\min}(\mathcal A_0)=2/5$ by \eqref{eq:m-min} and \cref{lem:uv}.
Its parameter $\beta_j=\pi/2-\pi/2n_{s-1}$ lies in $(0,\pi/2)$, so $\sin(2\beta_j)>0$ and $b_j>0$, and by \eqref{eq:new-lines}, $\mu_j<0$.
Third, condition (W1) for $s=0$ is \cref{cor:rightmost-L4}, and the transition from $s$ to $s+1$ is provided by \cref{lem:first-step,lem:right-wedges}.
Fourth, condition (W2): for $s=0$ the strip \eqref{eq:strip} is $|y|<m_{\min}(\mathcal A_0)/n_0=1/15$ by \cref{lem:uv}. By \cref{lem:r-triangles-below}, the four $R$-triangles other than $\Delta_6=RY_0L_4$ lie in the half-plane $y\le-1$ and therefore do not meet the strip, whereas $\Delta_6$ contains the point $(\delta,0)$ and meets it. The transition from $s$ to $s+1$ is provided by \cref{lem:r-triangle}.

The arrangement $\mathcal A_0$ has exactly two triple points $P$, $Q$ and no other multiple points (\cref{lem:seed}). By \cref{lem:doubling}, each doubling step preserves the multiple points together with their multiplicities and creates no new ones, and adding $R$ does not change multiplicities (\cref{lem:r-avoids}). Hence $\mathcal C_t$ has exactly two triple points.

For $0\le s\le t$, let $A_s=a_3(\mathcal A_s)$ and $r_s=r(\mathcal A_s\cup\{R\})$. By \cref{lem:seed},
\[
A_0=10,
\qquad
r_0=5.
\]
For $0\le s<t$, \cref{lem:doubling,lem:r-triangle} give
\[
A_{s+1}=A_s+n_s^2,
\qquad
r_{s+1}=r_s+\frac{n_s}{2}.
\]
Summing the geometric progressions
\[
\sum_{s=0}^{t-1}n_s^2=36\sum_{s=0}^{t-1}4^s=12\left(4^t-1\right),
\qquad
\sum_{s=0}^{t-1}\frac{n_s}2=3\sum_{s=0}^{t-1}2^s=3\left(2^t-1\right),
\]
and substituting $n_t^2=36\cdot4^t$, $n_t=6\cdot2^t$, we obtain
\[
A_t=10+12\left(4^t-1\right)=\frac{n_t^2}{3}-2,
\qquad
r_t=5+3\left(2^t-1\right)=\frac{n_t}{2}+2,
\]
and therefore, by \cref{lem:r-avoids},
\[
a_3(\mathcal C_t)=A_t+r_t
=\frac{n_t^2}{3}+\frac{n_t}{2}.
\]
Since $N_t=n_t+2$,
\[
a_3(\mathcal C_t)=
\frac{N_t(N_t-\frac52)+1}{3},
\]
which is \eqref{eq:a3-series}. \qed

\section{Multiple points and exhaustive search for pseudolines}\label{sec:enum}

We now turn from straight line arrangements to the general case of pseudoline arrangements.
An \emph{arrangement of $n$ pseudolines} is a collection of $n$ unbounded simple curves in the plane in which every two meet exactly once and cross there.
Every line arrangement is a pseudoline arrangement, so the maximum number $K(n)$ of triangular faces in an arrangement of $n$ lines does not exceed its pseudoline analogue $K_{\mathrm{ps}}(n)$:
\begin{equation}\label{eq:K-ps}
K(n)\le K_{\mathrm{ps}}(n)\overset{\mathrm{def}}{=}\max\{a_3(\mathcal A):\ \mathcal A\ \text{is an arrangement of } n \text{ pseudolines}\}.
\end{equation}

The only known upper bound for arrangements with multiple points is \eqref{eq:a3-bound}, and among the known even arrangements it is attained only for $n=4$.
We wrote a program that enumerates all arrangements with at least a given number $K$ of triangles, with vertices of any multiplicity.
For $n=6,8,10,12$ the enumeration finishes in reasonable time, and for $n=14$ the program shows that no arrangement has 55 or more triangles.
Together with line arrangements attaining them, this gives the exact values $K(n)=K_{\mathrm{ps}}(n)=7,15,25,38,54$, which equal the conjectured bound \eqref{eq:conj-bound}.
The program does not use the parity of $n$, and for $n=11$ it rules out 33 triangles, so $K(11)=32$.

\subsection{Representing arrangements by words}

The search works with a combinatorial representation of arrangements obtained by \emph{sweeping}.
Take a curve that meets every pseudoline exactly once and move it through the arrangement so that it passes the intersection points, or \emph{crossings}, one at a time.
The pseudolines cross the curve in some order, which is a permutation.
Passing a crossing changes this permutation, and the sequence of these changes encodes the arrangement.

A straight line sweeps any line arrangement, but for pseudolines a sweep has to be shown to exist.
By Goodman's classical theorem, every pseudoline arrangement is isomorphic to a wiring diagram, that is, admits such a sweep~\cite{goodman80}, \cite[Theorem~5.1.4]{fg17}.
However, the theorem concerns the projective plane: the isomorphism need not preserve the line at infinity, and with it the bounded faces and the number $a_3$.
We therefore use the sweeping theorem of Snoeyink and Hershberger~\cite[Theorem~3.1]{sh89}: every arrangement in the affine plane, simple or not, can be swept by a pseudoline.
Such a sweep takes place in the affine plane, so the wiring diagram it produces preserves the partition of the faces into bounded and unbounded ones.

We have already considered simple arrangements in detail in~\cite[Subsection~2.1]{pu26b}, where a sweep gives a reduced word for the longest permutation $w_0$, and each letter exchanges two adjacent pseudolines.
At a crossing of multiplicity $m$, the permutations of the allowable sequence~\cite{gp93} before and after it differ by the reversal of a block of $m$ consecutive elements instead of two.
We call the positions in these permutations \emph{slots} and number them from $0$ to $n-1$.
A pseudoline arrangement thus corresponds to a word in the letters
\begin{equation}\label{eq:letter}
(g,m):\qquad \text{reversal of the block in slots } g,\ldots,g+m-1,
\end{equation}
with one letter for each crossing.
A letter $(g,m)$ accounts for the $\binom m2$ pairs of pseudolines that cross at one point.

A \emph{word} is any sequence of letters, and a \emph{wiring diagram} is a drawing of a word. We call its curves \emph{wires}.
The letters act on the current permutation from left to right, starting with the identity permutation.
A word is \emph{partial} if no two wires cross twice in it, and \emph{complete} if every two wires cross exactly once.
Appending letters $V$ to a partial word $W$ gives an \emph{extension} $WV$ of $W$, and $V$ is a \emph{continuation} of $W$.
A complete extension is a \emph{completion} of $W$.

The verdicts of the search need only a word for every arrangement, and the Snoeyink--Hershberger sweep~\cite{sh89} provides one: if the search emits no words, no arrangement of the kind sought exists.
In the other direction, every emitted word gives a wiring diagram, which is an example of such an arrangement.

Two consecutive letters $(g,m)$ and $(g',m')$ of a partial word commute when they act on disjoint sets of wires, that is, when their blocks of slots do not overlap:
\begin{equation}\label{eq:mult-commut}
(g,m)(g',m') = (g',m')(g,m)\iff g'\geq g+m \lor g \geq g'+m'.
\end{equation}
A sweep is free to take commuting letters in either order, so we regard the drawings of words that differ by such swaps as the same wiring diagram.
The search removes this commutation freedom and emits one word from each commutation class.

\begin{lemma}[canonical form]\label{lem:canonical}
Each commutation class contains exactly one word with no pair of consecutive letters $(g,m)(g',m')$ such that $g'+m'\le g$, that is, with the block of the second letter entirely to the left of the block of the first. This word is the \emph{canonical form} of the class.
\end{lemma}

\begin{proof}
The proof follows~\cite[Lemma~4.7]{pu26b}. Swapping consecutive commuting letters \eqref{eq:mult-commut} either removes a forbidden pair or creates one. Remove forbidden pairs one at a time. Each such swap decreases by one the number of pairs of letters in which the later block lies entirely to the left of the earlier one, so the process terminates. Two forbidden pairs overlap only in a triple $abc$ where the block of $c$ is to the left of the block of $b$, which in turn is to the left of the block of $a$. Then the block of $c$ is to the left of the block of $a$, and either order of swaps leads to $cba$. By Newman's lemma~\cite{newman42}, all words of the class reach the same result.
\end{proof}

The search emits only canonical forms: after a letter $(g,m)$ it forbids the letters $(g',m')$ with $g'\le g-m'$.

One more freedom comes from the sweep rather than from the arrangement: the choice of the start, or ``direction''.
The start can be chosen so that the first letter $(g_0,m_0)$ begins in slot $g_0=0$.
This choice, the \emph{first-letter fixation}, is useful for negative verdicts, where it suffices to search over words with first letter $(0,m)$.

\begin{lemma}[first letter $g_0=0$]\label{lem:first-letter}
Every affine arrangement of $n\ge2$ pseudolines, with or without multiple points, has a word whose first letter is $(0,m)$ for some $2\le m\le n$.
The canonical form (\cref{lem:canonical}) of this word also begins with $(0,m)$.
\end{lemma}

\begin{proof}
Take any sweep of the arrangement $\mathcal A$.
Let $(g,m)$ be its first letter and $v$ the crossing it produces.
The sweep passes no crossings before $v$, so the initial curve can be pushed to infinity without touching any crossing. There it runs along $n$ consecutive ends of pseudolines, and the order of the wires on it is the order of these ends.

Any $n$ consecutive ends contain one end of each pseudoline, since the two ends of a pseudoline are separated by $n-1$ ends on each side.
Rotate the initial curve past $g$ ends.
The new curve $\gamma$ still crosses each pseudoline once, and its initial order is a cyclic shift of the old one, with the wires of $v$ in slots $0,\ldots,m-1$.
The two positions of the initial curve differ only near the ends of the $g$ transferred pseudolines, and none of these pseudolines passes through $v$.
Hence the wires of $v$ still have no other crossings between $\gamma$ and $v$, and $\gamma$ can be moved across $v$ to a curve $\gamma'$. This move is the first letter $(0,m)$.

Since $\gamma'$ crosses each pseudoline exactly once, $\mathcal A\cup\{\gamma'\}$ is a pseudoline arrangement, and by~\cite[Theorem~3.1]{sh89} it can be swept starting from $\gamma'$ towards the side of $\gamma'$ that contains all crossings of $\mathcal A$ other than $v$.
Such a sweep uses two kinds of operations: passing an empty triangle, and passing or taking the first ray~\cite[Lemma~3.1]{sh89}.
The curve $\gamma'$ already crosses all pseudolines, so it never takes a ray, and the crossings of $\mathcal A$ are passed by the first operation.
A triangle in~\cite{sh89} is not a face of the arrangement but the region bounded by the sweep curve and two wires consecutive on it, up to their common crossing.
A crossing of multiplicity $m$ has $m-1$ such regions, and the sweep passes all of them at once, with one letter~\cite[p.~356]{sh89}.
The result is a complete word with first letter $(0,m)$.

Bringing this word to canonical form removes forbidden pairs as in the proof of \cref{lem:canonical}.
The first letter stays in place: moving it would require the block of the second letter to lie entirely to the left of slot 0, which is impossible.
\end{proof}

\subsection{Counting faces}

Euler's formula gives budget constraints for the search.

\begin{lemma}\label{lem:faces}
Let $\mathcal A$ be an arrangement of $n\ge3$ pseudolines.
Then the number of its bounded faces is
\begin{equation}\label{eq:euler-F}
F_b(\mathcal A)=\binom{n-1}2-\sum_v\binom{m(v)-1}2,
\end{equation}
where the sum runs over all vertices $v$ and $m(v)$ is the multiplicity of $v$.
\end{lemma}

\begin{proof}
The $k_\ell$ vertices on a pseudoline $\ell$ divide it into $k_\ell-1$ bounded segments and two rays.
Each vertex $v$ lies on $m(v)$ pseudolines, so counting the incident pairs ``vertex, pseudoline'' in two ways gives $\sum_\ell k_\ell=\sum_v m(v)$. The number of bounded edges is therefore
\begin{equation}\label{eq:edges}
E_b=\sum_\ell(k_\ell-1)=\sum_v m(v)-n.
\end{equation}
The graph of vertices and bounded edges is connected: vertices $u\in\ell$ and $v\in\ell'$ are joined along $\ell$ and $\ell'$ through their common point.
The rays lie in the unbounded face of this graph, so its bounded faces are the bounded faces of $\mathcal A$, and by Euler's formula for a connected plane graph their number is $F_b=E_b-V+1$, where $V$ is the number of vertices.
Since every pair of pseudolines crosses exactly once, the total number of pairwise crossings is
\begin{equation}\label{eq:pair-cross}
P(\mathcal{A})=\sum_v\binom{m(v)}2=\binom n2,
\end{equation}
and $\binom m2-\binom{m-1}2=m-1$ gives
\[
\sum_v\bigl(m(v)-1\bigr)=\binom n2-\sum_v\binom{m(v)-1}2 .
\]
Substituting this into $E_b-V+1=\sum_v(m(v)-1)-n+1$ and using $\binom n2-n+1=\binom{n-1}2$, we obtain \eqref{eq:euler-F}.
\end{proof}

For a simple arrangement every vertex has $m(v)=2$, all terms of the sum vanish, and \eqref{eq:euler-F} gives $\binom{n-1}2$, the ceiling on the number of bounded faces of any arrangement. The sum shows how many faces multiple points take away: a vertex of multiplicity $m$ destroys $\binom{m-1}2$ bounded faces. Not all of the remaining bounded faces are triangles. Let $a_k(\mathcal A)$ be the number of bounded faces with $k$ sides. Bounded digons do not occur, since two pseudolines cross only once. Then
\begin{equation}\label{eq:a3-split}
a_3(\mathcal A)=\binom{n-1}2-\sum_v\binom{m(v)-1}2-(a_4+a_5+\ldots).
\end{equation}
Every multiple point and every bounded non-triangular face pushes the number of triangles below the ceiling $\binom{n-1}2$, and the target $a_3\ge K$ tolerates a shortfall of at most $\binom{n-1}2-K$ triangles.

\begin{definition}\label{def:budget}
We call $B=\binom{n-1}2-K$ the \emph{budget} of the target $K$ and the number
\begin{equation}\label{eq:spent}
S(\mathcal A)=\sum_v\binom{m(v)-1}2+a_4+a_5+\ldots
\end{equation}
the \emph{spending} of the arrangement: a vertex of multiplicity $m$ spends $\binom{m-1}2$ units of the budget, and each bounded non-triangular face spends one.
\end{definition}

\begin{corollary}\label{cor:budget}
Every arrangement satisfies $a_3(\mathcal A)=\binom{n-1}2-S(\mathcal A)$, so $a_3(\mathcal A)\ge K$ if and only if
\begin{equation}\label{eq:budget}
S(\mathcal A)\le B .
\end{equation}
\end{corollary}

The budget is shared between multiple points and non-triangular faces.
For instance, for $n=8$ and $K=16$ we have $B=5$; a vertex of multiplicity 4 takes 3 units, leaving 2 for all non-triangular faces of the arrangement.

A triple point costs only one unit of the budget, so it can pay for itself: for the price of one face and three segments, some of which may be unused, it lets neighboring triangles share a segment.\footnote{The condition $a_3\ge K$ can be written as two budget constraints, one from counting faces \eqref{eq:budget} and one from counting segments, $w\le b$ in the notation of \eqref{eq:waste-limit} and \eqref{eq:segm-balance}.
For an arrangement, that is, for a complete word, they are equivalent, but on a partial word they give different prunings, so the search uses both (see \cref{sec:segment-count}).}
Empirically, the cost of lost faces grows with the multiplicity faster than the possible gain, and the optima contain no points of higher multiplicity (see \cref{tab:triple-points}).

\subsection{Search algorithm}\label{sec:search}

The search builds a word from left to right, one letter per step, and
emits exactly the canonical forms with $a_3\ge K$. Define it as the following recursive
procedure:

\begin{verbatim}
Function Search(W), where W is a partial word:
   if S(W) > B, return;                    // a3 >= K is unreachable
   if P(W) = C(n,2), output W and return;  // W is complete
   else for each letter (g,m) admissible at W, call Search(W·(g,m)).
\end{verbatim}

\noindent
The search starts from the empty word.
Let $a[0],\ldots,a[n-1]$ be the permutation\footnote{Square brackets index slots.
In~\cite[Section~4]{pu26b} this permutation was written $a_i$;
the notation is changed to avoid confusion with the triangle count $a_3$.}
of wire labels obtained from the identity permutation by applying $W$; that is, $a[i]$
is the wire occupying slot $i$ after the crossings of $W$.
A letter $(g,m)$ reverses the block $a[g],\ldots,a[g+m-1]$.
The number of wire pairs already crossed in $W$ is
$P(W)=\sum_v\binom{m(v)}2$, as in \eqref{eq:pair-cross}.
The word $W$ is complete when $P(W)=\binom n2$, and then $a[i]=n-1-i$ is the longest permutation $w_0$.
The spending $S(W)$ is defined by \eqref{eq:spent} over the vertices of $W$
and the bounded faces that $W$ has already closed.

A letter $(g,m)$ is \emph{admissible} at $W$ if it satisfies all of the following conditions.

(i) $0\le g\le n-m$ and $a[g]<a[g+1]<\cdots<a[g+m-1]$: no two wires of the block have crossed yet, so the letter crosses each pair of them exactly once and produces all $\binom m2$ pairwise crossings at once.

(ii) $S(W)+\binom{m-1}2\le B$: the new vertex of multiplicity $m$ produced by $(g,m)$ destroys $\binom{m-1}2$ bounded faces (\cref{lem:faces}), and its cost, added to the current spending, has to fit into the budget \eqref{eq:budget}.
The same condition automatically bounds $m$, so the exhaustive search needs no separate cap on the multiplicity.

(iii) If $W$ is nonempty and its last letter is $(\widehat g,\widehat m)$, then $g\ge \widehat g-m+1$. Otherwise the block of $(g,m)$ lies entirely to the left of the block of $(\widehat g,\widehat m)$, and $W\cdot(g,m)$ contains a pair forbidden by \cref{lem:canonical}, so none of its extensions is canonical.

Condition (i) is responsible for correctness. Without it, the search would emit words that correspond to no
wiring diagram, since some pair of wires would cross twice.
Condition (iii) leaves exactly one word, the canonical form, from each commutation class (\cref{lem:canonical}).
Condition (ii) is the \emph{budget pruning} (P1). The letter adds at least $\binom{m-1}2$ to the spending, and appending letters to a partial word $W$ can only increase $S(W)$. So once $S(W)>B$, every completion of $W$ contradicts \cref{cor:budget}, and the discarded branch contains no words with $a_3\ge K$.

To count the second term of the spending, the search has to know the number of sides of each face.
A bounded face of a wiring diagram receives its vertices from letters.
We call the current face between the wires in slots $i$ and $i+1$ \emph{face $i$}.
A letter that crosses the wires in slots $i$ and $i+1$, that is, whose block contains both slots, closes face $i$ and opens a new face $i$.
Letters whose block ends at slot $i$ or begins at slot $i+1$ touch only one of these wires and give face $i$ intermediate vertices.
So the number of vertices of every face, and in particular the numbers of triangles and non-triangular faces, can be counted directly from the partial word $W$.

\subsection{Prunings from properties of wiring diagrams}

Conditions (i)--(iii) alone do not yet give a practical exhaustive search.

\begin{lemma}[face boundedness]\label{lem:bounded-face}
If $\min(a[0],\ldots,a[i])<n-1-i$, that is, the set of wire labels in slots $0,\ldots,i$ is not yet final, then face $i$ closes in every completion of the partial word $W$. Otherwise it closes in no completion.
In particular, a face $i$ already opened by some letter is bounded if and only if $\min(a[0],\ldots,a[i])<n-1-i$.
\end{lemma}

\begin{proof}
A letter whose block does not contain both slots $i,i+1$ permutes wires within slots $0,\ldots,i$ or within slots $i+1,\ldots,n-1$ and does not change the set of labels on the left. Only a letter that closes face $i$ can change this set.
In the final permutation, slots $0,\ldots,i$ hold the labels $n-1,\ldots,n-1-i$.
If the labels on the left are not yet final, their set changes in every completion, that is, face $i$ closes.
If they are final, then every label on the left exceeds every label on the right.
In particular, $a[i]>a[i+1]$, now and after every letter that does not close face $i$, so by condition (i) face $i$ never closes (cf.~\cite[Remark~B.8]{pu26b}).
\end{proof}

The set of labels $\{a[0],\ldots,a[i]\}$ changes only when face $i$ closes, so whether this face is bounded is decided at the moment it opens.
A bounded face $i$ that already has at least two intermediate vertices will have at least four vertices when it closes.
Let $\delta S(W)$ be the number of such faces.
They are not yet counted in $S(W)$, and each of them adds one to $S$ when it closes.
This gives a pruning by \emph{early non-triangle count} (P2): the search discards a branch as soon as $S(W)+\delta S(W)>B$, since by \cref{cor:budget} the branch contains no words with $a_3\ge K$.

\emph{Forced descent of a suffix-maximal wire} (P3). In~\cite[Appendix~B.1]{pu26b} we described a pruning based on the forced descent of a wire whose label exceeds the labels of all wires to its right in the permutation. The idea carries over to multiple points, but now only the last slot of the block of the next letter is forced ($g+m-1=\widehat g$).

\begin{lemma}[forced descent]\label{lem:forced-descent}
Let $(\widehat g,\widehat m)$ be the last letter of a partial word $W$ in canonical form. Suppose the wire $w=a[\widehat g]$ is not yet in its final slot $n-1-w$, and all labels to its right in the permutation ($a[\widehat g+1],a[\widehat g+2],\ldots$) are smaller than $w$. Then in every canonical completion of $W$ the next letter has the form $(\widehat g-m+1,m)$, that is, its block ends at slot $\widehat g$.
\end{lemma}

\begin{proof}
The wire $w$ is not yet in its final slot, so some later letter moves it. Let $\alpha$ be the first such letter.
The block of $\alpha$ contains slot $\widehat g$ but not $\widehat g+1$: all labels to the right of $w$ are smaller than $w$, while by condition (i) the block of $\alpha$ must be increasing. Hence this block ends at $\widehat g$. The letters between $(\widehat g,\widehat m)$ and $\alpha$ do not touch slot $\widehat g$, so each of them lies entirely to the left or entirely to the right of it. A letter to the left of $\widehat g$ cannot come right after $(\widehat g,\widehat m)$, since the two would form a pair forbidden by \cref{lem:canonical}. So if there are letters in between, the first of them lies to the right of $\widehat g$, and the last letter on the right is followed by a letter to the left of $\widehat g$ or by $\alpha$ itself, and these two letters again form a forbidden pair. Hence $\alpha$ comes right after $(\widehat g,\widehat m)$.
\end{proof}

After $\alpha$ the wire $w$ occupies the initial slot of its block, and all labels to its right are again smaller than $w$. By induction, the lemma applies to this wire until it reaches its final slot. The multiplicity of the next letter remains arbitrary, so the letters of the descent chain are not determined in advance: for each possible multiplicity $m$ the search keeps a single block start $\widehat g-m+1$.

\subsection{Counting segments and the perturbation principle}\label{sec:segment-count}

In the introduction we wrote down the segment balance \eqref{eq:simple-segm} for the simple case. We treated this balance in more detail in~\cite[Subsection~2.3]{pu26b}.
Each of the $n(n-2)$ bounded segments is a side of at most one triangle,
since otherwise either simplicity or the requirement that two pseudolines cross only once would be violated (see~\cref{lem:double-segment}).
We now consider what changes in the presence of multiple points.

Roudneff~\cite{roudneff96} described a perturbation principle for multiple points.
Multiple points fall into two groups. Points of the first group admit a deformation of the pseudolines near the point that lowers its multiplicity and does not decrease the total number of triangles. For points of the second group, every such deformation decreases this number.
Only the latter are of interest, since only they can give more triangles than simple crossings.
All such points are reduced in the sense of the following definition.

\begin{definition}\label{def:reduced}
The \emph{sectors} of a vertex $v$ of multiplicity $m$ are the $2m$ faces incident to it.
A vertex of multiplicity $3$ or $4$ is \emph{reduced} if its sectors that are not bounded triangles form one of the types of~\cite[Fig.~3.5]{roudneff96}.
For $m=3$, each of the two triples of alternate sectors contains at most one such sector.
For $m=4$, there are at most two such sectors, and if there are two, they are opposite.
An arrangement is \emph{reduced} if all its vertices of multiplicity $\ge3$ are reduced.
\end{definition}

For example, under a deformation a triple point turns into a new triangle, and the faces in three alternate sectors each gain an extra side.
The direction of the deformation determines which triple of sectors gains them.
At a reduced triple point each triple contains at least two triangular sectors, so every direction of deformation destroys at least these two triangles and creates only one, and the total number of triangles decreases.

\begin{lemma}[perturbation]\label{lem:perturbation}
If a vertex $v$ of multiplicity $m\ge3$ of an arrangement $\mathcal A$ is not reduced, then there is an arrangement $\mathcal A'$ of the same number of pseudolines with $a_3(\mathcal A')\ge a_3(\mathcal A)$ and more vertices.
\end{lemma}

\begin{proof}
We carry over Roudneff's proof in the projective plane~\cite[\S3]{roudneff96} to the affine setting.
In a small disk $D$ around $v$, replace the $m$ arcs through $v$ by an arrangement $U$ of $m$ arcs with the same endpoints, any two of which cross once.
The faces that meet $\partial D$ correspond to the sectors of $v$ and do not change outside $D$. So boundedness is preserved, and the number of sides can only grow: a sector remains a triangle exactly when its local face in $U$ is a wedge of two rays with a common vertex.
The new faces lie in $D$ and are bounded.
Therefore $a_3$ does not decrease if $U$ has at least as many bounded triangles as there are destroyed triangular sectors, and the number of vertices grows, since $U$ has more than one vertex.

For $m=3,4,5$ take $U$ from~\cite[Fig.~3.3, 3.4, 3.1]{roudneff96}. For $m=5$ the number of triangles never decreases, and for $m=3,4$ it decreases only at vertices of the types of~\cite[Fig.~3.5]{roudneff96}, that is, at reduced vertices.
Roudneff's case analysis uses only which sectors are marked as non-triangular.
We apply it with the sectors that are not bounded triangles marked, as in \cref{def:reduced}; an unbounded sector is then never destroyed as a triangle, and the affine count is never worse.

For $m\ge6$ take the arrangement $\mathcal B_m$ of F\"uredi and Pal\'asti~\cite[Example~2, Property~3]{fp84}. Its number of projective triangles is $p_3(\mathcal B_m)=m(m-3)/3+6-2\varepsilon/3$, where $\varepsilon=0$ for $m\equiv0\pmod3$ and $\varepsilon=2$ otherwise, and $p_3(\mathcal B_m)\ge2m$ for all $m\ge6$ (with equality for $m=6,7$).
Cut the projective plane along a pseudoline in general position. Each triangle of $\mathcal B_m$ that this pseudoline crosses gives a wedge, and the triangles it misses remain bounded triangles.
Hence there are at least $2m$ wedges and triangles in total, and there are no more destroyed sectors than new triangles.
\end{proof}

Multiple points change the segment balance in two ways. At a vertex of multiplicity $m$, each of the $m$ pseudolines through it meets $m-1$ others at once and so loses $m-2$ bounded segments. Moreover, a bounded segment with a multiple endpoint can be a side of two triangles (\cref{lem:double-segment}). We call such segments \emph{double} and denote their number by $d(\mathcal{A})$. The segment balance (\cref{lem:segment-balance}) then reads
\begin{equation}\label{eq:non-simple-segm}
3a_3(\mathcal{A})+u(\mathcal{A})+3t_3(\mathcal{A})+8t_4(\mathcal{A})+\ldots-d(\mathcal{A})=n(n-2),
\end{equation}
where $t_m(\mathcal{A})$ is the number of vertices of multiplicity $m$, and $d(\mathcal{A})$ compensates for double segments counted twice in $3a_3(\mathcal{A})$.

This is how multiple points get around the bound \eqref{eq:blanc-even} for simple arrangements: the lost segments may be the unused ones that are unavoidable in the simple case, and double segments let triangles share sides. By \cref{lem:perturbation}, only reduced multiple points can gain triangles in these ways.

\emph{Segment pruning} (P4). A discharging argument~\cite{cw17} applied to the segment balance \eqref{eq:non-simple-segm} for crossings of arbitrary multiplicity leads to the inequality $u(\mathcal{A})\le n(n-2)-3K+\tfrac32t_3(\mathcal{A})$ (\cref{cor:unused-bound}). The search abandons a subtree as soon as this inequality can no longer hold.

\emph{State table pruning} (P5). A hash table stores the state of the search once a subtree has been traversed without finding any word with the target number of triangles $a_3\geq K$. The search skips a node whose state is \emph{dominated} by a stored one: it has the same permutation and the same $\widehat g$, and every completion of it gives at most as many triangles as the same completion of the stored word. Such nodes occur often: for $n=10$ the table alone reduces the search tree by a factor of about 400. In the search for perfect arrangements~\cite[Section~4]{pu26b}, a similar table removed few nodes and made each node more expensive, so we did not use it there.

\emph{Doomed-face count} (P6). \Cref{app:doomed} extends the early non-triangle count (P2) to all \emph{doomed} faces, that is, faces that can no longer become bounded triangles.

The justification of these three prunings is considerably longer and subtler than that of the others, and they are harder to implement.
We therefore give their details in \cref{app:prunings} and enable them in the search by a separate run parameter, so that numerical results independent of \cref{app:prunings} remain available.
The parameter does not change the output, as we checked on the exhaustive enumeration for $n=10$.
The theorems of this section for $n<12$ rest only on runs without these prunings.

\emph{Pruning of non-reduced arrangements} (P7). For a negative verdict, \cref{lem:perturbation} allows the search to consider only arrangements with the largest number of vertices.

\begin{corollary}\label{cor:reduced}
Among arrangements of $n$ pseudolines with $a_3\ge K$, every arrangement with the largest number of vertices is reduced.
It is represented by a word in canonical form with first letter $(0,2)$.
\end{corollary}

\begin{proof}
If such an arrangement were not reduced, \cref{lem:perturbation} would give an arrangement with $a_3\ge K$ and more vertices.
By \cref{lem:first-letter}, the arrangement has a word with first letter $(0,m)$, and the canonical form of this word also begins with $(0,m)$.
The vertex produced by the first letter is the first one on each of its $m$ pseudolines, so $m$ consecutive rays meet at it, and $m+1$ unbounded sectors are incident to them.
Hence this vertex has at least $m+1$ non-triangular sectors, while a reduced vertex has at most two.
Therefore $m\ge3$ is impossible: in a reduced arrangement all vertices of multiplicity $\ge3$ are reduced.
\end{proof}

In particular, the pruning of non-reduced arrangements can be combined with the first-letter fixation, and then only $(0,2)$ remains.

With the pruning of non-reduced arrangements, the search uses only letters $(g,m)$ with $m\le4$ and discards a branch as soon as the known non-triangular sectors of some vertex of multiplicity 3 or 4 contradict \cref{def:reduced}.
A sector counts as non-triangular once its face has closed as a non-triangle or has become doomed (\cref{app:doomed}). The other sectors are either actual triangles or can still become triangles.
This rule discards no words of reduced arrangements.
Under this restriction the search no longer emits all words.
Nevertheless, the output contains all words of maximal length, so it is empty if and only if no arrangement with $a_3\ge K$ exists.
If the state table pruning is used as well, its entries must also store the word length (\cref{lem:state-table-reduced}).

\subsection{Computational results}

We have described the search algorithm and justified its correctness. Its implementation, run instructions, output words and run logs are available in the public repository~\cite{nonsimple-arrangements}. It remains to describe the results.

For convenience, we group the restrictions and prunings into four modes (\cref{tab:modes}), each extending the previous one.
Each further mode shrinks the search tree but relies on more proofs, and in modes III and IV the search no longer emits all words.
A word has one letter per vertex, so the longest words represent the arrangements with the largest number of vertices; mode IV still emits all the longest words of mode III.
Running several modes lets us cross-check the outputs and find the shortest chain of proofs behind each result.

\begin{table}[tp]
\centering
\small
\begingroup
\setlength{\tabcolsep}{4pt}
\caption{Search modes. The search emits words with $a_3\ge K$ only; the last column shows the subset guaranteed to be emitted.}
\label{tab:modes}
\renewcommand{\arraystretch}{1.2}
\begin{tabular}{l>{\raggedright\arraybackslash}p{5.3cm}>{\raggedright\arraybackslash}p{2.6cm}>{\raggedright\arraybackslash}p{2.65cm}}
\toprule
mode & restrictions and prunings & justification & guaranteed output \\
\midrule
I & canonical form, budget pruning (P1), early non-triangle count (P2), forced~descent~(P3) &\cref{lem:canonical}, \cref{cor:budget}, Lemmas~\ref{lem:bounded-face},~\ref{lem:forced-descent} & all canonical words \\
II & I + segment pruning (P4), state~table~pruning~(P5), doomed-face count (P6) &\cref{app:prunings}, \ref{lem:segment-balance}--\ref{lem:state-table} & same as I \\
III & II + first-letter fixation ($g_0=0$) & \cref{lem:first-letter} & a word of each arrangement \\
IV & III + pruning of non-reduced arrangements (P7) &Lemmas~\ref{lem:perturbation}, \ref{lem:state-table-reduced}, \cref{cor:reduced} & all the longest words of mode III \\
\bottomrule
\end{tabular}
\endgroup

\bigskip
\caption{Search runs. For $n=11$ ($K=33$) and $n=14$ the target $K$ exceeds the best known $a_3$ by one, and the output is empty; in the other rows $K$ is the maximum.
In mode IV the total number of words depends on the state table size, but the number of the longest words does not.
$^*$Parallel run on 32 vCPUs with a shared state table; its node count is nondeterministic.}
\label{tab:runs}
\bigskip
\begin{tabular}{lllrrrr}
\toprule
$n$ & $K$ & mode & words & longest words & nodes & time \\
\midrule
6 & 7 & I & 112 & 32 & 10\,222 & $<$1 s \\
 & & II & 112 & 32 & 3\,000 & $<$1 s \\
 & & III & 40 & 10 & 1\,028 & $<$1 s \\
 & & IV & 22 & 10 & 643 & $<$1 s \\
\midrule
8 & 15 & I & 32 & 32 & 8\,321\,398 & $<$1 s \\
 & & II & 32 & 32 & 203\,056 & $<$1 s \\
 & & III & 12 & 12 & 74\,465 & $<$1 s \\
 & & IV & 12 & 12 & 48\,619 & $<$1 s \\
\midrule
10 & 25 & I & 2\,848 & 1\,408 & 79\,036\,703\,322 & 41 min \\
 & & II & 2\,848 & 1\,408 & 61\,604\,535 & 12 s \\
 & & III & 1\,002 & 408 & 23\,781\,671 & 5 s \\
 & & IV & 728 & 408 & 16\,326\,051 & 4 s \\
\midrule
11 & 32 & II & 220 & 220 & 597\,776\,192 & 1.9 min \\
 & & III & 80 & 80 & 229\,384\,410 & 42 s \\
 & & IV & 80 & 80 & 157\,657\,333 & 42 s \\
\midrule
11 & 33 & I + ($g_0=0$) & 0 & 0 & 627\,129\,383\,709 & 5.8 h \\
 & & II & 0 & 0 & 204\,761\,604 & 33 s \\
 & & III & 0 & 0 & 78\,929\,213 & 12 s \\
 & & IV & 0 & 0 & 51\,840\,705 & 12 s \\
\midrule
12 & 38 & II & 912 & 384 & 38\,999\,559\,942 & 2.0 h \\
 & & III & 360 & 148 & 11\,972\,835\,436 & 38 min \\
 & & IV & 360 & 148 & 8\,931\,540\,630 & 39 min \\
\midrule
14 & 55 & IV & 0 & 0 & 6\,641\,676\,901\,608\rlap{$^*$} & 22.6 h $\times$ 32\rlap{$^*$} \\
\bottomrule
\end{tabular}
\end{table}

Modes I and II emit identical sets of words for $n=6,8,10$. The gain from the prunings of \cref{app:prunings} grows and reaches a factor of 1283 in the number of tree nodes for $n=10$.
Without mode II the case $n=12$ is out of reach: extrapolating the speed of mode I gives months.
However, as $n$ increases, the growth of the tree should make the state table less efficient, since its entries are constantly evicted.

In modes III and IV, the first-letter fixation $g_0=0$ shrinks the tree by another factor of 2.6--3.3, and the pruning of non-reduced arrangements by a factor of 1.3--1.6.
For $n\le12$ mode IV gains almost nothing in running time, because checking reduced vertices makes each node more expensive, but sampled estimates of subtree sizes for $n=14$ predict a severalfold gain.

The run for $n=14$ used 32 vCPUs of a cloud server.
The task was split by 8-letter word prefixes into 107\,732 subtrees, and their parallel traversal with a shared 32~GB state table took about 22.6 hours (700 vCPU-hours in total).
The other runs in \cref{tab:runs} were performed single-threaded on a laptop with a 4~GB state table.\footnote{The server has an AMD EPYC 9845 processor, and the laptop an AMD Ryzen AI 9 HX 370. In single-threaded runs on a partially loaded server, a server vCPU was slower than a laptop core by a factor of 1.4 in mode I, where the speed is limited by the processor, and by a factor of 1.15 in modes II--IV, where it is limited by memory access.}

The node counts of single-threaded runs are reproducible.
In the parallel run for $n=14$, the total number of nodes (about $6.6\cdot10^{12}$, \cref{tab:runs}) depends on the order in which the parallel traversal fills the shared state table, so it is nondeterministic. This does not affect the verdict (\cref{rem:shared-table}, which carries over to \cref{lem:state-table-reduced}).

\begin{theorem}\label{thm:small-n}
For even $6\le n\le14$, the maximum number of triangles in an arrangement of $n$ lines with vertices of any multiplicity equals the conjectured bound \eqref{eq:conj-bound}, and the same holds for pseudolines:
\begin{equation}\label{eq:small-n}
\begin{array}{c|ccccc}
n & 6 & 8 & 10 & 12 & 14\\ \hline
K(n)=K_{\mathrm{ps}}(n) & 7 & 15 & 25 & 38 & 54
\end{array}
\end{equation}
\end{theorem}

\begin{proof}
\emph{Upper bound.} The search (\cref{tab:runs}) rules out pseudoline arrangements with $a_3\ge K(n)+1$.
For $n=6,8,10$ this is shown by mode I, and for $n=12$ by mode II: these modes emit all canonical words with $a_3\ge K(n)$, and all emitted words have $a_3=K(n)$.
For $n=14$ this is shown by mode IV: with the target $a_3\ge55$ its output is empty, while by \cref{lem:state-table-reduced} it would be nonempty if such an arrangement existed.
Hence $K_{\mathrm{ps}}(n)$ does not exceed the values \eqref{eq:small-n}, and by \eqref{eq:K-ps} neither does $K(n)$.

\emph{Lower bound.} The arrangements of 8 and 14 lines with 15 and 54 triangles are the members of the series of \cref{thm:series} for $t=0$ and $t=1$.
Line arrangements for $n=6,10,12$ are shown in \cref{app:figures} (\cref{fig:n6-all,fig:n10-all,fig:n12-all}). Some of them were known before (for $n=10,12$ see~\cite{khodulev99}).

The cases $n=6,8,10$ use only mode I and depend neither on \cref{app:prunings} nor on \cref{lem:perturbation,lem:first-letter}. The case $n=12$ depends on \cref{app:prunings}, and the case $n=14$ also on \cref{cor:reduced}, which rests on both lemmas.
\end{proof}

The algorithm does not use the parity of $n$ in any way.
For odd $n=5,7,9$ it reproduces the optimal arrangements of our previous exhaustive enumeration~\cite[Section~7]{pu26b}.
We also checked the odd value $n=11$.
It is of interest because the known arrangements of 11 lines have at most 32 triangles, while the bound \eqref{eq:a3-bound} gives $a_3(\mathcal A)\leq33$.
For odd~$n\equiv3,5\pmod 6$, perfect arrangements (those attaining the upper bound \eqref{eq:a3-bound} without rounding) exist up to $n=93$, with $n=11$ the only known exception~\cite[Tables~1 and~4]{pu26b}.
To our knowledge, the absence of arrangements with $a_3=33$ has previously been checked only for simple arrangements~\cite{blanc11,savchuk25}.
Our exhaustive search covers vertices of arbitrary multiplicity.

\begin{theorem}\label{thm:n11}
No arrangement of 11 pseudolines with vertices of any multiplicity has $a_3\ge33$. Consequently, $K(11)=K_{\mathrm{ps}}(11)=32$.
\end{theorem}

\begin{proof}
\emph{Upper bound.} The search in mode I with the first-letter fixation $g_0=0$ and the target $a_3\ge33$ emitted no words (\cref{tab:runs}).
By \cref{lem:first-letter}, every arrangement is represented by a canonical form with $g_0=0$, and mode I emits all such words with $a_3\ge K$.
Hence there are no arrangements of 11 pseudolines with $a_3\ge33$, and by \eqref{eq:K-ps} no such line arrangements either.

\emph{Lower bound.} Arrangements of 11 lines with 32 triangles are known, for example the one found by Kabanovich~\cite{khodulev99}.

The proof relies on \cref{lem:first-letter} but not on \cref{app:prunings} or \cref{lem:perturbation}.
\end{proof}

Modes II--IV give the same empty output about three orders of magnitude faster than mode I (\cref{tab:runs}).
With the target $a_3\ge32$, mode II emits 220 words representing five arrangements, and all of them are simple, as are the optima found for $n=5,7,9$. So for odd $n\le11$ all optima are simple.

\subsection{Optimal arrangements}

Passing from words to arrangements with multiple points requires extending to the non-simple case the $O$-matrix~\cite[Definition~5.1]{pu26b}, which encodes wiring diagrams (cf.~Savchuk's table notation~\cite{savchuk25}).
The bijection between commutation classes of words and $O$-matrices~\cite[Lemma~5.2]{pu26b} and the bijection between $O$-matrices and wiring diagrams representing arrangements~\cite[Remark~5.3]{pu26b} have to be extended as well.
We adapted our deduplication and symmetry analysis programs~\cite[Section~5]{pu26b} to arrangements with multiple points and present their results as observations.
The formal extension requires no new ideas but is lengthy, and we omit it.

The same arrangement is represented by many wiring diagrams, since the direction of the sweep is arbitrary, and each diagram is encoded by a commutation class of words (see~\cite[Observation~2.1]{pu26b}).
A labeling of the pseudolines is determined by three choices: which pseudoline is $0$, an end of it, and a direction around infinity.
An affine arrangement has $4n$ such labelings.
Each of them gives a wiring diagram, and only one word of its commutation class is in canonical form.
An arrangement with symmetry group $H_E$~\cite[Definition~5.11]{pu26b} gives $4n/|H_E|$ canonical words.

To eliminate duplicate words that encode the same arrangement, we repeat the deduplication procedure of~\cite[Subsection~7.2]{pu26b}.
Each word gives an $O$-matrix, to which all Euclidean rotations and reflections and all projective rotations are applied.
The lexicographic minima of the resulting $O$-matrices are canonical representatives of the Euclidean and projective equivalence classes ($E$- and $P$-classes)~\cite[Definitions~5.7 and~5.8]{pu26b}. The result is shown in \cref{tab:symmetries}.

\begin{table}[tp]
\centering
\small
\setlength{\tabcolsep}{4pt}
\caption{Optimal arrangements and their symmetries. $\#E$ and $\#P$ are the numbers of Euclidean and projective classes, $H_E$ is the symmetry group of an affine arrangement, $G_P$ that of its projective closure~\cite[Definition~5.10]{pu26b}; $4C_1$ means four classes with group $C_1$.}
\label{tab:symmetries}
\begin{tabular}{lrrrlrll}
\toprule
$n$ & words & words, $g_0=0$ & $\#E$ & $H_E$ & $\#P$ & $G_P$ & $\sum_E 4n/|H_E|$ \\
\midrule
6 & 112 & 40 & 6 & $4C_1,\ 2C_3$ & 6 & $2C_1,\ 2C_2,\ 2C_3$ & $24\cdot4+8\cdot2$ \\
8 & 32 & 12 & 1 & $C_1$ & 1 & $C_1$ & $32\cdot1$ \\
10 & 2\,848 & 1\,002 & 72 & $71C_1,\ C_5$ & 70 & $63C_1,\ 6C_2,\ C_5$ & $40\cdot71+8\cdot1$ \\
12 & 912 & 360 & 19 & $19C_1$ & 19 & $19C_1$ & $48\cdot19$ \\
\bottomrule
\end{tabular}
\end{table}

Each $E$-class turns out to have $4n/|H_E|$ representatives\footnote{In the odd case the number of representatives was half as large, $2n/|H_E|$~\cite[Lemma~5.13]{pu26b}: for odd $n$ the reflection $g\mapsto n-2-g$ changes the parity of the generators $(g,2)$, and by construction the search of~\cite[Section~4]{pu26b} emitted only words in which odd generators are in the majority.}, that is, an $E$-class is exactly an affine arrangement with unlabeled pseudolines. The sums of the numbers of representatives over all $E$-classes are 112, 32, 2848 and 912, and, as they should, they equal the total numbers of emitted words.
Words with first letter $g_0=0$ form only part of the output, but by \cref{lem:first-letter} they represent all arrangements, and their deduplication gives the same $E$-classes.

Only three optima have affine symmetries: two with $H_E=C_3$ for $n=6$ (\cref{fig:n6-all}) and one with $H_E=C_5$ for $n=10$ (\cref{fig:n10-all}).
The projective closure can be more symmetric: for $n=6$ and $n=10$ some arrangements without affine symmetries have $G_P=C_2$.
For $n=10$, two projective classes each contain two optimal affine arrangements, so the 72 $E$-classes give 70 $P$-classes.

Apart from symmetries, the optimal arrangements differ in their multiple points.

\begin{table}[H]
\centering
\small
\setlength{\tabcolsep}{4pt}
\caption{Multiple points of the optima. All of them are triple points. The table gives the number of $E$-classes with a given number $t_3$ of triple points.}
\label{tab:triple-points}
\begin{tabular}{lcccccccccccc}
\toprule
$t_3$ & 0 & 1 & 2 & 3 & 4 & 5 & 6 & 7 & 8 & 9 & 10 & total \\
\midrule
$n=6$ & 2 & 1 & 1 & 2 &  &  &  &  &  &  &  & 6 \\
$n=8$ &  &  & 1 &  &  &  &  &  &  &  &  & 1 \\
$n=10$ & 36 & 8 & 3 & 5 & 3 & 3 & 1 & 4 & 3 & 2 & 4 & 72 \\
$n=12$ &  &  & 8 & 4 & 5 & 1 & 1 &  &  &  &  & 19 \\
\bottomrule
\end{tabular}
\end{table}

All multiple points of the optima are triple points, and their number varies widely: for $n=10$ it takes all values from 0 to 10 (\cref{fig:n10,fig:n10-all}). Simple optima exist for $n=6$ and $n=10$, where \eqref{eq:blanc-even} coincides with \eqref{eq:conj-bound}. For $n=8$ the bound \eqref{eq:blanc-even} already rules them out, since $8\equiv 2\pmod 6$, and the unique optimum, which has no nontrivial symmetries (\cref{tab:symmetries}), contains exactly two triple points.

For $n=12$ the bound \eqref{eq:blanc-even} also allows 38 triangles, but none of the 19 optima is simple. In the next section we consider ways to construct arrangements of an even number of pseudolines by adding one pseudoline. In particular, we recall that a pseudoline can be added to any perfect arrangement of $n-1$ pseudolines so that the result is a simple arrangement attaining \eqref{eq:blanc-even}. As noted above, perfect arrangements are known for all $n-1\equiv3,5\pmod6$ up to 93, except $n-1=11$. Hence, among $n\equiv0,4\pmod6$, $n=12$ is the only value in the range $n\le94$ without a simple arrangement attaining \eqref{eq:blanc-even}.

\section{Augmenting odd arrangements}\label{sec:augment}

We continue with ways to obtain arrangements of an even number of pseudolines, now by adding one pseudoline to a triangle-maximal odd arrangement.

\subsection{Augmentation by an outermost line}

\begin{figure}[H]
\centering
\tikzset{every picture/.append style={yscale=0.71}}
\begin{minipage}{0.85\linewidth}
\input{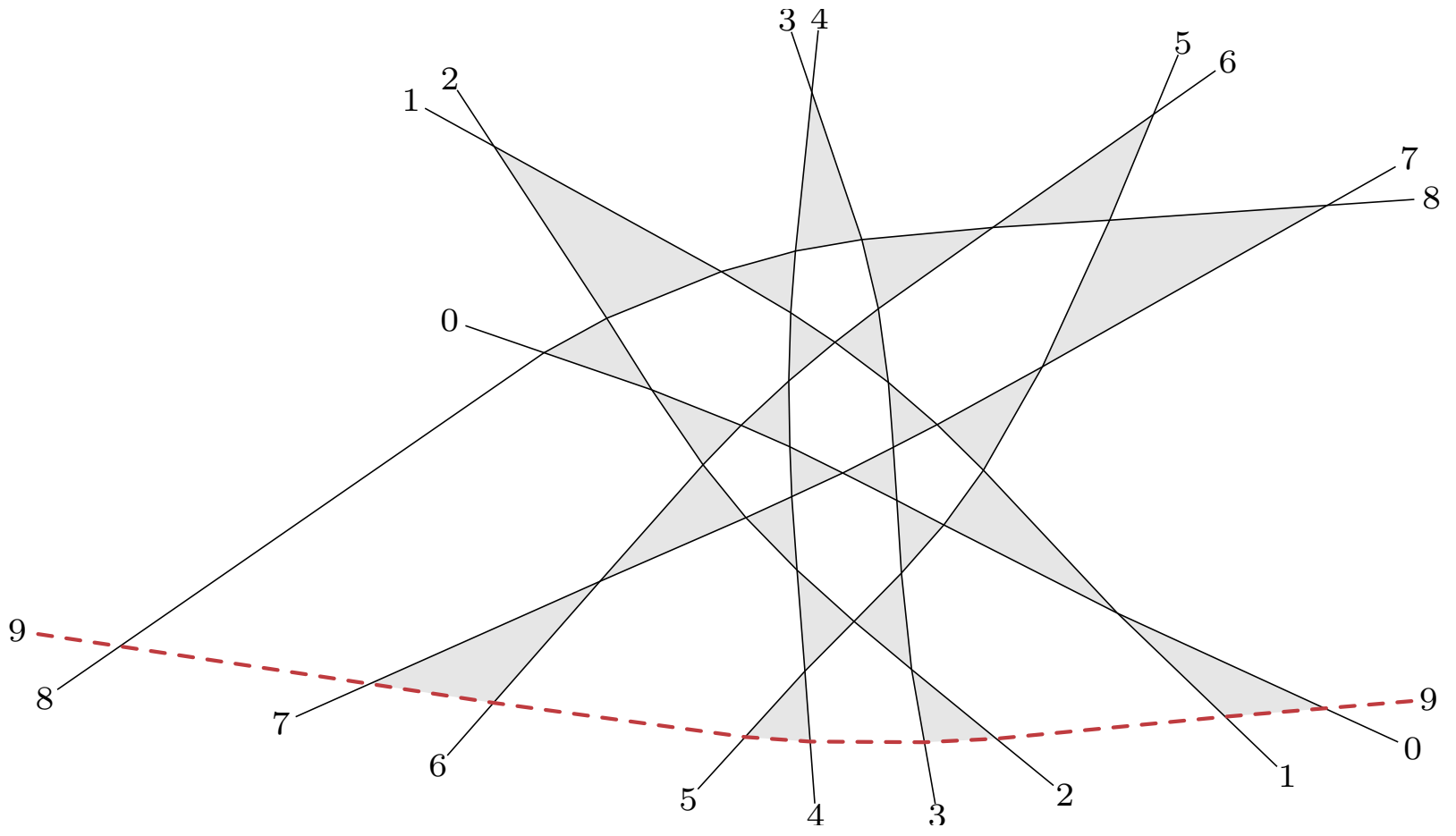}
\end{minipage}
\caption{A perfect arrangement of 9 pseudolines, labeled 0 to 8, and the outermost pseudoline 9 (dashed).
The resulting arrangement of 10 pseudolines has $a_3=25$ triangles.}
\label{fig:outermost}
\end{figure}

Even arrangements of straight lines with the largest known number of triangles are often obtained by adding one line that passes far from all intersection points~\cite{bbl,blanc11,savchuk25}. We call such a line \emph{outermost}. A perfect arrangement of $n$ lines, $n\equiv3,5\pmod 6$, has $n(n-2)/3$ triangles and $n$ external digons~\cite[Lemma~2.3]{pu26b}.
The outermost line crosses both sides of $(n-1)/2$ external digons and cuts off a new triangle from each of them (see \cref{fig:outermost} for pseudolines).
The new simple arrangement $\mathcal{A}$ of $N=n+1$ lines therefore has
\[
a_3(\mathcal{A})=\frac{n(n-2)}3+\frac{n-1}2=\frac{N(2N-5)}6
\]
triangles, so it attains the bound \eqref{eq:blanc-even}.

\subsection{Augmentation patterns for multiple points}\label{sec:augment-patterns}

For $n\equiv1\pmod 6$ there are no perfect arrangements. Here the triangle-maximal arrangements belong to the wider family of \emph{2-defective arrangements}. This family can be introduced as follows. Pass to the projective closure by adding the line at infinity. It consists of an even number $n+1$ of pseudolines, so its faces admit a checkerboard coloring, in which faces sharing a side get different colors. Call black the color of at least half of the faces. The black faces are the same for all affine charts of the projective arrangement, whichever of its pseudolines is sent to infinity. The non-triangular projective black faces are called \emph{defects}. For $n\equiv1\pmod 6$ the number of black faces is maximal precisely when there are either two quadrilateral defects or one pentagonal defect; such arrangements are called \emph{2-defective}~\cite[Lemma~2.7 and Definition~2.8]{pu26b}. An affine arrangement is triangle-maximal if and only if it is 2-defective and all its defects are external (unbounded)~\cite[Lemma~2.9]{pu26b}.

We now consider how 2-defective arrangements can be augmented.
An outermost line can still be added~\cite[Proposition~5.1]{blanc11}, but its position must be chosen so that it does not cross both unbounded sides of an external defect. Otherwise it cuts off fewer than $(n-1)/2$ new triangles. The result is again a simple arrangement, and the same bound \eqref{eq:blanc-even} holds for it.
We found two augmentation patterns for 2-defective arrangements that can exceed this bound.

\begin{figure}[p]
\centering
\tikzset{every picture/.append style={yscale=0.86}}
\captionsetup[sub]{skip=-8pt}
\begin{subfigure}[b]{0.75\linewidth}
\centering
\input{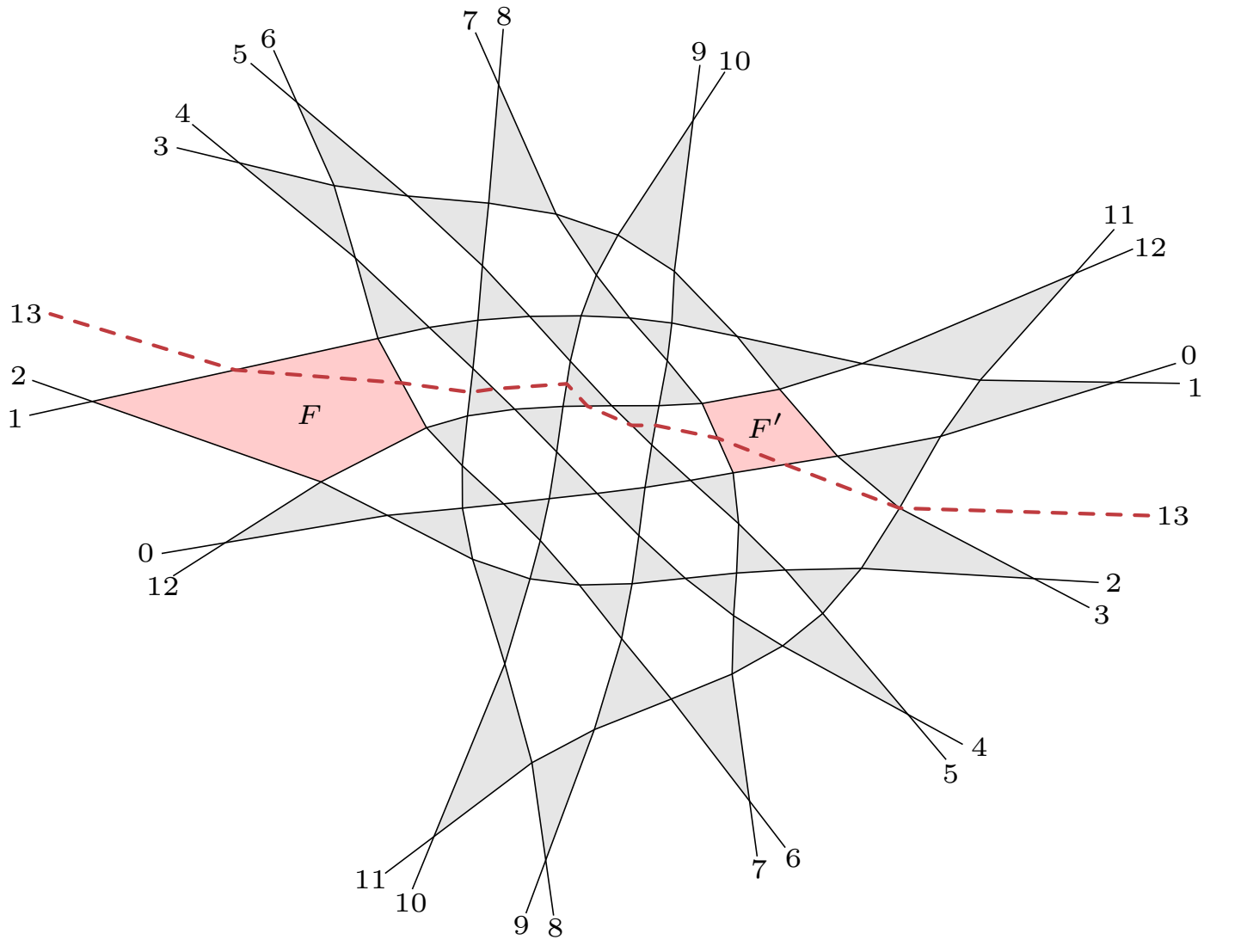}
\caption{}
\end{subfigure}

\vspace{2\bigskipamount}
\begin{subfigure}[b]{0.75\linewidth}
\centering
\input{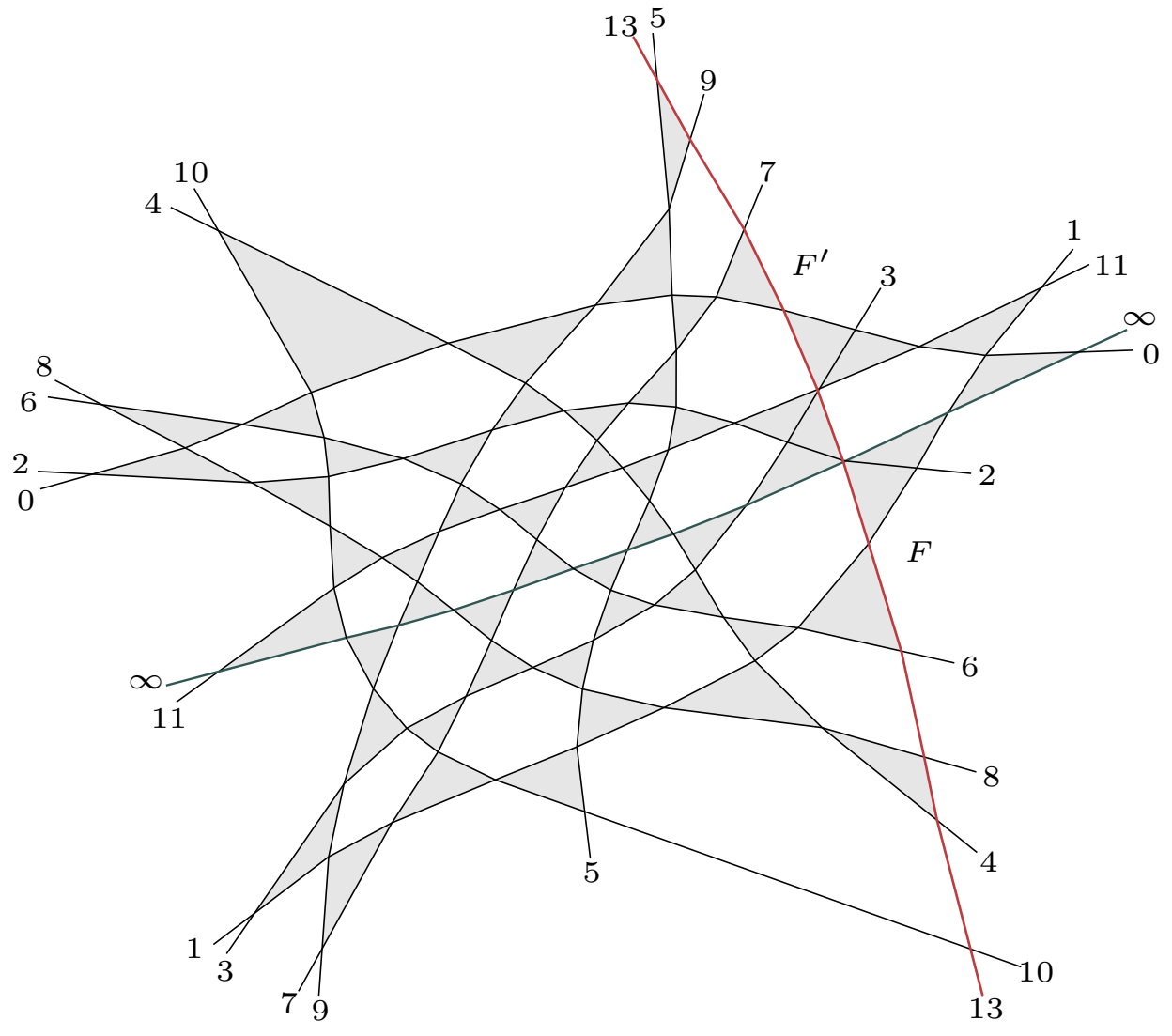}
\caption{}
\end{subfigure}
\caption{Direct augmentation. (a) A 2-defective arrangement of 13 pseudolines, labeled 0 to 12, with defects $F$, $F'$. The new pseudoline 13 is added to it (dashed). The pseudolines are not relabeled, so their labels are not in order at infinity. (b) The arrangement after a change of affine chart. Pseudoline 12 has moved to infinity. It contained sides of both defects, so both of them have become external. The former line at infinity, labeled $\infty$, has become affine. In this affine chart the arrangement has $a_3=54$ triangles.}
\label{fig:direct}
\end{figure}

In the first pattern, which we call \emph{direct augmentation}, the new pseudoline is drawn directly through the original 2-defective arrangement (\cref{fig:direct}). It must cut off triangles from both quadrilateral defects and close up through other sides and vertices. After a possible change of affine chart, the number of triangles can reach the value given by \eqref{eq:conj-bound}.

\begin{figure}[p]
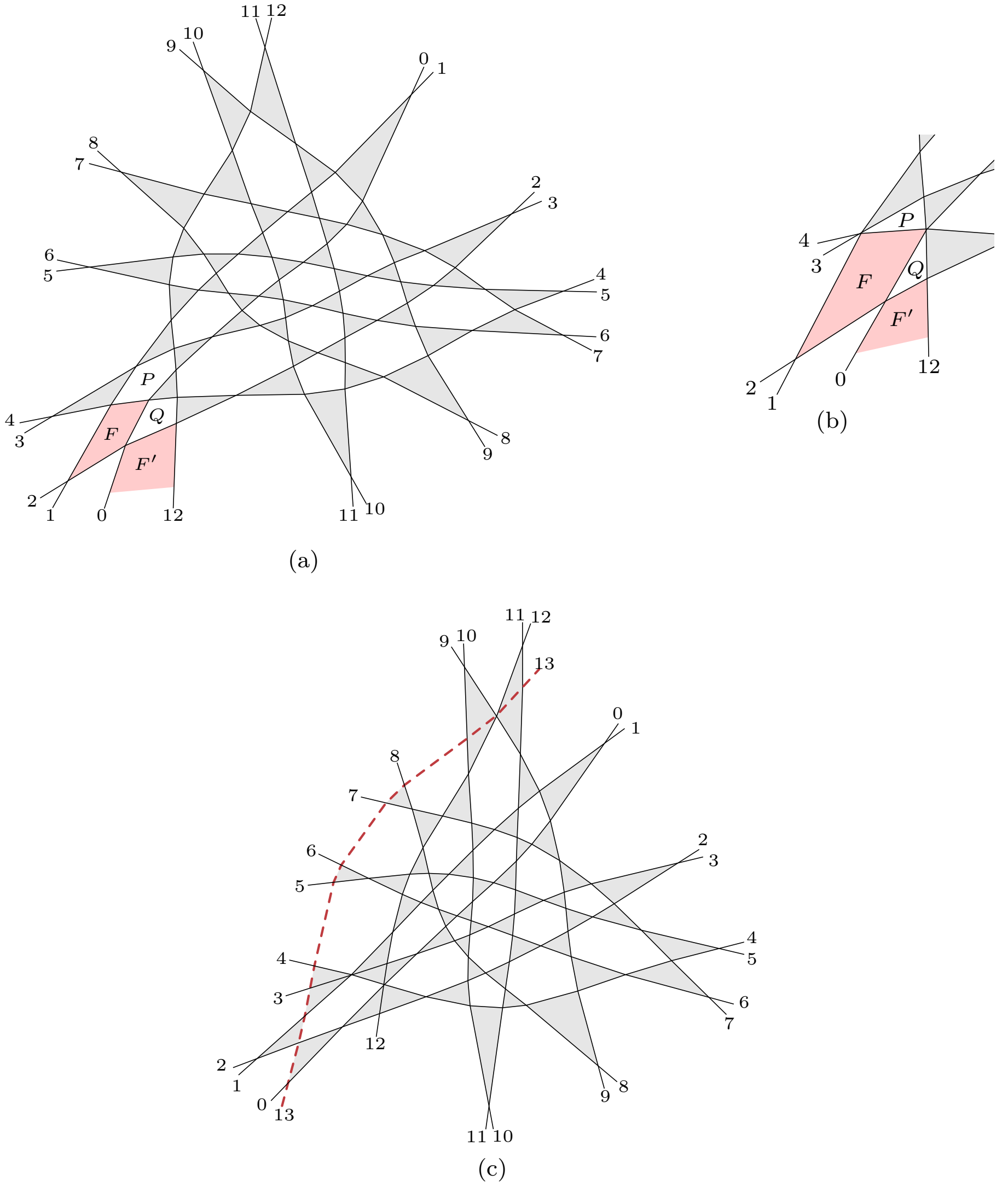

\centering
\tikzset{every picture/.append style={yscale=0.9}}
\captionsetup[sub]{skip=-8pt}
\begin{subfigure}[c]{0.64\linewidth}
\centering
\input{figures/13-before-collapse}
\caption{}
\end{subfigure}\hfill
\begin{subfigure}[c]{0.35\linewidth}
\centering
\clipbox{0pt 0pt {0.65\width} {0.6\height}}{\begin{minipage}[b]{2.5\linewidth}\input{figures/13-after-collapse}\end{minipage}}
\caption{}
\end{subfigure}

\bigskip
\begin{subfigure}[b]{0.64\linewidth}
\centering
\input{figures/13+1-after-collapse}
\caption{}
\end{subfigure}
\caption{Collapse augmentation. (a) The original arrangement. $F$ and $F'$ are the internal and the external defect. Next to them are the white quadrilateral $Q$ and the pentagon $P$. (b) The triangles with vertices at the endpoints of the common side of $P$ and $F$ are collapsed into triple points. Then $P$ and $Q$ become triangles. The common triple point of $P$ and $Q$ is already surrounded by four triangles, and the other triple point by only three. The third of these triangles is projective, bounded by pseudolines $3,4,\infty$. (c) Pseudoline 13 is added so that the second triple point gets a fourth triangle, bounded by pseudolines $1,3,13$. The resulting arrangement of 14 pseudolines has $a_3=54$ triangles and a third triple point.}
\label{fig:collapse}
\end{figure}

In the second pattern, \emph{collapse augmentation}, two particular triangles are first collapsed into triple points, and then the new pseudoline is drawn. This pattern needs an original arrangement in which the two defects $F$ and $F'$ share a vertex (\cref{fig:collapse}), and one of the two white faces at this vertex is a quadrilateral $Q$. In addition, the defect $F$ must share a side with a white pentagon $P$. The two triangles with vertices at the endpoints of the common side of $F$ and $P$ are collapsed into triple points. The pentagon $P$ and the quadrilateral $Q$ then become triangles. The total number of triangles does not change yet: two disappear and two appear. However, the triple point between $P$ and $Q$ is now surrounded by four triangles and is already reduced. The second triple point is surrounded by only three triangles. The new pseudoline must cut off a fourth triangle from the non-triangular face at the second triple point, so that this point becomes reduced as well. Otherwise, by \cref{lem:perturbation}, this triple point gives no advantage.

\subsection{Augmentation algorithm}

We wrote a program that adds one pseudoline to an arrangement in all possible ways to reach $a_3\ge K$ for a given $K$, by default \eqref{eq:conj-bound}. The source code and the output words are published in the repository~\cite{nonsimple-arrangements}. The new pseudoline follows a route through the faces of the projective closure and crosses each pseudoline once. It leaves a face through a side or through a simple vertex, which then becomes a triple point.
Working in the projective closure lets the program treat external defects as quadrilaterals, as both patterns require.

Augmentation patterns restrict the routes. For direct augmentation, the route must pass through a defect. For collapse augmentation, the program first finds suitable triangles to turn into triple points, and the route must cut off a triangle at the triple point surrounded by three triangles. If there are no quadrilateral defects, all routes are searched, including the outermost line.
The search has two stages. A dynamic program over the faces bounds from above, for every affine chart simultaneously, the number of triangles that a route can give. Then a depth-first search enumerates the routes, pruned by this bound.

After the new pseudoline has been drawn, the program tests all affine charts (that is, all choices of the pseudoline of the projective arrangement that goes to infinity) and looks for one in which $a_3\ge K$.
If the pseudoline sent to infinity passes through a triple point, the other two pseudolines of the point do not cross in the affine chart, like parallel lines. Such charts are allowed.

Because of the projective nature of the computations, it is enough to take as input one representative of each projective equivalence class ($P$-class,~\cite[Definition~5.8]{pu26b}), not all affine arrangements with an odd number of pseudolines.
For every input and every route found, the program emits one word, in an affine chart where $a_3\ge K$.

\subsection{Computational results}

In~\cite{pu26b} we obtained exhaustive enumerations of 2-defective arrangements for $n\equiv1\pmod 6$ with $n=7,13,19$, and of perfect arrangements for $n\le27$. As input we took the representatives of their $P$-classes~\cite[Subsection~7.2]{pu26b}. The results for 2-defective inputs are summarized in \cref{tab:augmentation}.

\begin{table}[tb]
\centering
\caption{Augmentation of 2-defective arrangements, by $P$-classes. The target $a_3$ is the bound \eqref{eq:conj-bound}; an input is productive if at least one of its augmentations reaches it.}
\label{tab:augmentation}
\begin{tabular}{lrrrrrrrr}
\toprule
 & & & \multicolumn{3}{c}{productive inputs} & \multicolumn{3}{c}{output $P$-classes} \\
\cmidrule(lr){4-6}\cmidrule(lr){7-9}
$n\to N$ & $a_3$ & inputs & direct & collapse & total & direct & collapse & total \\
\midrule
$7\to8$   & 15  & 1    & 1  & 1  & 1   & 1   & 1   & 1   \\
$13\to14$ & 54  & 6    & 3  & 2  & 4   & 10  & 4   & 14  \\
$19\to20$ & 117 & 1312 & 64 & 75 & 137 & 275 & 168 & 443 \\
\bottomrule
\end{tabular}
\end{table}

At $n=7$ the single input class is productive under both patterns, and they give the same output $P$-class.
At $n=13$ one of the four productive input classes is productive under both patterns, and at $n=19$ two of the 137.
At $N=14,20$ the patterns share no output class: each of the 14 output classes at $N=14$ and each of the 443 at $N=20$ comes from exactly one input class by exactly one pattern.
For many input arrangements, the program finds no route that reaches \eqref{eq:conj-bound}.\footnote{Consider direct augmentation.
To exceed \eqref{eq:blanc-even}, the new pseudoline must pass through vertices of the original arrangement, which are all simple, so it creates triple points.
At least one of them must be reduced: otherwise, by \cref{lem:perturbation}, the triple points can be removed without losing triangles, and the result is a simple arrangement.
A reduced triple point is an endpoint of double segments (\cref{sec:segment-count}). Unused and double segments form the phase seam, which contains a closed curve that is non-contractible in $\mathbb{RP}^2$ (\cref{lem:seam}; \cref{fig:phase-seam} shows it for an arrangement obtained by direct augmentation).
The local coloring differs from the checkerboard coloring of the original arrangement only in the faces cut by the new pseudoline, so the seam lies within these faces and runs along the new pseudoline, which is non-contractible as well.
We expect the criterion for reaching \eqref{eq:conj-bound} to lie in the relative position of these two curves, but we have not found it.}
For larger $n$, augmentation also reaches \eqref{eq:conj-bound}: we found examples at $n=25,31,37,43$~\cite{nonsimple-arrangements}.

As mentioned above, a 2-defective arrangement may contain one pentagonal defect. Such arrangements first appear at $n=25$. For all 2324 $P$-classes of such arrangements of 25 pseudolines, the program finds no augmentation reaching \eqref{eq:conj-bound}.

Applied to perfect arrangements, the program produces simple and non-simple arrangements, but none exceeds \eqref{eq:blanc-even}: multiple points give no advantage there.

Both augmentation patterns work for pseudolines and do not guarantee stretchability.
We stretched all 14 augmented $P$-classes for $N=14$, 353 of the 443 for $N=20$, and 25 for $N=32$; they are published in the online gallery~\cite{gallery}.
Some of them contain a pair of parallel lines; a small rotation of one line of the pair  makes them cross without decreasing the number of triangles.
Examples are shown in \cref{fig:n20,fig:n20-collapse,fig:n32}, and an arrangement of 44 pseudolines in \cref{fig:n44}.

Augmentation does not produce all known optima. A new pseudoline added to 9 pseudolines creates at most 4 triple points, and collapse adds at most 2, so the arrangements with $N=10$ and 10 triple points (\cref{fig:n10}) cannot be obtained this way. For $N=12$, there is no perfect arrangement of 11 pseudolines, and for the optimal arrangements of 11 pseudolines the program finds no augmentation reaching \eqref{eq:conj-bound}; all optima for $N=12$ are non-simple (\cref{tab:triple-points}).

\section{Discussion}\label{sec:discussion}

\subsection{Our results}

\Cref{thm:small-n} gives the maxima 7, 15, 25, 38, 54 for $n=6,8,10,12,14$, and all of them attain \eqref{eq:conj-bound}. For $n\equiv2\pmod6$, the series of \cref{sec:series} and the augmented arrangements of \cref{sec:augment} up to $n=44$ exceed \eqref{eq:blanc-even} by exactly one. For $n\equiv0,4\pmod6$, augmenting perfect arrangements attains \eqref{eq:blanc-even}. None of the arrangements we found exceeds \eqref{eq:conj-bound}.

The doubling mechanism of \cref{sec:series}, with one line temporarily removed and then restored, may also work for starting arrangements of other series, e.g.\ $2\cdot2^t+2$, $14\cdot2^t+2$, $18\cdot2^t+2$. The first two, however, contain arrangements with the same number of triangles as the perfect arrangements of the series $2\cdot2^t+1$, $14\cdot2^t+1$~\cite{bbl} with an outermost line added, so they are of little interest as new results.

Other series of line arrangements with $n\equiv2\pmod6$ attaining \eqref{eq:conj-bound}, such as $18\cdot2^t+2$, would be of interest, since they cannot be obtained by adding an outermost line. For pseudolines, such series exist.
Doubling pseudolines needs only the combinatorial part (D3) of the hypotheses of \cref{lem:doubling}.
The removed pseudoline $R$ can be drawn so that it cuts off triangles from the combinatorial analogs of the right wedges, as in \cref{lem:r-triangle}.
With the increments $n^2$ and $n/2$ the result attains \eqref{eq:conj-bound} again.
The arrangements of \cref{fig:n20-collapse,fig:n32} give the series $18\cdot2^t+2$ and $30\cdot2^t+2$ in pseudolines; their members up to $t=3$, all attaining \eqref{eq:conj-bound}, are in the repository~\cite{nonsimple-arrangements}. In contrast, the arrangement of \cref{fig:n44} violates (D3) for every choice of $Y_0$: every doubling of it gives 86 pseudolines with at least three triangles fewer than \eqref{eq:conj-bound}.

The first three members of the pseudoline series $18\cdot2^t+2$, with $n=20,38,74$, are stretchable (\cref{fig:n20-collapse,fig:n38,tab:n20-collapse-lines,tab:n38-lines,tab:n74-lines}). Stretched members do not yet give a series of lines: doubling lines also needs the metric conditions (D1) and (D2), and the removed line has to be movable arbitrarily far away, like $R$ in \eqref{eq:base-lines}.
We have a candidate construction of a series $18\cdot2^t+2$ of straight line arrangements that starts from another arrangement of 20 lines and avoids the last requirement. As $t$ grows, the point $R\cap Y_0$ approaches the point where $Y_0$ meets another line of the arrangement, rather than the point at infinity of $Y_0$. After the last doubling we change the affine chart so that $R$ becomes the line at infinity and the former line at infinity becomes a line of the arrangement. In this chart the arrangement attains \eqref{eq:conj-bound}. Its members with $n=20$ and $n=38$ are verified in rational interval arithmetic. We will present the construction separately.

\subsection{Upper bounds and future work}\label{sec:upper-bounds}

For non-simple arrangements the only applicable upper bound is \eqref{eq:a3-bound}.
It is often referred to as Tamura's bound~\cite{cb07}.
For even $n\ge6$ it is attained by no known arrangement, from the literature or from this paper.
The gap is at least $\lfloor n/6\rfloor$, and we conjecture that this holds in general.

\begin{conjecture}\label{conj:bound}
For even $n$, the number of triangular faces of an arrangement $\mathcal{A}$ of $n$ pseudolines satisfies
\begin{equation}\label{eq:conj-bound-2}
a_{3}(\mathcal{A})\le\left\lfloor\dfrac{n(n - 2)}{3}\right\rfloor-\left\lfloor\dfrac{n}{6}\right\rfloor=\begin{cases}
\dfrac{n(n - \frac52)}{3}, & n\equiv0,4\pmod 6 \\
\dfrac{n(n - \frac52) + 1}{3}. &  n\equiv2\pmod 6
\end{cases}
\end{equation}
\end{conjecture}

This bound is convenient to restate in terms of the \emph{waste}, defined as
\begin{equation}\label{eq:waste}
w(\mathcal{A})=n(n-2)-3a_3(\mathcal{A}).
\end{equation}
The waste counts the segment units that do not go into triangles and reflects the gap between the actual number of triangular faces of the arrangement and \eqref{eq:a3-bound}. The conjecture \eqref{eq:conj-bound-2} then takes the same form for all even $n$:
\begin{equation}\label{eq:conj-waste}
w(\mathcal{A})\overset{?}{\ge} n/2-1.
\end{equation}
The definition \eqref{eq:waste} gives $w\equiv n(n-2)\pmod3$. Hence $w$ can equal $n/2-1$ only for $n\equiv2\pmod 6$, and for $n\equiv0,4\pmod 6$ the inequality \eqref{eq:conj-waste} automatically strengthens to $w(\mathcal{A})\ge n/2$. So the split into cases is due to rounding alone, and the uniform bound \eqref{eq:conj-waste} allows a triangle over \eqref{eq:blanc-even} only for $n\equiv2\pmod 6$.

Why do we believe that the bound holds?
\begin{enumerate}
\item The series $6\cdot2^t+2$ of \cref{sec:series} and the series $2\cdot2^t+1$ with an outermost line added~\cite{for,bbl} attain the bound for infinitely many even $n$ in every residue class modulo 6.
\item The bound has no exceptions at small $n$ (unlike \eqref{eq:p3-bound}): it holds for all even $n$ up to 14 by \cref{thm:small-n}, and no counterexample is known for larger $n$.
\item It is natural to expect the bound to be polynomial in each residue class, as in the odd case~\cite{blanc11}. Tamura's bound \eqref{eq:a3-bound} and the known series fix the leading term as $n^2/3$. The exact values of \cref{sec:enum} for $n=6,8,10,12,14$ and the trivial case $n=4$ give two points in each residue class of $n$ modulo 6 and thereby fix the linear and the constant term.
\item The exact values for $n\le14$ alone rule out some other candidates for a sharp bound, for example the bound of Cl\'ement and Bader~\cite{cb07} with the additional term $-1$ for some residues: the gap between Tamura's bound \eqref{eq:a3-bound} and the computed optima is exactly $\lfloor n/6\rfloor$ for these $n$. This is why it was important for us to check $n=12$ computationally, where the gap first reaches 2. For this case we had to develop the prunings of \cref{app:prunings} and prove that they lose no words.
\end{enumerate}

We anticipate a natural question from the reader. The series of \cref{sec:series} starts at $n=8$, and each of its arrangements has only two triple points. What guarantees that for larger $n$ additional multiple points will not absorb the unused segments and break the conjectured bound? Perhaps this effect is absent for small $n$ only because there is not enough room in the arrangement? We have no such guarantee, and this is exactly what the conjecture asserts. But we can show what would have to happen and how the growth of $n$ affects it.

From here on we omit the argument $\mathcal{A}$ where this causes no confusion. In the simple case the waste $w$ equals the number $u$ of unused segments \eqref{eq:simple-segm}. More generally, if the only multiple points are $t_3$ triple points, \eqref{eq:non-simple-segm} gives
\begin{equation}\label{eq:waste-eq}
w=3t_3+u-d.
\end{equation}
Vertices of higher multiplicity add $\sum_{m\ge4}m(m-2)t_m=8t_4+15t_5+\cdots$ to the waste. We assume that triple points suffice to minimize the waste and omit higher multiplicities in this section.

\emph{Can the unused segments be absorbed completely?} The answer is given by the identity \eqref{eq:waste-eq}, which reads as a conservation law for charge. We assign the charge $+1$ to an unused segment, $+3$ to a triple point (it takes three segment units), and $-1$ to a double segment. A simple arrangement has no double segments, so negative charge can appear only in the presence of triple points. The total charge $3t_3+u-d$ equals $w$ and is determined by $a_3$ alone.

The question of what multiple points gain by absorbing unused segments is therefore the question whether the negative charge $d$ can exceed the positive charge $3t_3$ in absolute value. The computations give no example where it can. The inequality
\begin{equation}\label{eq:compensation}
d\overset{?}{\le}3t_3,
\end{equation}
which by \eqref{eq:waste-eq} is equivalent to $w\ge u$, holds for all arrangements from the exhaustive enumerations and from augmentation, and is strict for $t_3>0$: $d\le3t_3-1$.

The waste decreases only if the multiple points destroy more unused segments than the remainder $3t_3-d$ of their own charge. Otherwise they only shift the waste between the terms. For example, in the extreme case $n=10$, each of the two optima of \cref{fig:n10} has $t_3=10$ triple points, no unused segments at all ($u=0$), but $d=25$. All five unused segments of the simple optimum went into the remainder $3t_3-d=5$, and the total charge $w=5=n/2$ did not change. For $n\equiv2\pmod6$ the multiple points do gain, but in every arrangement we have found the gain is at most one triangle.

We could prove only a weaker form of \eqref{eq:compensation}: the discharging argument of \cref{cor:unused-bound} gives $u\le w+\tfrac32t_3$, that is, $d\le\tfrac92t_3$.

\begin{figure}[p]
\centering
\captionsetup[sub]{skip=-2pt}
\begin{subfigure}[t]{\linewidth}
\input{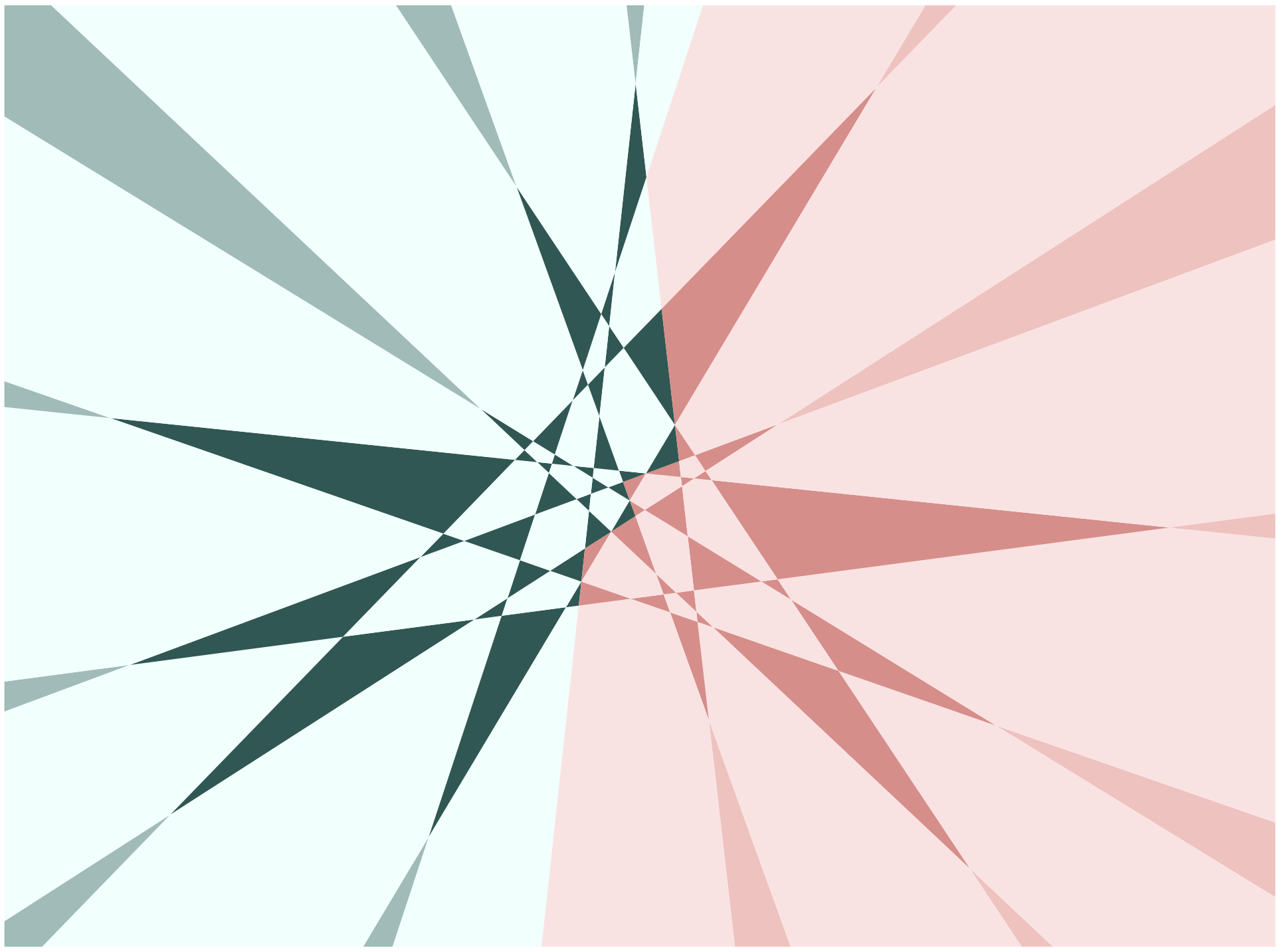}
\caption{An optimal arrangement of $14$ lines in the affine plane, from direct augmentation (\cref{sec:augment}).}
\label{fig:phase-seam-plane}
\end{subfigure}

\vspace{18pt}
\begin{subfigure}[t]{\linewidth}
\centering
\begin{minipage}{0.65\linewidth}
\input{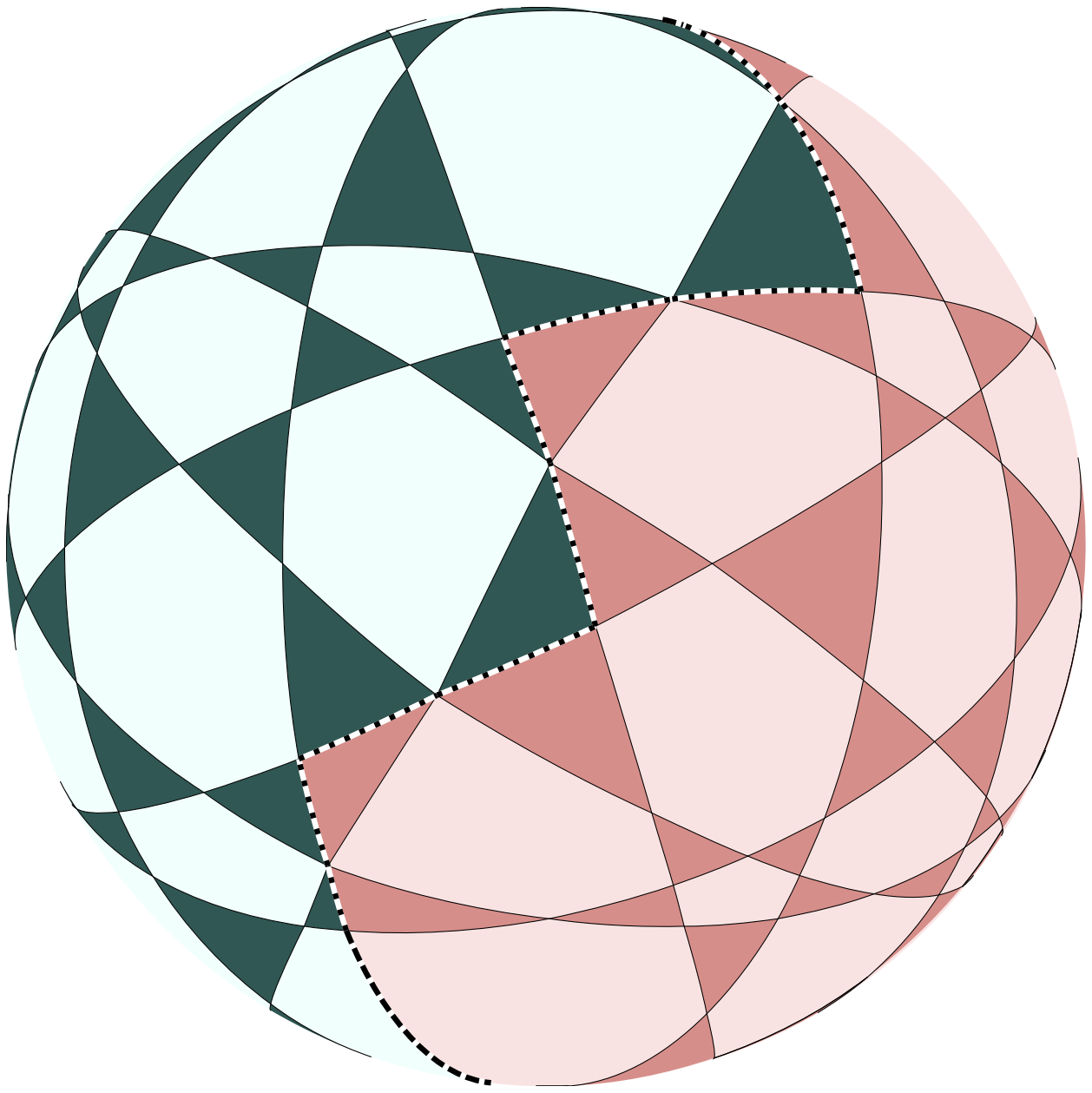}
\end{minipage}
\caption{The same arrangement with the line at infinity on the sphere (phase seam shown dashed).}
\label{fig:phase-seam-sphere}
\end{subfigure}
\caption{Phases of the local coloring and the phase seam.}
\label{fig:phase-seam}
\end{figure}

\emph{Can the growth of $n$ help the waste disappear?} In the even case the waste is not a local phenomenon confined to small neighborhoods: it forms a global structure. To see this, pass from the arrangement $\mathcal{A}$ of $n$ pseudolines to its projective closure $\mathcal{A}^+$ of $N=n+1$ pseudolines, which contains the line at infinity. The notions of unused and double segments extend to $\mathcal{A}^+$, and we denote their numbers by $u^+$ and $d^+$. Unlike in \cref{sec:augment-patterns}, here $N$ is odd, so there is no checkerboard coloring, and we use the local coloring instead: triangles are black and all other faces are white. Two faces sharing a segment then have different colors if it is an ordinary used segment, and the same color if it is an unused or double segment. So if there are few unused and double segments, the arrangement has large regions (phases) where the local coloring is a checkerboard coloring. The phases are separated by a phase seam made of unused and double segments (\cref{fig:phase-seam}).

\begin{definition}\label{def:seam}
The \emph{phase seam} $Z$ is the set of unused and double segments of the projective closure $\mathcal A^+$.
\end{definition}

\begin{lemma}\label{lem:seam}
Let $n$ be even. Then every vertex of $\mathcal A^+$ has even degree in $Z$, and exactly one connected component of $Z$ contains a closed curve that is non-contractible in $\mathbb{RP}^2$.
\end{lemma}

\begin{proof}
A vertex is an endpoint of an even number of segments. On a full turn around it, the color changes at every segment outside $Z$ and returns to the starting one, so the number of segments outside $Z$ is even, and so is the number in $Z$.

In the sphere model, take a path from a face to its antipode that avoids the vertices. It crosses each of the $n+1$ projective pseudolines an odd number of times, and since $n+1$ is odd, the total is odd. The antipode is the same projective face, with the same color, so the color changes an even number of times along the path, and the path crosses $Z$ an odd number of times.
In $\mathbb{RP}^2$ the path becomes a non-contractible closed curve $\gamma$. Since all degrees in $Z$ are even, the segments of $Z$ split into closed curves. A contractible closed curve meets $\gamma$ an even number of times, so at least one of these curves is non-contractible.
Only one component of $Z$ contains such curves: curves in different components are disjoint, while any two non-contractible closed curves in $\mathbb{RP}^2$ meet, like two projective lines.
\end{proof}

The proof uses the parity of $n$. For odd $n$ the conclusion fails: either there is no seam at all (perfect case) or the seam consists of the boundaries of the defects (see~\cite[Subsection~2.5]{pu26b}), and then it does not go around the projective plane.

Thus the charged elements are not scattered arbitrarily over the arrangement: they lie on closed contours, and one of these contours goes around the projective plane.
Even if no unused segments are left in the affine arrangement, as in \cref{fig:n10}, the seam does not disappear: in $\mathcal A^+$ it runs along the double segments of $\mathcal A$ and along the unused and double segments of $\mathcal A^+$ that lie on the former rays and on the line at infinity.
The seam is a topological obstruction. It does not depend on $n$ being small or on the lack of room in the arrangement. Multiple points change the composition of the seam (see \eqref{eq:waste-eq}) but cannot destroy it.

Unfortunately, $u$ and $d$ enter the ``length'' of the seam with the same sign, $|Z|=u^++d^+$, and \eqref{eq:waste-eq} with opposite signs. Therefore a lower bound on the length of the seam does not yield a bound on $w$. For the same reason, considerations of parity or of the checkerboard coloring do not help either: they do not distinguish $+d$ from $-d$. Moreover, the contribution of a multiple point is not determined by its neighborhood alone. It also depends on the neighborhoods of the opposite endpoints of the double segments. To prove the bound, one therefore has to pass from the statement ``the seam is unavoidable'' to the statement \eqref{eq:conj-waste} ``the waste is at least $n/2-1$''.

Blanc's bound \eqref{eq:blanc-even} does not carry over to non-simple arrangements directly. His argument is a count over pseudolines in two stages, demand and sharing. Along a pseudoline with an odd number of vertices, the alternation of triangles cannot close up, so each pseudoline requires an adjacent unused segment~\cite[Proposition~2.3]{blanc11}. This is the demand. One segment serves at most two pseudolines, and this is the sharing. Multiple points break it in two ways. A pseudoline through triple points escapes the argument. At a triple point the alternation of triangles is thrown out of step, and an odd number of such points also changes the parity of the number of vertices, so the alternation is allowed to close up. A pseudoline free of triple points still has an odd number of vertices and still makes its demand, but the demand can be met from another source. Besides an adjacent unused bounded segment, as in Blanc's argument, it can be met by a double segment with one endpoint on the pseudoline and a triple point at the other (a ``simple leg'' of that triple point, see \cref{lem:simple-legs}). The demand is then met by the three units that the triple point contributes to \eqref{eq:waste-eq}. The arrangements of \cref{fig:n10} show the first mechanism in the limit: in each of them all ten lines pass through triple points, the demand disappears completely, and $u=0$.

However, meeting the demand through $u$ or $d$ does not yield a bound on $w$. The reason is the same as for the length of the seam: double segments enter the demand with a plus sign and \eqref{eq:waste-eq} with a minus sign.

We arrive at the following conclusions.
In the even case, non-simple arrangements exceed the simple bound \eqref{eq:blanc-even} for arbitrarily large $n$ (\cref{sec:series}).
The phase seam, which contains the support of the waste, is unavoidable for topological reasons.
Every pseudoline free of triple points makes a demand for waste, and the demand is met by an unused segment or by the units of a triple point transferred through a double segment.
What is missing for a proof is a quantitative statement of how much waste the triple points must transfer to the pseudolines with a demand and leave on the pseudolines through them.

\appendix
\numberwithin{table}{section}
\numberwithin{figure}{section}
\section{Intersection orders in the starting arrangement}\label{app:seed}

We determine the coordinates and the order of the intersection points of the starting arrangement, the bounding box of these points, and the number of faces. We also treat the first doubling step ($n=6$) separately, because it partly differs from the general step.

\subsection{Parameters}

The lines $R,Y_0,L_1,\ldots,L_6$ are given by the equations \eqref{eq:base-lines}, and the parameters
$\varepsilon$, $\eta$, $\delta$ by \eqref{eq:params}. Since $n_t\ge6$, we have
\begin{equation}\label{eq:A-num}
\varepsilon\le\frac1{24},
\qquad
\eta\ge3\sqrt3\cdot6^6>2\cdot10^5,
\qquad
\delta\ge2\cdot6^6>9\cdot10^4.
\end{equation}

\begin{lemma}\label{lem:uv}
The slopes $u$ and $v$ \eqref{eq:star-slopes} of the lines $L_5$ and $L_6$ satisfy
\begin{equation}\label{eq:uv}
u<0<v,
\qquad
\frac95<|u|<2,
\qquad
3<v<4 .
\end{equation}
Thus the slopes of $L_1,L_4,L_5,L_6$ do not exceed 4 in absolute value, and the
minimum absolute value of a slope among $L_1,\ldots,L_6$ is $m_{\min}(\mathcal A_0)=2/5$.
\end{lemma}

\begin{proof}
We estimate the term $\varepsilon\sqrt3\,(2\eta+1)$ in the denominator of $u$:
\[
\varepsilon\sqrt3\,(2\eta+1)>2\varepsilon\sqrt3\,\eta
=\frac{\sqrt3}{4n_t}\cdot6\sqrt3\,n_t^6
=\frac92\,n_t^5\ \ge\ \frac92\cdot6^5\ >\ 3 ,
\]
so the denominator in \eqref{eq:star-slopes} is positive and greater than $2\eta$. Hence $u<0$ and
\[
|u|=\frac{4\eta}{2\eta-3+\varepsilon\sqrt3(2\eta+1)}<\frac{4\eta}{2\eta}=2 .
\]
On the other hand, the same denominator equals $2\eta(1+\varepsilon\sqrt3)-(3-\varepsilon\sqrt3)$,
and $3-\varepsilon\sqrt3>0$ by \eqref{eq:A-num}. So the denominator is less than $2\eta(1+\varepsilon\sqrt3)$, which gives
\[
|u|>\frac{4\eta}{2\eta(1+\varepsilon\sqrt3)}
=\frac2{1+\varepsilon\sqrt3}
\ \ge\ \frac2{1+\sqrt3/24}
\ >\ \frac95 .
\]
For $v$:
\[
v=\frac{4\sqrt3}{\sqrt3+9\varepsilon}<\frac{4\sqrt3}{\sqrt3}=4,
\qquad
v\ \ge\ \frac{4\sqrt3}{\sqrt3+9/24}
=\frac{4\sqrt3}{\sqrt3+3/8}
\ >\ 3 .
\]
The absolute values of the slopes of $L_1,\ldots,L_6$ are $1/2$, $\eta$, $\eta$, $2/5$, $|u|$, $v$, and the smallest of them is $2/5$.
\end{proof}

\subsection{Coordinates of the intersection points}

The lines of the starting arrangement are pairwise nonparallel, and none of them is vertical. Therefore the order of the intersection points along a line is the order of their $x$-coordinates. Below we list the $x$-coordinates of all intersection points, including the triple points
\[
P=L_1\cap L_4\cap L_6,
\qquad
Q=L_1\cap L_3\cap L_5.
\]
The orders \eqref{eq:seed-orders} in the next subsection are obtained from them by sorting.

Two lines $y=m(x-a)$ and $y=m'(x-a')$ intersect at the point with coordinates
\begin{equation}\label{eq:pair-point}
x=\frac{ma-m'a'}{m-m'} ,
\qquad
y=\frac{mm'(a-a')}{m-m'} .
\end{equation}

We compute the coordinates of the triple points. By \eqref{eq:pair-point}, the lines $L_1$ and $L_4$ intersect at the point
\[
P=\left(-\frac{\sqrt3}9,\ -\frac{4\sqrt3}9\right),
\]
and the line $L_6:y=v(x-\varepsilon)$ passes through $P$ if and only if
\[
v=\frac{y(P)}{x(P)-\varepsilon}=\frac{-4\sqrt3/9}{-\sqrt3/9-\varepsilon}=\frac{4\sqrt3}{\sqrt3+9\varepsilon} ,
\]
as \eqref{eq:star-slopes} requires.
The lines $L_1$ and $L_3$ intersect at the point
\begin{equation}\label{eq:Q-coords}
Q=\left(\frac{\eta/\sqrt3-\sqrt3/2}{\eta+1/2},\
-\frac{2\eta}{\sqrt3\,(\eta+1/2)}\right),
\end{equation}
and the line $L_5:y=u(x+\varepsilon)$ passes through $Q$ if and only if
\[
\begin{aligned}
u&=\frac{y(Q)}{x(Q)+\varepsilon}
=\frac{-2\eta/\sqrt3}{\eta/\sqrt3-\sqrt3/2+\varepsilon(\eta+1/2)}\\
&=\frac{-2\eta}{\eta-3/2+\varepsilon\sqrt3(\eta+1/2)}
=-\frac{4\eta}{2\eta-3+\varepsilon\sqrt3(2\eta+1)} ,
\end{aligned}
\]
that is, again under condition \eqref{eq:star-slopes}.
Thus $P$ and $Q$ are triple points. There are no other multiple points, since on each line the $x$-coordinates of all other intersection points are distinct (\cref{lem:seed-orders}).

\emph{(i) Intersection points of $Y_0,L_1,L_4,L_5,L_6$.} Apart from $Q$, which we assign to group (iii), these five lines meet at $P$ and at the following points:
\[
\begin{gathered}
Y_0\cap L_1,\ Y_0\cap L_5,\ Y_0\cap L_6,\ Y_0\cap L_4:\quad
x=-\sqrt3,\ -\varepsilon,\ \varepsilon,\ \sqrt3,\\
L_4\cap L_5:\ x=\frac{2\sqrt3/5+u\varepsilon}{2/5-u},
\qquad
L_5\cap L_6:\ x=\frac{\varepsilon\,(u+v)}{v-u} .
\end{gathered}
\]
By \eqref{eq:uv}, the numerator for $L_4\cap L_5$ lies between $\tfrac{2\sqrt3}5-2\varepsilon$ and $\tfrac{2\sqrt3}5$, and the denominator between $\tfrac25+\tfrac95$ and $\tfrac25+2$, so
\[
x(L_4\cap L_5)>\frac{2\sqrt3/5-1/12}{12/5}>\frac14,
\qquad
x(L_4\cap L_5)<\frac{2\sqrt3/5}{11/5}=\frac{2\sqrt3}{11}<\frac13 .
\]
For $L_5\cap L_6$, \eqref{eq:uv} gives $v-u=v+|u|>0$, $u+v=v-|u|>3-2>0$, and
$(v-u)-(u+v)=-2u>0$, so
\[
0<x(L_5\cap L_6)<\varepsilon .
\]

\emph{(ii) Points on $L_2$.} The line $L_2:y=-\eta(x+1/\sqrt3)$ is almost vertical: by \eqref{eq:A-num}, the absolute value of its slope exceeds $2\cdot10^5$.
By \eqref{eq:pair-point}, it intersects a line $D:y=m(x-a)$ at the point with $x$-coordinate
\[
x(L_2\cap D)=-\frac1{\sqrt3}+\frac{g_D}{m+\eta},
\qquad\text{where}\qquad
g_D=m\left(a+\frac1{\sqrt3}\right).
\]
The denominators $m+\eta$ of different lines are almost equal, since $|m|\ll\eta$. Therefore the order of the points on $L_2$ is determined primarily by the numerators $g_D$:
\[
g_{L_5}=u\left(\frac1{\sqrt3}-\varepsilon\right),
\quad
g_{Y_0}=0,
\quad
g_{L_1}=\frac1{\sqrt3},
\quad
g_{L_4}=\frac8{5\sqrt3},
\quad
g_{L_6}=v\left(\frac1{\sqrt3}+\varepsilon\right),
\]
and by \eqref{eq:uv} and \eqref{eq:A-num}
\[
-\frac2{\sqrt3}<g_{L_5}<-\frac95\left(\frac1{\sqrt3}-\frac1{24}\right)<-\frac{24}{25},
\qquad
\sqrt3<g_{L_6}<4\left(\frac1{\sqrt3}+\frac1{24}\right)<\frac52 .
\]
Sorting these values gives
\[
g_{L_5}<-\frac{24}{25}<g_{Y_0}=0<g_{L_1}=\frac1{\sqrt3}
<g_{L_4}=\frac8{5\sqrt3}<\sqrt3<g_{L_6}.
\]
Replacing the denominators $m+\eta$ by $\eta$ does not change this order.
Indeed, $|g_D|<5/2$ for all these lines and $|m|\le4$, so $m+\eta\ge\eta-4>10^5$ by \eqref{eq:A-num} and
\[
\left|\frac{g_D}{m+\eta}-\frac{g_D}{\eta}\right|
=\frac{|g_D|\,|m|}{\eta\,|m+\eta|}
\le\frac{(5/2)\cdot4}{\eta\,(\eta-4)}
<\frac{10^{-4}}\eta .
\]
The difference between consecutive values of $g_D$ is bounded below by $g_{L_4}-g_{L_1}=(8/5-1)/\sqrt3=\sqrt3/5>1/3$, so consecutive values of $g_D/\eta$ differ by at least $1/3\eta$.
Each of them shifts by less than $10^{-4}/\eta$.
Therefore the order of the points on $L_2$ is the same: $L_5,Y_0,L_1,L_4,L_6$. The points themselves lie in the strip
\[
\left|x+\frac1{\sqrt3}\right|\le\frac{5/2}{\eta-4}<\frac1{10^{4}}.
\]

\emph{(iii) Points on $L_3$.} The line $L_3:y=\eta(x-1/\sqrt3)$ is also almost vertical, and similarly
\[
x(L_3\cap D)=\frac1{\sqrt3}+\frac{h_D}{\eta-m},
\qquad\text{where}\qquad
h_D=m\left(\frac1{\sqrt3}-a\right).
\]
As in (ii), the left-to-right order on $L_3$ is therefore determined primarily by the numerators $h_D$:
\[
h_{L_1}=-\frac2{\sqrt3},
\quad\!
h_{L_5}=u\Bigl(\varepsilon+\frac1{\sqrt3}\Bigr),
\quad\!
h_{L_4}=-\frac4{5\sqrt3},
\quad\!
h_{Y_0}=0,
\quad\!
h_{L_6}=v\Bigl(\frac1{\sqrt3}-\varepsilon\Bigr).
\]
Here the lines $L_1$ and $L_5$ give the same point $Q$ (this is what the choice of $u$ achieves), although the values $h_{L_1}$ and $h_{L_5}$ themselves differ.
The denominators $\eta+1/2$ and $\eta-u$ compensate for the difference.
Therefore only four values remain to be compared:
\[
h_{L_1}=-\frac2{\sqrt3}<h_{L_4}=-\frac4{5\sqrt3}<h_{Y_0}=0<\frac32<h_{L_6},
\]
where $h_{L_6}>3\left(\tfrac1{\sqrt3}-\tfrac1{24}\right)>\tfrac32$ and $h_{L_6}<4/\sqrt3<5/2$.
Consecutive values differ by at least $-h_{L_4}=\tfrac4{5\sqrt3}>\tfrac13$, and all of them are less than $5/2$ in absolute value.
The same estimate of the denominators as in (ii) gives the order of the points on $L_3$: first $Q$, then $L_4$, $Y_0$, $L_6$.
All these points lie in the strip $|x-1/\sqrt3|<10^{-4}$.

\emph{(iv) The point $V=L_2\cap L_3$.}
By \eqref{eq:pair-point},
\[
V=\left(0,\ -\frac\eta{\sqrt3}\right)=\left(0,\ -3n_t^6\right).
\]

\emph{(v) Points on $R$.} Since $\eta/\sqrt3=3\delta/2$ by \eqref{eq:params}, formula \eqref{eq:pair-point} gives
\[
\begin{aligned}
x(R\cap L_1)&=\frac{2\delta-\sqrt3}3, & x(R\cap L_2)&=-\frac{\delta}{2(1+\eta)}, & x(R\cap L_3)&=\frac{\delta}{2(\eta-1)},\\
x(R\cap L_4)&=\frac{5\delta-2\sqrt3}3, & x(R\cap L_5)&=\frac{\delta-|u|\varepsilon}{1+|u|}, & x(R\cap L_6)&=-\frac{\delta-v\varepsilon}{v-1},\\
x(R\cap Y_0)&=\delta.
\end{aligned}
\]
This yields the order of the $x$-coordinates on $R$:
\[
x(R\cap L_6)<x(R\cap L_2)<x(R\cap L_3)<x(R\cap L_5)
<x(R\cap L_1)<x(R\cap Y_0)<x(R\cap L_4).
\]
The chain is verified by the following inequalities, which use the estimates \eqref{eq:uv} and \eqref{eq:A-num}:
\[
x(R\cap L_6)=-\frac{\delta-v\varepsilon}{v-1}<-\frac{\delta-1/6}3<-\frac\delta4<-\frac{\sqrt3}9
\qquad\left(v\varepsilon<\frac16,\ v-1<3\right),
\]
\[
-\frac{\sqrt3}9=-\frac\delta{3\sqrt3\,\delta}<-\frac\delta{2(1+\eta)}=x(R\cap L_2)<0
\qquad\left(2(1+\eta)>2\eta=3\sqrt3\,\delta\right),
\]
\[
0<x(R\cap L_3)=\frac\delta{2(\eta-1)}<\frac\delta{5\delta}=\frac15<\frac\delta4
\qquad\left(2(\eta-1)=3\sqrt3\,\delta-2>5\delta\right),
\]
\[
\frac\delta4<\frac{\delta-1/12}3<x(R\cap L_5)=\frac{\delta-|u|\varepsilon}{1+|u|}
<\frac\delta{1+9/5}=\frac{5\delta}{14}
\qquad\left(|u|\varepsilon<\frac1{12}\right),
\]
\[
\begin{gathered}
\frac{5\delta}{14}<x(R\cap L_1)=\frac{2\delta-\sqrt3}3<\delta=x(R\cap Y_0)
<\frac{5\delta-2\sqrt3}3=x(R\cap L_4)\\
\left(13\delta>14\sqrt3,\ \delta>\sqrt3\right).
\end{gathered}
\]
Thus the points $R\cap L_2$ and $R\cap L_3$ have $x$-coordinates in $(-\sqrt3/9,\ 1/5)$, and all other points on $R$ have $|x|>\delta/4$.

The line $R$ has positive slope, so the $y$-coordinates of the points on it are ordered in the same way as the $x$-coordinates, and
\begin{equation}\label{eq:y-neg}
y(R\cap D)\le y(R\cap L_1)=-\frac{\delta+\sqrt3}3<-\frac\delta3
\quad\text{for}\quad
D\in\{L_1,L_2,L_3,L_5,L_6\} .
\end{equation}

\subsection{Intersection orders}

The order of the intersection points along a line is the \emph{local sequence} of the line~\cite{gp84}.
For a simple arrangement, the local sequences of all its lines form its $O$-matrix~\cite[Section~5.1]{pu26b}. Table \eqref{eq:seed-orders} below is an analog of this matrix for an arrangement with triple points.

\begin{lemma}\label{lem:seed-orders}
The left-to-right order of the intersection points on each line is as follows (the entry $L_i/L_j$ denotes a triple point shared with the lines $L_i$ and $L_j$):
\begin{equation}\label{eq:seed-orders}
\begin{array}{c|ccccccc}
R   & L_6 & L_2 & L_3 & L_5 & L_1 & Y_0 & L_4 \\
Y_0 & L_1 & L_2 & L_5 & L_6 & L_3 & L_4 & R \\
L_1 & Y_0 & L_2 & L_4/L_6 & L_3/L_5 & R \\
L_2 & L_5 & Y_0 & L_1 & L_4 & L_6 & R & L_3 \\
L_3 & L_2 & R & L_1/L_5 & L_4 & Y_0 & L_6 \\
L_4 & L_2 & L_1/L_6 & L_5 & L_3 & Y_0 & R \\
L_5 & L_2 & Y_0 & L_6 & L_4 & L_1/L_3 & R \\
L_6 & R & L_2 & L_1/L_4 & L_5 & Y_0 & L_3
\end{array}
\end{equation}
\end{lemma}

\begin{proof}
Each row is obtained by sorting the $x$-coordinates from the previous subsection.
Groups (ii) and (iii) give the rows $L_2$ and $L_3$ without the points $R\cap L_2$, $R\cap L_3$, and $V$.
The points on $R$ other than $R\cap L_2$ and $R\cap L_3$ have $|x|>\delta/4$, while all points of groups (i)--(iv) have $|x|\le\sqrt3<\delta/4$. Therefore these points on $R$ stand at the ends of their rows.
The remaining comparisons are as follows.

\begin{itemize}
\item Row $Y_0$ is given by the order of the zeros of the right-hand sides of \eqref{eq:base-lines}, taking into account $\varepsilon<\tfrac1{\sqrt3}$.
\item Rows $L_5$ and $L_6$: $x(L_5\cap L_6)<\varepsilon<\tfrac14<x(L_4\cap L_5)$.
\item Rows $L_1$, $L_4$, $L_6$: the point $P$ lies between group (ii) and zero, since $-\tfrac1{\sqrt3}+10^{-4}<x(P)=-\tfrac{\sqrt3}9<0$.
\item Rows $L_2$ and $L_3$: $-\tfrac{\sqrt3}9<x(R\cap L_2)<0<x(R\cap L_3)<\tfrac15$, so both points lie to the right of group (ii) and to the left of group (iii), and the point $V$ with $x=0$ separates them.
\item Rows $L_4$ and $L_5$: $x(L_4\cap L_5)<\tfrac13<\tfrac1{\sqrt3}-10^{-4}$.\qedhere
\end{itemize}
\end{proof}

\begin{corollary}\label{cor:rightmost-L4}
The rightmost intersection point of $\mathcal A_0$ on the line $L_4$ is $L_4\cap Y_0$. It is also the rightmost one on $Y_0$.
\end{corollary}

\begin{proof}
Rows $L_4$ and $Y_0$ of table \eqref{eq:seed-orders} with $R$ removed.
\end{proof}

\begin{corollary}\label{cor:rightmost-L2}
Among the intersection points of $\mathcal A_0$ on the line $L_2$, the point $V=L_2\cap L_3$ is the rightmost, and its left neighbor is $L_2\cap L_6$.
Among those on $L_3$, the point $V$ is the leftmost, and its right neighbor is the triple point $Q$.
\end{corollary}

\begin{proof}
Rows $L_2$ and $L_3$ of table \eqref{eq:seed-orders} with $R$ removed.
\end{proof}

\subsection{Bounding box of the intersection points}

\begin{lemma}\label{lem:seed-box}
All intersection points of $\mathcal A_0$ except $V=L_2\cap L_3$ lie in the rectangle
\[
|x|\le\sqrt3,
\qquad
|y|\le8\sqrt3<14 ,
\]
and the point $V$ itself has coordinates $(0,-3n_t^6)$. Moreover, $y(L_2\cap L_6)<-1$ and $y(Q)<-1$.
\end{lemma}

\begin{proof}
The largest $x$-coordinate is $\sqrt3$ (the point $Y_0\cap L_4$), and the smallest is $-\sqrt3$ (the point $Y_0\cap L_1$): groups (ii) and (iii) are within $10^{-4}$ of $\mp1/\sqrt3$, and the other $x$-coordinates of group (i) do not exceed $\sqrt3$ in absolute value.
Each point other than $V$ lies on at least one of the lines $Y_0,L_1,L_4,L_5,L_6$, for which $|m|\le4$ by \cref{lem:uv} and $|a|\le\sqrt3$. Therefore
\[
|y|=|m|\,|x-a|\le4\left(\sqrt3+\sqrt3\right)=8\sqrt3 .
\]
The point $V=(0,-3n_t^6)$ itself is given in (iv).
Finally, by (ii), $x(L_2\cap L_6)<-\tfrac1{\sqrt3}+10^{-4}<-\tfrac12$, so
\[
y(L_2\cap L_6)=v\bigl(x(L_2\cap L_6)-\varepsilon\bigr)<3\cdot\left(-\frac12\right)=-\frac32 ,
\]
and for the triple point $Q$, by \eqref{eq:Q-coords},
\[
y(Q)=-\frac{2\eta}{\sqrt3\,(\eta+1/2)}
=-\frac2{\sqrt3}+\frac1{\sqrt3\,(\eta+1/2)}
<-\frac2{\sqrt3}+10^{-5}<-1 ,
\]
since $\eta>2\cdot10^5$ by \eqref{eq:A-num}.
\end{proof}

\begin{lemma}\label{lem:r-triangles-below}
Each of the four $R$-triangles $RL_1L_5$, $RL_2L_3$, $RL_2L_6$, $RL_3L_5$ lies entirely in the half-plane $y\le-1$, whereas $RY_0L_4$ contains the point $(\delta,0)$.
\end{lemma}

\begin{proof}
By \eqref{eq:y-neg}, all intersection points on $R$ except $R\cap Y_0=(\delta,0)$ and $R\cap L_4$ have $y$-coordinates less than $-\delta/3<-1$.
The third vertices of these four triangles are $Q$, $V$,
$L_2\cap L_6$, and again $Q$, respectively, and their $y$-coordinates are also less than $-1$ by \cref{lem:seed-box}.
Hence all three vertices of each of these triangles lie in the half-plane $y\le-1$, and so does the
whole triangle.

The triangle $RY_0L_4$ has the point $R\cap Y_0=(\delta,0)$ as a vertex.
\end{proof}

\subsection{Counting faces}\label{app:seed-faces}

\begin{lemma}\label{lem:seed-faces}
Three lines of the arrangement bound a triangular face if and
only if, on each of them, the two corresponding intersection points
are distinct and consecutive.
\end{lemma}

\begin{proof}
Denote the three lines by $D_1,D_2,D_3$.
Necessity: the side of the face on $D_1$ is the segment with endpoints $D_1\cap D_2$ and $D_1\cap D_3$, with no intersection points in its interior.
Sufficiency: a line $E\neq D_1,D_2,D_3$ entering the interior of the triangle crosses its boundary at two points, and these cannot both be vertices, since otherwise $E$ would coincide with one of the $D_i$.
Hence $E$ crosses the interior of a side, say the one on $D_1$, that is, between $D_1\cap D_2$ and $D_1\cap D_3$, which contradicts their being consecutive.
\end{proof}

\begin{proof}[Proof of \cref{lem:seed}]
By \cref{lem:seed-faces}, triangular faces correspond to triples of lines whose intersection points are consecutive on all three of these lines.
From \eqref{eq:seed-orders} we obtain exactly 15 such triples:
\[
\begin{gathered}
RY_0L_4,\quad RL_1L_5,\quad RL_2L_3,\quad
RL_2L_6,\quad RL_3L_5,\\
Y_0L_1L_2,\quad Y_0L_2L_5,\quad Y_0L_3L_4,\quad
Y_0L_3L_6,\quad Y_0L_5L_6,\\
L_1L_2L_4,\quad L_1L_4L_5,\quad L_2L_4L_6,\quad
L_3L_4L_5,\quad L_4L_5L_6.
\end{gathered}
\]
The first five contain $R$, and the others belong to $\mathcal A_0$.
Removing $R$ creates no new triangles: the new adjacencies in rows $L_2$ and $L_3$ give no consecutive triples.
By \cref{lem:seed-orders}, $P$ and $Q$ are the only multiple points of $\mathcal A_0$.
\end{proof}

\subsection{Wedges at the first step}

At the first doubling step, $n=n_0=6$. This is the only value in the series with $n\equiv2\pmod4$.
\Cref{lem:right-wedges} lists the right wedges for $n\equiv0\pmod4$ and does not extend to $n=6$.
We treat this case separately.

\begin{lemma}\label{lem:first-step}
Let $n=6$, let the arrangement $\mathcal A$ satisfy \textup{(D1)--(D3)} and \textup{(W1)}, let the slope of the line $L_4$ be positive, and let
\[
\mathcal B=\mathcal A\cup\{M_1,\ldots,M_6\}.
\]
Then the right wedges of $\mathcal B$ that involve at least one new line are exactly the pairs
\begin{equation}\label{eq:wedges-first-step}
(M_1,M_3),\quad (L_4,M_2),\quad (M_4,M_5),\quad (Y_0,M_6).
\end{equation}
There are $4=n/2+1$ of them, and each of the lines $M_1,\ldots,M_6$ belongs to exactly one pair.
Moreover, $\mathcal B$ satisfies \textup{(W1)} with $12$ in place of $6$ and $M_6$ in the role of $L_{10}$.
\end{lemma}

\begin{proof}
Compared with the general case in \cref{lem:right-wedges}, the listing changes in two places.
The auxiliary angles are
\[
\beta_1,\ldots,\beta_6
=-\frac{5\pi}{12},\ -\frac\pi4,\ -\frac\pi{12},\
\frac\pi{12},\ \frac\pi4,\ \frac{5\pi}{12},
\]
and condition \eqref{eq:Par1} $i+j=n/2+1=4$ holds for $i=1$, $j=3$, and also for $i=j=2$ ($\beta_1 + \beta_3 = 2\beta_2=-\pi/2$), so the line $M_2$ has no almost parallel partner.
Similarly, \eqref{eq:Par2} reduces to $i+j=10$ and gives the index pair $(4,6)$, with $M_5$ left without a partner at $i=j=5$.
In addition, $3n/4=9/2$ is not an integer, so the pair $(L_{n-2},M_{3n/4})$ of \eqref{eq:wedges-old} does not exist for $n=6$.
Its place is taken by the pair $(L_4,M_2)$, and the slope of $L_4$ is positive rather than negative as in the general case.

The rightmost intersection points are found by case analysis, as in \cref{lem:right-wedges}. The lines $M_1$ and $M_3$ satisfy \eqref{eq:Par1}, so
\[
x(M_1\cap M_3)=-n^6\sin(2\beta_1)=\frac{n^6}2>n^5.
\]
The lines $M_4$ and $M_5$ form a pair from \eqref{eq:wedges-high} with $i+j=3n/2=9$. Their auxiliary angles sum to $\pi/2-\pi/6=\pi/3$, and by \eqref{eq:M-M-coord} and \eqref{eq:theta}
\[
x(M_4\cap M_5)=\tan\frac\pi3+\theta=a_4+\theta,
\qquad \frac3{n^8}<\theta<\frac1{n^2}.
\]
The line $M_6$ intersects $Y_0$ at $x=b_6=\cot(\pi/{12})=a_4+2$.

It remains to list the points on $M_2$. For its intersections with $M_3,M_4,M_5,M_6$, the sums of auxiliary angles are $\beta_2+\beta_j=-\pi/3,-\pi/6,0,\pi/6$. By \eqref{eq:M-M-coord} and \eqref{eq:theta}, these intersection points satisfy
\[
x(M_2\cap M_j)=\tan(\beta_2+\beta_j)+\theta_j\le\tan\frac\pi6+\frac1{n^2}<1,
\qquad |\theta_j|<\frac1{n^2}
\]
(exactly $x=0$ for $j=5$). The line $M_2$ intersects $Y_0$ at the point $x=b_2=-1$, and by \eqref{eq:L-M-coord} the old lines $L_k\neq L_4$ at points with
\[
x(M_2\cap L_k)\le a_3+\frac3{n^9}<1.
\]
For the point $M_1\cap M_2$ the sum is $-2\pi/3$, and the term $\theta$ in \eqref{eq:M-M-coord} is negative by \eqref{eq:theta} (all three factors $\beta_1+\beta_2$, $\beta_1$, $\beta_2$ are negative), so
\[
x(M_1\cap M_2)<\tan\Bigl(-\frac{2\pi}3\Bigr)-\frac3{n^8}=a_4-\frac3{n^8}.
\]

The intersection point $L_4\cap M_2$, by contrast, lies to the right of $a_4$.
We argue as for the pair $(L_{n-2},M_{3n/4})$ in \cref{lem:right-wedges}, but now the slope $\ell$ of $L_4$ is positive, so $L_4$ lies below each line $M_j$ to the left of their common point and above it to the right.
Hence $L_4\cap M_2$ is the rightmost point on $L_4$ among the intersection points with the new lines if at its $x$-coordinate the line $M_2$ is the highest of the new lines. By \eqref{eq:M-height}, as in \eqref{eq:F-gap}, this holds if $F_k-F_2>12/n^9$ for all $k\neq2$.
By \eqref{eq:G-j} with $n=6$, the leading term is $G_j=-2\sin(\pi j/3)-1$. It attains its smallest value $-1-\sqrt3$ for two values, $j=1$ and $j=2$, which gives a tie (in the general case, by contrast, the maximum is separated from the closest values by $2\sin(\pi/n)$).
The term $H_j$ from the proof of \cref{lem:right-wedges} breaks the tie:
\[
\begin{gathered}
H_1=-\frac{2(\sqrt3-1)}{n^6},\qquad H_2=-\frac{\sqrt3+1}{n^6},\\
F_1-F_2=\frac{3-\sqrt3}{n^{6}}>2.7\cdot10^{-5}>20\cdot\frac{12}{n^9}.
\end{gathered}
\]
For the other $k$ the gap is at least $1$ already in the leading term.
By (W1), there are no old intersection points on $L_4$ to the right of $a_4$. Moreover, $F_2<0$, because its leading term is $G_2=-1-\sqrt3$, so $M_2$ passes above the point $(a_4,0)$ of $L_4$, and $L_4\cap M_2$ lies to the right of $a_4$. Hence it is the rightmost point on $L_4$. It is also the rightmost point on $M_2$, since all other intersection points of $M_2$ lie to the left of $a_4$.
Thus the right wedges with a new line are exactly the pairs \eqref{eq:wedges-first-step}. There are four of them, and each index $1,\ldots,6$ is covered once.
In \textup{(W1)}, the role of $L_4$ with $a_4=\cot(\pi/6)$ passes to $M_6$ with $b_6=\cot(\pi/12)$.
\end{proof}

\section{Arrangements attaining the conjectured bound}\label{app:figures}

This appendix collects the arrangements referred to in the main text, in increasing order of $n$.
The optimal arrangements for $n=6,10,12$ (\cref{fig:n6-all,fig:n10,fig:n10-all,fig:n12-all}) give the lower bounds of \cref{thm:small-n}; the two non-stretchable optima for $n=10$ are drawn with pseudolines (\cref{fig:n10-nonstretchable}).
\Cref{fig:n26} shows the member $\mathcal C_2$ of the series of \cref{thm:series}.
Augmentation (\cref{sec:augment}) gives the arrangements of 20, 32 and 44 pseudolines in \cref{fig:n20,fig:n20-collapse,fig:n32,fig:n44}.
Doubling the arrangement of \cref{fig:n20-collapse} gives the pseudoline series $18\cdot2^t+2$ discussed in \cref{sec:discussion}, whose members with $n=38,74$ are stretchable (\cref{fig:n38,tab:n74-lines}), and doubling that of \cref{fig:n32} gives the series $30\cdot2^t+2$.
\Cref{tab:n20-collapse-lines,tab:n32-lines,tab:n38-lines,tab:n74-lines} give the equations of the stretched arrangements.

Different shades correspond to different phases of the local coloring (\cref{fig:phase-seam}); unbounded triangles are lighter.
The boundaries between the phases are thus the unused and double segments, that is, the part of the phase seam (\cref{def:seam}) visible in the affine chart.
\Cref{fig:n44} is not colored.

The coefficients $m_i,b_i$ in the tables are exact.
We checked in rational arithmetic that for these values the triple points are exact, the order of the crossings on every line agrees with the word of the arrangement, and the number of triangles is as stated.

For brevity, letters in the words are written as their initial slot followed by $m-2$ asterisks: a simple crossing $(g,2)\mapsto g$, a triple point $(g,3)\mapsto g*$, a crossing of multiplicity 4 $(g,4)\mapsto g**$, and so on.

\newcommand\ngap{\hspace{0.01\linewidth}}
\newcommand\nsix[1]{\begin{minipage}{0.326\linewidth}\input{figures/n6/6-#1}\end{minipage}}

\begin{figure}[t]
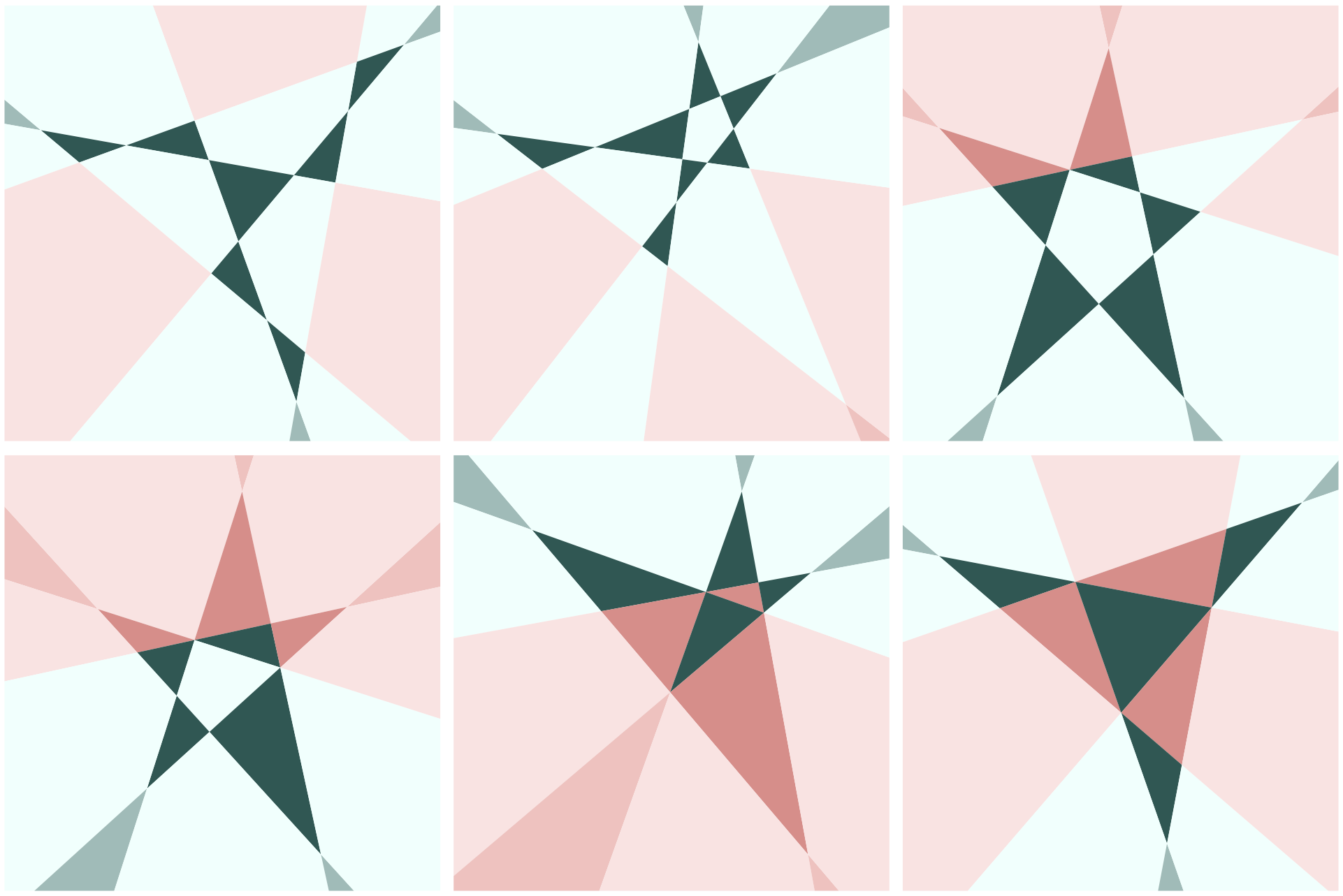

\centering
\begin{minipage}{0.326\linewidth}\input{figures/n6/6-1330sgdjqf8p8}\end{minipage}\ngap\begin{minipage}{0.326\linewidth}\input{figures/n6/6-1nevwww2r41us}\end{minipage}\ngap\begin{minipage}{0.326\linewidth}\input{figures/n6/6-3asa5zegavnkz}\end{minipage}\par\vspace{1mm}
\begin{minipage}{0.326\linewidth}\input{figures/n6/6-3pkh1zclcef98}\end{minipage}\ngap\begin{minipage}{0.326\linewidth}\input{figures/n6/6-2d5kqooycn0g3}\end{minipage}\ngap\begin{minipage}{0.326\linewidth}\input{figures/n6/6-3aiagxuwsbcdh}\end{minipage}
\caption{All 6 optimal arrangements for $n=6$ ($a_3=7$), with $t_3=0,0,1$ in the first row and $t_3=2,3,3$ in the second.}
\label{fig:n6-all}
\end{figure}

\begin{figure}[p]
\centering
\begin{minipage}{0.685\linewidth}\input{figures/n10/27}\end{minipage}\par\vspace{1mm}
\begin{minipage}{0.685\linewidth}\input{figures/n10/65}\end{minipage}
\caption{Optimal arrangements for $n=10$ with 10 triple points. Top: \texttt{0 2 4 6 5 4 3 1* 3* 6 5 8 6* 4* 2* 0* 8 6* 4* 2* 1 4 3 5* 7}. Bottom: \texttt{0 2 4 6 8 7 6 5 3* 1* 0 5* 3* 1* 8 7 5* 3* 2 1 7* 5* 3* 5 7}.}
\label{fig:n10}
\end{figure}

\begin{figure}[p]
\centering
\begin{minipage}{0.785\linewidth}\input{figures/n10/non-stretchable-1}\end{minipage}\par\vspace{1mm}
\begin{minipage}{0.785\linewidth}\input{figures/n10/non-stretchable-2}\end{minipage}
\caption{Non-stretchable optimal arrangements for $n=10$ with 10 triple points, drawn with pseudolines. Top: \texttt{0 1 3 6 7 5* 4 2* 0* 4* 2* 6* 4* 8 6* 5 3* 1* 0 3 4 6 5 4 6}. Bottom: \texttt{0 1 4 3 6 8 7 6 4* 2* 0* 6 4* 2* 1 8 6* 4* 2* 1 8 6* 4* 3 6}.
  Non-stretchability was verified by the method of biquadratic final polynomials~\cite{br90}; the certificates are in the repository~\cite{nonsimple-arrangements}.}
\label{fig:n10-nonstretchable}
\end{figure}

\newcommand\nten[1]{\begin{minipage}{0.192\linewidth}\input{figures/n10/#1}\end{minipage}}

\begin{figure}[p]
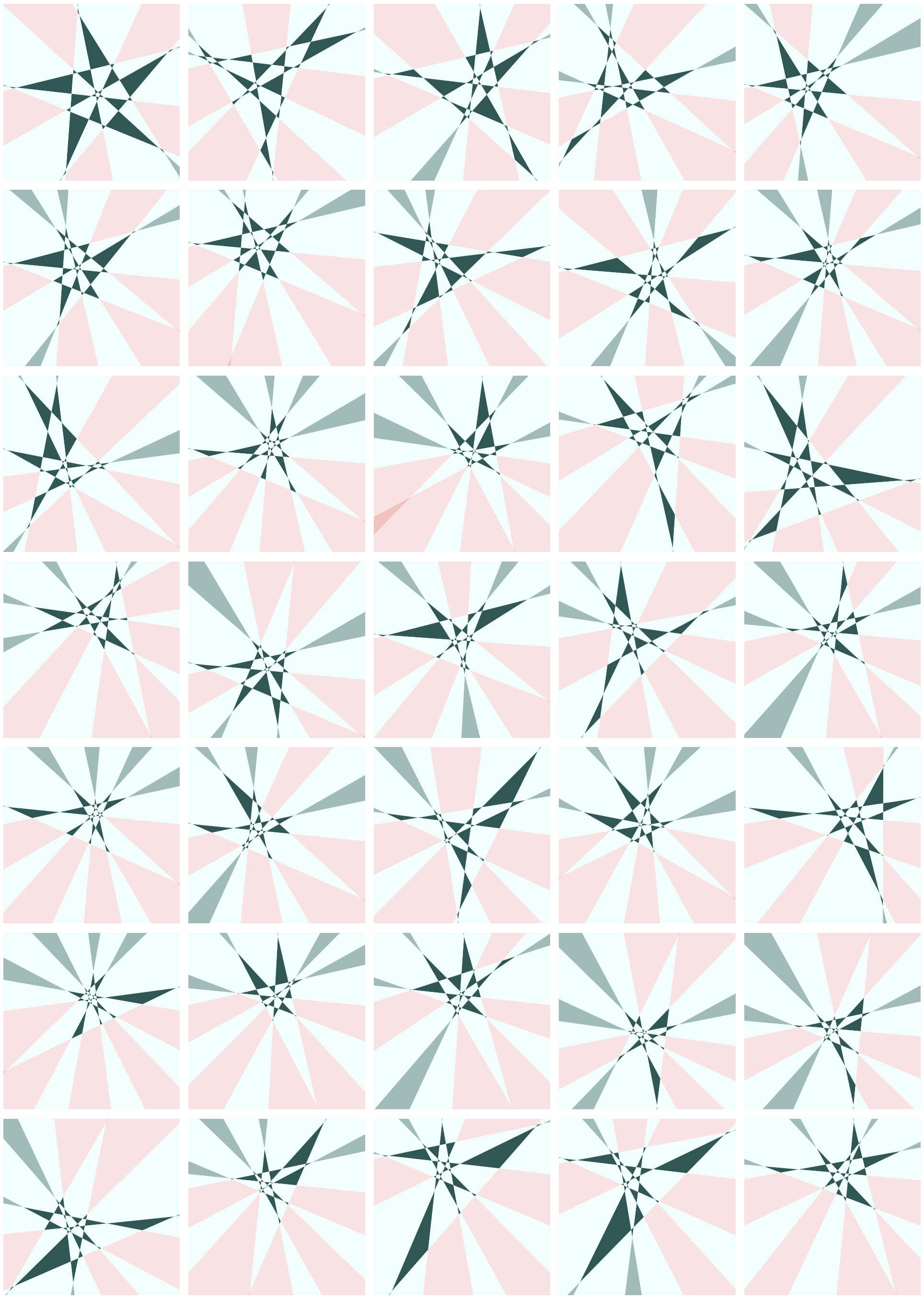

\noindent\begin{minipage}{0.192\linewidth}\input{figures/n10/10}\end{minipage}\ngap\begin{minipage}{0.192\linewidth}\input{figures/n10/14}\end{minipage}\ngap\begin{minipage}{0.192\linewidth}\input{figures/n10/15}\end{minipage}\ngap\begin{minipage}{0.192\linewidth}\input{figures/n10/16}\end{minipage}\ngap\begin{minipage}{0.192\linewidth}\input{figures/n10/24}\end{minipage}\par\vspace{1mm}
\noindent\begin{minipage}{0.192\linewidth}\input{figures/n10/25}\end{minipage}\ngap\begin{minipage}{0.192\linewidth}\input{figures/n10/26}\end{minipage}\ngap\begin{minipage}{0.192\linewidth}\input{figures/n10/33}\end{minipage}\ngap\begin{minipage}{0.192\linewidth}\input{figures/n10/37}\end{minipage}\ngap\begin{minipage}{0.192\linewidth}\input{figures/n10/38}\end{minipage}\par\vspace{1mm}
\noindent\begin{minipage}{0.192\linewidth}\input{figures/n10/40}\end{minipage}\ngap\begin{minipage}{0.192\linewidth}\input{figures/n10/41}\end{minipage}\ngap\begin{minipage}{0.192\linewidth}\input{figures/n10/42}\end{minipage}\ngap\begin{minipage}{0.192\linewidth}\input{figures/n10/44}\end{minipage}\ngap\begin{minipage}{0.192\linewidth}\input{figures/n10/45}\end{minipage}\par\vspace{1mm}
\noindent\begin{minipage}{0.192\linewidth}\input{figures/n10/47}\end{minipage}\ngap\begin{minipage}{0.192\linewidth}\input{figures/n10/48}\end{minipage}\ngap\begin{minipage}{0.192\linewidth}\input{figures/n10/49}\end{minipage}\ngap\begin{minipage}{0.192\linewidth}\input{figures/n10/50}\end{minipage}\ngap\begin{minipage}{0.192\linewidth}\input{figures/n10/51}\end{minipage}\par\vspace{1mm}
\noindent\begin{minipage}{0.192\linewidth}\input{figures/n10/52}\end{minipage}\ngap\begin{minipage}{0.192\linewidth}\input{figures/n10/54}\end{minipage}\ngap\begin{minipage}{0.192\linewidth}\input{figures/n10/56}\end{minipage}\ngap\begin{minipage}{0.192\linewidth}\input{figures/n10/57}\end{minipage}\ngap\begin{minipage}{0.192\linewidth}\input{figures/n10/58}\end{minipage}\par\vspace{1mm}
\noindent\begin{minipage}{0.192\linewidth}\input{figures/n10/59}\end{minipage}\ngap\begin{minipage}{0.192\linewidth}\input{figures/n10/60}\end{minipage}\ngap\begin{minipage}{0.192\linewidth}\input{figures/n10/61}\end{minipage}\ngap\begin{minipage}{0.192\linewidth}\input{figures/n10/62}\end{minipage}\ngap\begin{minipage}{0.192\linewidth}\input{figures/n10/63}\end{minipage}\par\vspace{1mm}
\noindent\begin{minipage}{0.192\linewidth}\input{figures/n10/64}\end{minipage}\ngap\begin{minipage}{0.192\linewidth}\input{figures/n10/66}\end{minipage}\ngap\begin{minipage}{0.192\linewidth}\input{figures/n10/67}\end{minipage}\ngap\begin{minipage}{0.192\linewidth}\input{figures/n10/68}\end{minipage}\ngap\begin{minipage}{0.192\linewidth}\input{figures/n10/69}\end{minipage}
\caption{Optimal arrangements for $n=10$ ($a_3=25$): the stretchable arrangements other than those in \cref{fig:n10}, in increasing order of the number $t_3$ of triple points. Here $t_3=0$.}
\label{fig:n10-all}
\end{figure}

\begin{figure}[p]
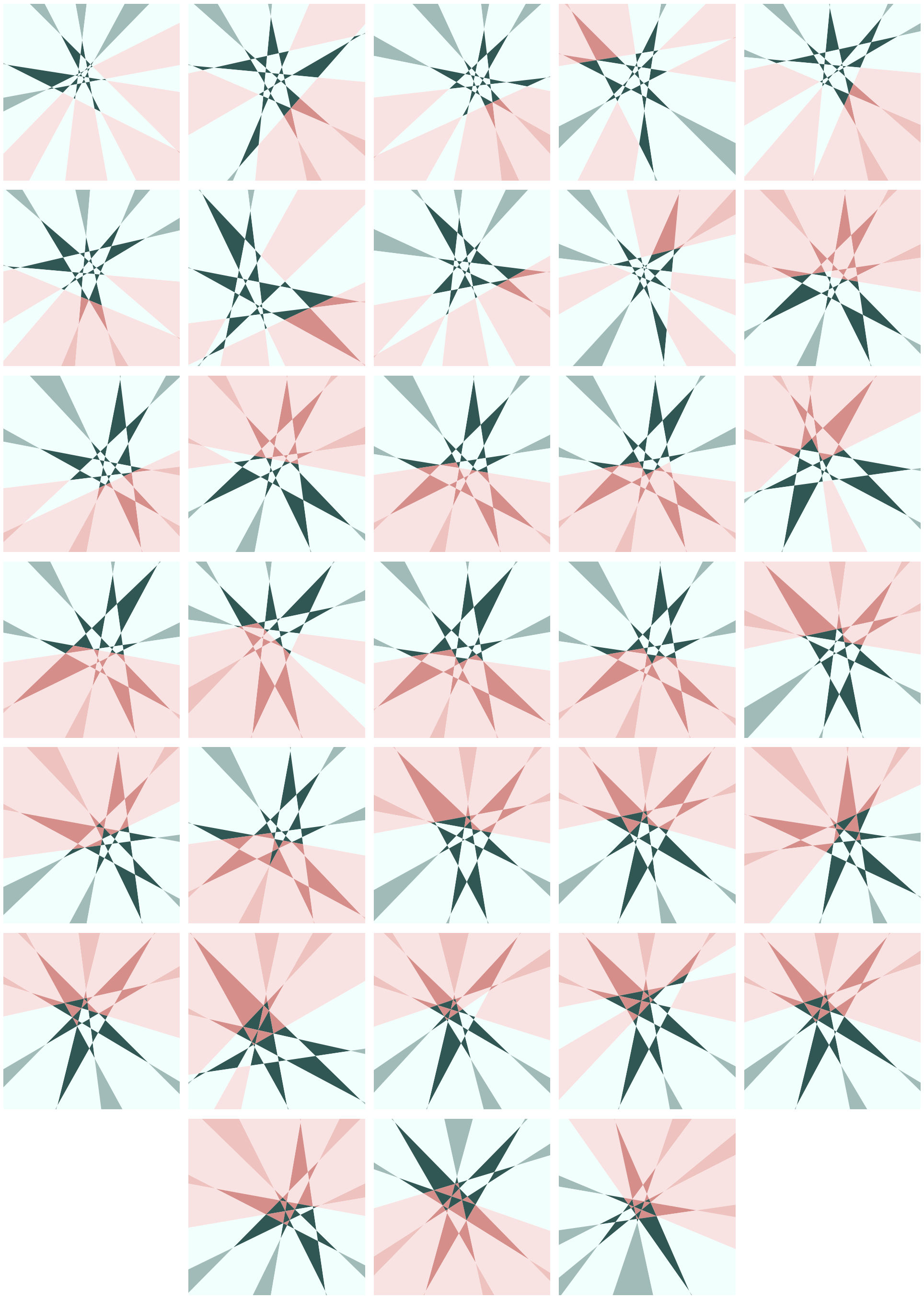

\ContinuedFloat
\noindent\begin{minipage}{0.192\linewidth}\input{figures/n10/70}\end{minipage}\ngap\begin{minipage}{0.192\linewidth}\input{figures/n10/07}\end{minipage}\ngap\begin{minipage}{0.192\linewidth}\input{figures/n10/17}\end{minipage}\ngap\begin{minipage}{0.192\linewidth}\input{figures/n10/18}\end{minipage}\ngap\begin{minipage}{0.192\linewidth}\input{figures/n10/21}\end{minipage}\par\vspace{1mm}
\noindent\begin{minipage}{0.192\linewidth}\input{figures/n10/23}\end{minipage}\ngap\begin{minipage}{0.192\linewidth}\input{figures/n10/30}\end{minipage}\ngap\begin{minipage}{0.192\linewidth}\input{figures/n10/39}\end{minipage}\ngap\begin{minipage}{0.192\linewidth}\input{figures/n10/55}\end{minipage}\ngap\begin{minipage}{0.192\linewidth}\input{figures/n10/04}\end{minipage}\par\vspace{1mm}
\noindent\begin{minipage}{0.192\linewidth}\input{figures/n10/13}\end{minipage}\ngap\begin{minipage}{0.192\linewidth}\input{figures/n10/22}\end{minipage}\ngap\begin{minipage}{0.192\linewidth}\input{figures/n10/03}\end{minipage}\ngap\begin{minipage}{0.192\linewidth}\input{figures/n10/05}\end{minipage}\ngap\begin{minipage}{0.192\linewidth}\input{figures/n10/12}\end{minipage}\par\vspace{1mm}
\noindent\begin{minipage}{0.192\linewidth}\input{figures/n10/28}\end{minipage}\ngap\begin{minipage}{0.192\linewidth}\input{figures/n10/29}\end{minipage}\ngap\begin{minipage}{0.192\linewidth}\input{figures/n10/02}\end{minipage}\ngap\begin{minipage}{0.192\linewidth}\input{figures/n10/11}\end{minipage}\ngap\begin{minipage}{0.192\linewidth}\input{figures/n10/35}\end{minipage}\par\vspace{1mm}
\noindent\begin{minipage}{0.192\linewidth}\input{figures/n10/01}\end{minipage}\ngap\begin{minipage}{0.192\linewidth}\input{figures/n10/06}\end{minipage}\ngap\begin{minipage}{0.192\linewidth}\input{figures/n10/19}\end{minipage}\ngap\begin{minipage}{0.192\linewidth}\input{figures/n10/53}\end{minipage}\ngap\begin{minipage}{0.192\linewidth}\input{figures/n10/09}\end{minipage}\par\vspace{1mm}
\noindent\begin{minipage}{0.192\linewidth}\input{figures/n10/20}\end{minipage}\ngap\begin{minipage}{0.192\linewidth}\input{figures/n10/34}\end{minipage}\ngap\begin{minipage}{0.192\linewidth}\input{figures/n10/43}\end{minipage}\ngap\begin{minipage}{0.192\linewidth}\input{figures/n10/08}\end{minipage}\ngap\begin{minipage}{0.192\linewidth}\input{figures/n10/36}\end{minipage}\par\vspace{1mm}
\noindent\begin{minipage}{0.192\linewidth}\input{figures/n10/46}\end{minipage}\ngap\begin{minipage}{0.192\linewidth}\input{figures/n10/31}\end{minipage}\ngap\begin{minipage}{0.192\linewidth}\input{figures/n10/32}\end{minipage}
\caption[]{(continued) The first picture is the last one with $t_3=0$. It is followed, row by row, by $t_3=1$ (8 arrangements), $2$ (3), $3$ (5), $4$ (3), $5$ (3), $6$ (1), $7$ (4), $8$ (3), $9$ (2).}
\end{figure}

\newcommand\ntw[1]{\begin{minipage}{0.2425\linewidth}\input{figures/n12/#1}\end{minipage}}

\begin{figure}[p]
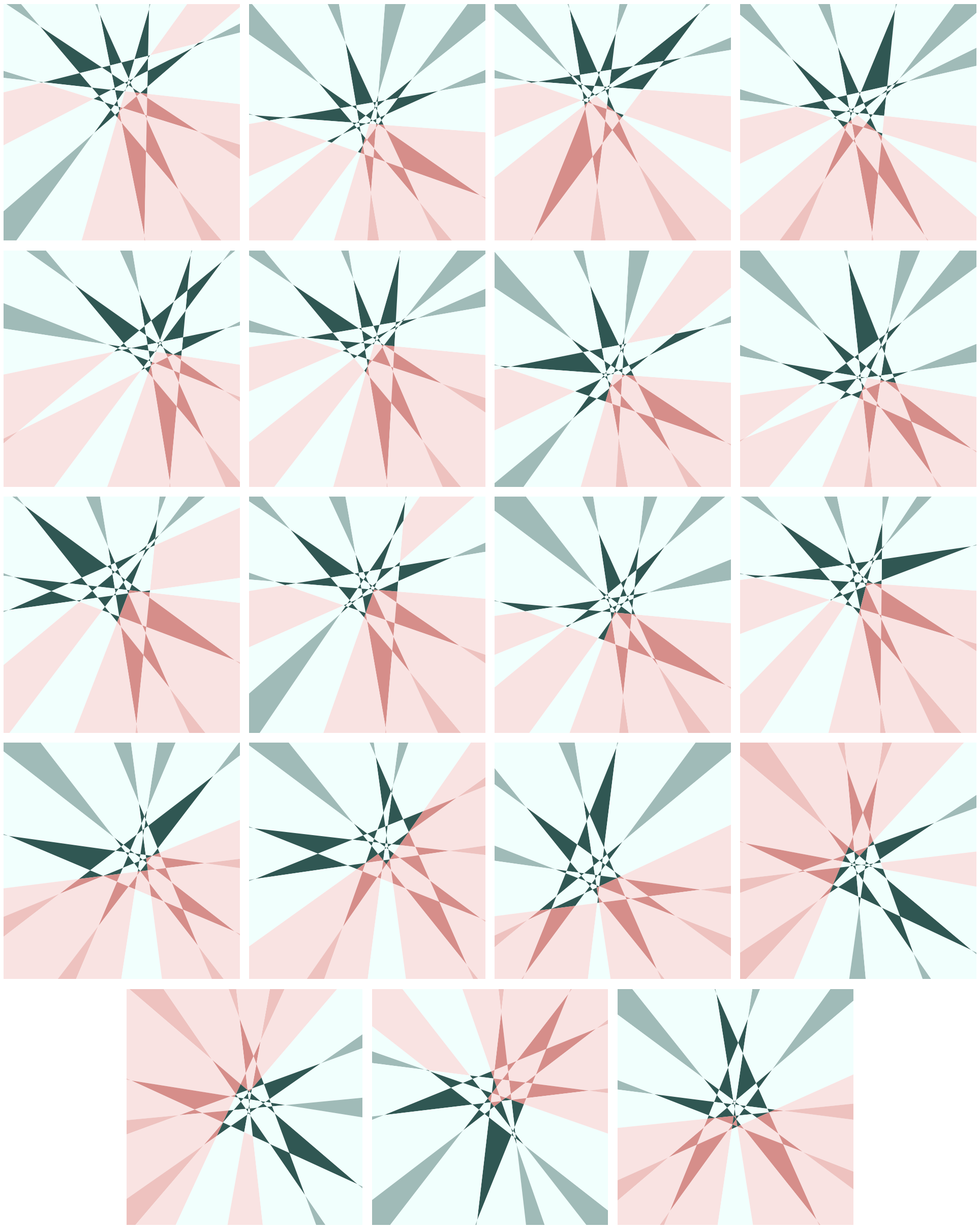

\noindent\begin{minipage}{0.2425\linewidth}\input{figures/n12/05}\end{minipage}\ngap\begin{minipage}{0.2425\linewidth}\input{figures/n12/06}\end{minipage}\ngap\begin{minipage}{0.2425\linewidth}\input{figures/n12/09}\end{minipage}\ngap\begin{minipage}{0.2425\linewidth}\input{figures/n12/13}\end{minipage}\par\vspace{1mm}
\noindent\begin{minipage}{0.2425\linewidth}\input{figures/n12/14}\end{minipage}\ngap\begin{minipage}{0.2425\linewidth}\input{figures/n12/15}\end{minipage}\ngap\begin{minipage}{0.2425\linewidth}\input{figures/n12/17}\end{minipage}\ngap\begin{minipage}{0.2425\linewidth}\input{figures/n12/18}\end{minipage}\par\vspace{1mm}
\noindent\begin{minipage}{0.2425\linewidth}\input{figures/n12/02}\end{minipage}\ngap\begin{minipage}{0.2425\linewidth}\input{figures/n12/08}\end{minipage}\ngap\begin{minipage}{0.2425\linewidth}\input{figures/n12/16}\end{minipage}\ngap\begin{minipage}{0.2425\linewidth}\input{figures/n12/19}\end{minipage}\par\vspace{1mm}
\noindent\begin{minipage}{0.2425\linewidth}\input{figures/n12/01}\end{minipage}\ngap\begin{minipage}{0.2425\linewidth}\input{figures/n12/03}\end{minipage}\ngap\begin{minipage}{0.2425\linewidth}\input{figures/n12/04}\end{minipage}\ngap\begin{minipage}{0.2425\linewidth}\input{figures/n12/11}\end{minipage}\par\vspace{1mm}
\noindent\begin{minipage}{0.2425\linewidth}\input{figures/n12/12}\end{minipage}\ngap\begin{minipage}{0.2425\linewidth}\input{figures/n12/10}\end{minipage}\ngap\begin{minipage}{0.2425\linewidth}\input{figures/n12/07}\end{minipage}
\caption{All 19 optimal arrangements for $n=12$ ($a_3=38$) in increasing order of the number of triple points: $t_3=2$ (the first two rows), $3$ (the third row), $4$ (the fourth row and the first picture of the fifth), $5$ and $6$.}
\label{fig:n12-all}
\end{figure}

\begin{figure}[p]
\noindent\input{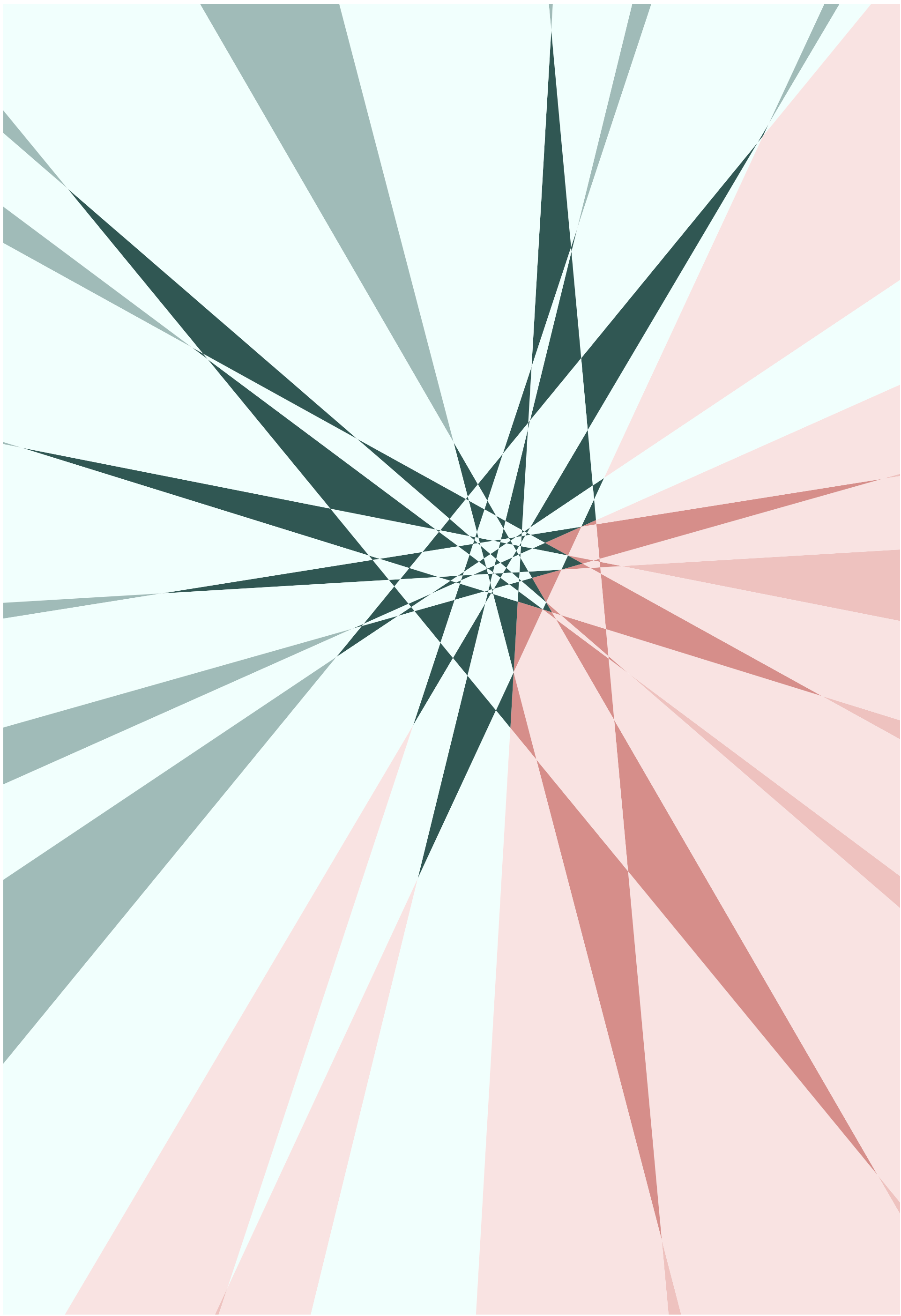}
\caption{An arrangement of 20 lines with $a_3=117$ triangles and 6 triple points, obtained by direct augmentation (\cref{sec:augment}).}
\label{fig:n20}
\end{figure}

\begin{figure}[p]
\noindent\input{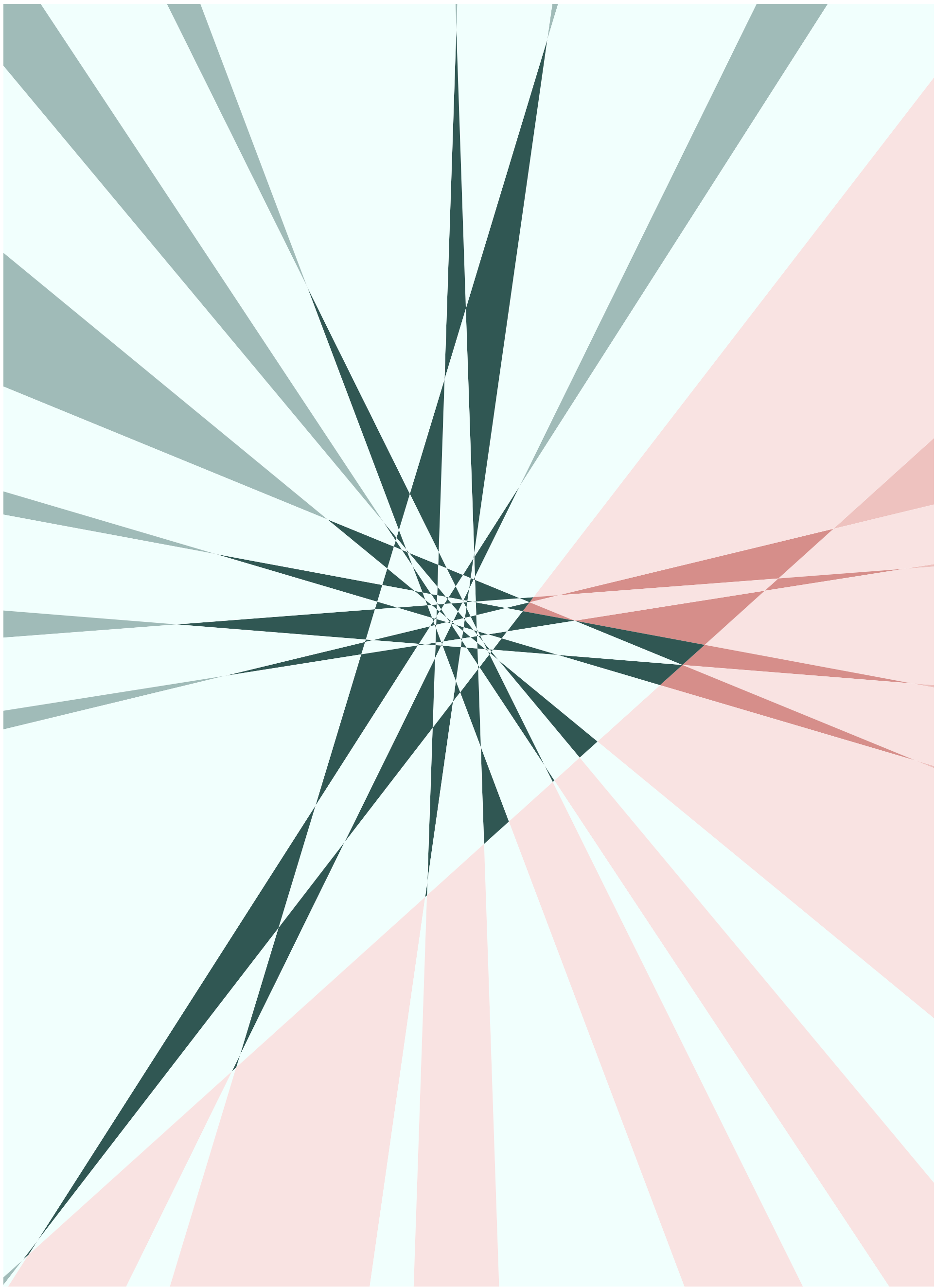}
\caption{An arrangement of 20 lines with $a_3=117$ triangles and 3 triple points, obtained by collapse augmentation (\cref{sec:augment}).}
\label{fig:n20-collapse}
\end{figure}

\begin{table}[t]
\centering
\small
\caption{Equations $y=m_ix+b_i$ of the lines of the arrangement in \cref{fig:n20-collapse}. The triple points lie on line 13 together with the pairs of lines $\{5,7\}$, $\{6,10\}$, $\{8,11\}$.}
\label{tab:n20-collapse-lines}
\renewcommand{\arraystretch}{1.25}
\begin{tabular}{rrr@{\qquad}rrr}
\toprule
$i$ & $m_i$ & $b_i$ & $i$ & $m_i$ & $b_i$ \\
\midrule
0 & $-0.034$ & $-0.8668$ & 10 & $12.4722$ & $-49.01273808$ \\
1 & $-0.1424$ & $0.7687$ & 11 & $5.3816$ & $-8.8336011$ \\
2 & $-0.3026$ & $-3.2348$ & 12 & $3.3902$ & $-9.8084$ \\
3 & $-0.4915$ & $1.1772$ & 13 & $2.4386$ & $-0.9237$ \\
4 & $-0.639$ & $0.6319$ & 14 & $1.2157$ & $-3.5408$ \\
5 & $-0.7707$ & $4.87422138$ & 15 & $0.833$ & $-1.7987$ \\
6 & $-1.1086$ & $16.07732016$ & 16 & $0.6615$ & $-2.1437$ \\
7 & $-4.1366$ & $10.95505632$ & 17 & $0.4958$ & $-0.5564$ \\
8 & $-6.343$ & $22.67860632$ & 18 & $0.3773$ & $-2.1232$ \\
9 & $-13.0923$ & $24.1398$ & 19 & $0.034$ & $0.8634$ \\
\bottomrule
\end{tabular}
\end{table}

\begin{figure}[p]
\noindent\input{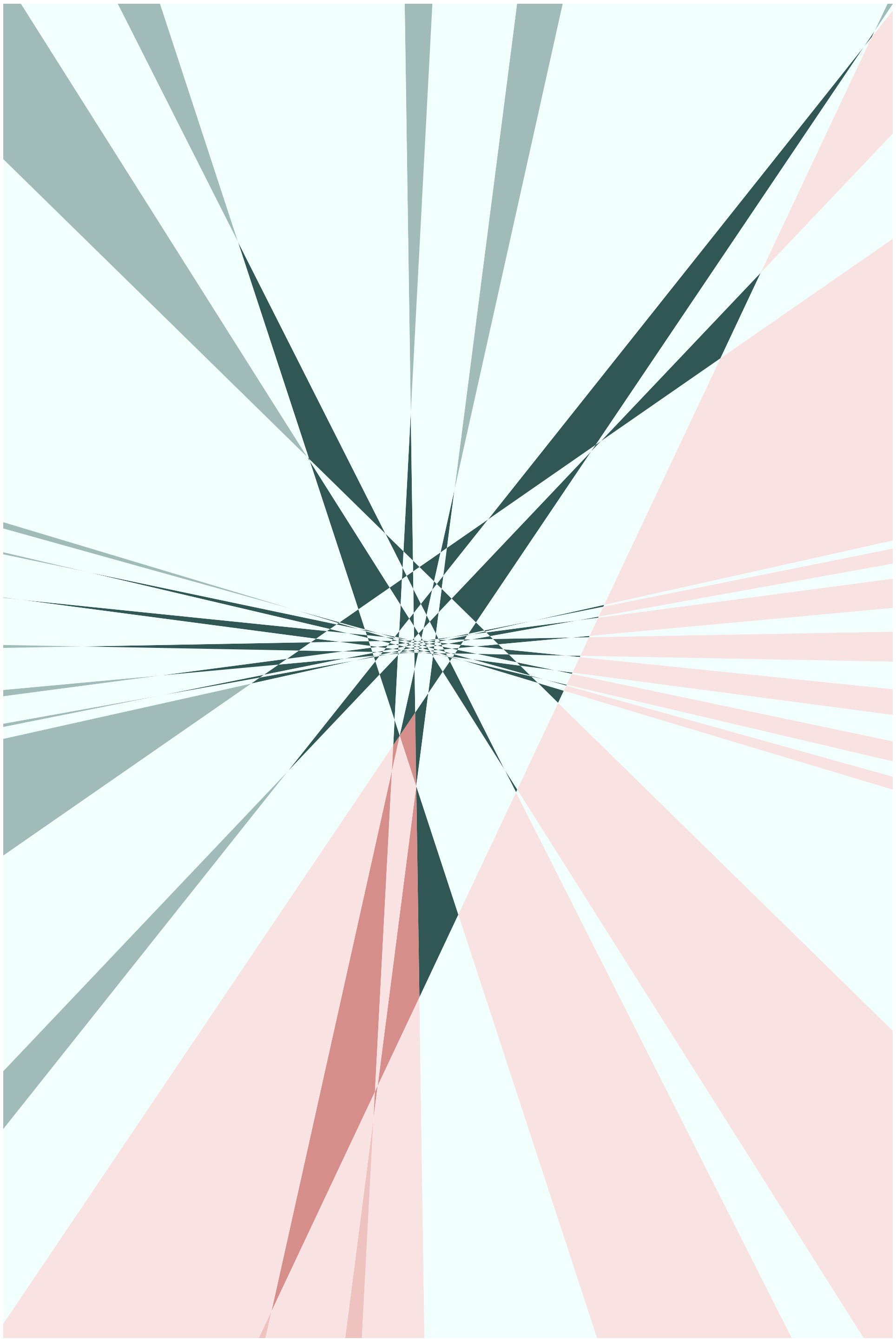}
\caption{An arrangement of 26 lines with $a_3=204$ triangles and 2 triple points, combinatorially equivalent to $\mathcal C_2$ of \cref{thm:series}. The positions of the lines are changed for readability.}
\label{fig:n26}
\end{figure}

\begin{figure}[p]
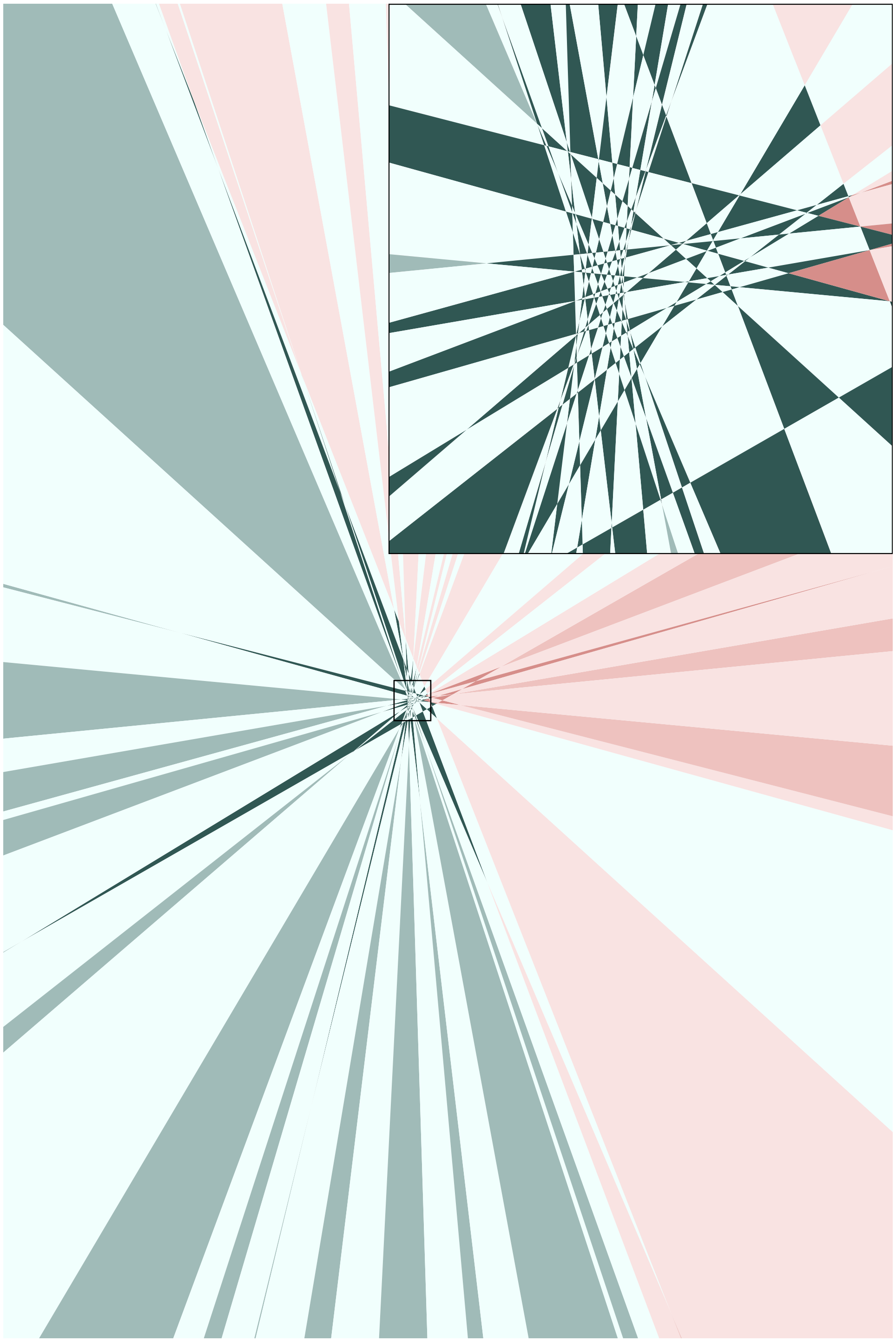

\noindent\input{figures/n32/32-31imvjadxi50t}\unskip\llap{\raisebox{0.882\linewidth}{\begin{minipage}[b]{0.5665\linewidth}\input{figures/n32/32-31imvjadxi50t-zoom}\end{minipage}}}
\caption{An arrangement of 32 lines with $a_3=315$ triangles and 4 triple points, obtained by direct augmentation (\cref{sec:augment}). Upper right: the framed region, enlarged.}
\label{fig:n32}
\end{figure}

\begin{table}[t]
\centering
\small
\caption{Equations $y=m_ix+b_i$ of the lines of the arrangement in \cref{fig:n32}. The triple points lie on line 5 together with the pairs of lines $\{0,2\}$, $\{1,31\}$, $\{25,27\}$, $\{28,30\}$.}
\label{tab:n32-lines}
\renewcommand{\arraystretch}{1.15}
\begin{tabular}{rrr@{\qquad}rrr}
\toprule
$i$ & $m_i$ & $b_i$ & $i$ & $m_i$ & $b_i$ \\
\midrule
0 & $-0.06476$ & $-0.7560626464$ & 16 & $4.11657$ & $12.85756$ \\
1 & $-0.1772$ & $4.0055813476$ & 17 & $2.86853$ & $15.03237$ \\
2 & $-0.19072$ & $1.096562072$ & 18 & $2.80636$ & $9.09436$ \\
3 & $-0.62491$ & $0.88298$ & 19 & $2.33795$ & $10.55365$ \\
4 & $-1.61302$ & $-8.02093$ & 20 & $2.12481$ & $7.49822$ \\
5 & $-1.74892$ & $24.01463$ & 21 & $1.85599$ & $7.50283$ \\
6 & $-1.82634$ & $6.28754$ & 22 & $1.16509$ & $-2.18024$ \\
7 & $-1.91825$ & $-10.33657$ & 23 & $0.59231$ & $0.48538$ \\
8 & $-2.06335$ & $-8.8103$ & 24 & $0.54137$ & $-2.54997$ \\
9 & $-2.09887$ & $-12.3846$ & 25 & $0.41874$ & $-1.8998753596$ \\
10 & $-3.7305$ & $-13.36831$ & 26 & $0.39469$ & $-10.62239$ \\
11 & $-5.84564$ & $-35.30939$ & 27 & $0.26104$ & $-0.0145623976$ \\
12 & $-7.7291$ & $-25.79848$ & 28 & $0.19652$ & $-1.9637835856$ \\
13 & $-22.151$ & $-144.70811$ & 29 & $0.18978$ & $0.89411$ \\
14 & $13.67108$ & $42.57195$ & 30 & $0.11893$ & $-0.9276862965$ \\
15 & $5.69142$ & $32.4434$ & 31 & $0.06777$ & $0.8869491177$ \\
\bottomrule
\end{tabular}
\end{table}

\begin{figure}[p]
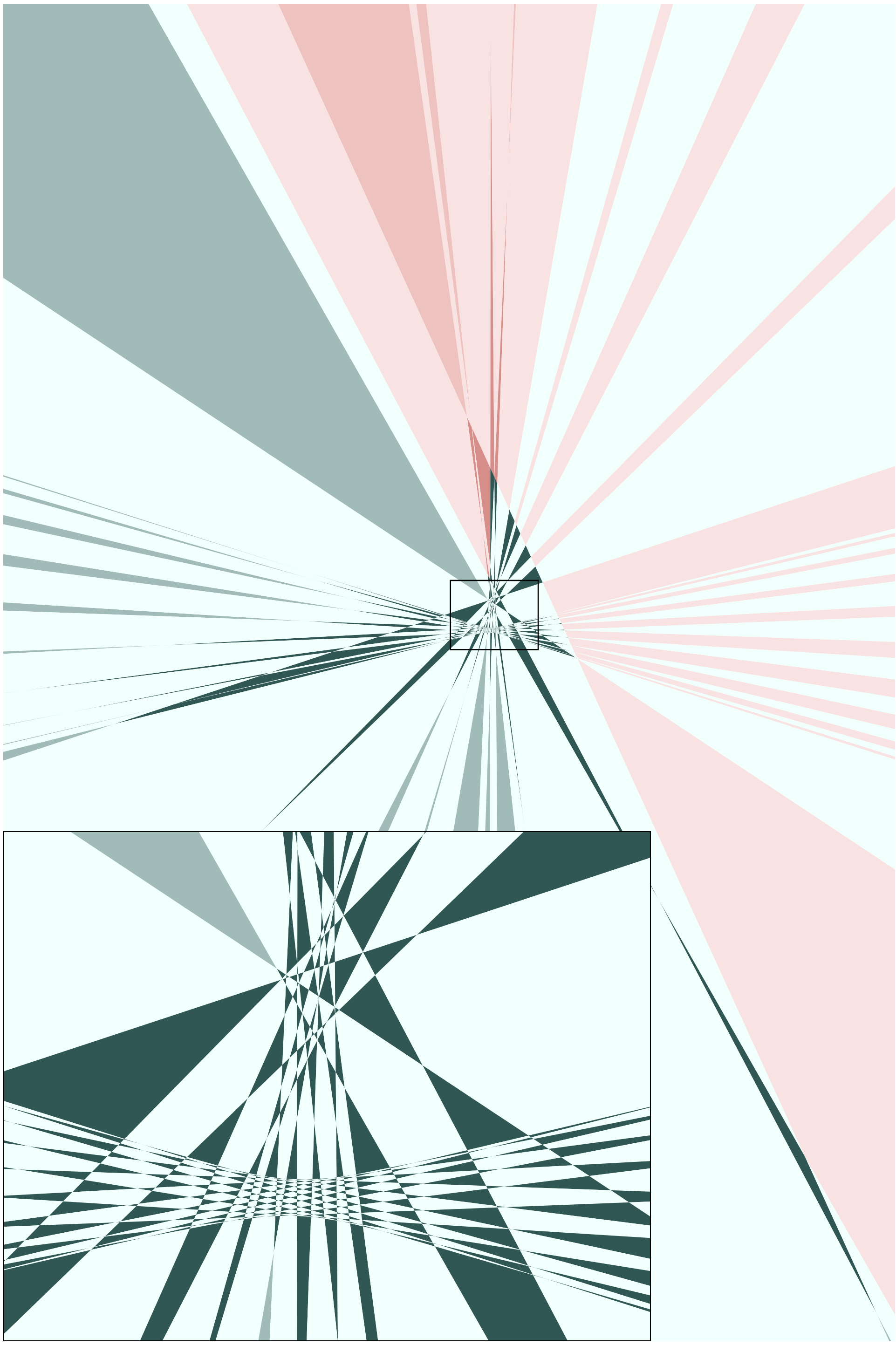

\noindent\input{figures/n38/38-1l8x1qx7on4je}\unskip\llap{\begin{minipage}[b]{0.726\linewidth}\input{figures/n38/38-1l8x1qx7on4je-zoom}\end{minipage}\hspace{0.274\linewidth}}
\caption{An arrangement of 38 lines with $a_3=450$ triangles and 3 triple points: the doubling of \cref{fig:n20-collapse} in pseudolines, stretched. Lower left: the framed region, enlarged.}
\label{fig:n38}
\end{figure}

\begin{table}[t]
\centering
\small
\caption{Equations $y=m_ix+b_i$ of the lines of the arrangement in \cref{fig:n38}. The triple points lie on line 20 together with the pairs of lines $\{12,14\}$, $\{13,17\}$, $\{15,18\}$.}
\label{tab:n38-lines}
\renewcommand{\arraystretch}{1.15}
\begin{tabular}{rrr@{\qquad}rrr}
\toprule
$i$ & $m_i$ & $b_i$ & $i$ & $m_i$ & $b_i$ \\
\midrule
0 & $-0.04793$ & $0.6523$ & 19 & $39.45675$ & $-3.68082$ \\
1 & $-0.0852$ & $-1.02173$ & 20 & $30.00493$ & $32.33556$ \\
2 & $-0.17455$ & $0.54753$ & 21 & $7.85177$ & $6.5065$ \\
3 & $-0.21998$ & $-0.99642$ & 22 & $4.83331$ & $9.25869$ \\
4 & $-0.29208$ & $0.36745$ & 23 & $4.55319$ & $7.42156$ \\
5 & $-0.32849$ & $-0.87598$ & 24 & $3.07899$ & $10.56713$ \\
6 & $-0.38084$ & $0.14532$ & 25 & $2.65403$ & $6.21769$ \\
7 & $-0.4004$ & $-0.69737$ & 26 & $1.39238$ & $10.9832$ \\
8 & $-0.42956$ & $-0.12128$ & 27 & $1.32243$ & $7.00528$ \\
9 & $-0.43625$ & $-0.4151$ & 28 & $0.45188$ & $10.11089$ \\
10 & $-0.90749$ & $9.22924$ & 29 & $0.34272$ & $0.05092$ \\
11 & $-2.40272$ & $7.63847$ & 30 & $0.3251$ & $0.28583$ \\
12 & $-2.59757$ & $15.581461275$ & 31 & $0.31861$ & $-0.22649$ \\
13 & $-2.99387$ & $52.766107044$ & 32 & $0.27275$ & $0.46675$ \\
14 & $-9.7393$ & $11.9113976453$ & 33 & $0.26314$ & $-0.5276$ \\
15 & $-10.3901$ & $17.0924954295$ & 34 & $0.18314$ & $0.59326$ \\
16 & $-12.70956$ & $5.44985$ & 35 & $0.17127$ & $-0.76651$ \\
17 & $-170.97458$ & $156.7680040263$ & 36 & $0.07339$ & $0.66667$ \\
18 & $1406.69858$ & $551.8309088275$ & 37 & $0.05342$ & $-0.93754$ \\
\bottomrule
\end{tabular}
\end{table}

\begin{table}[p]
\centering
\small
\caption{Equations $y=m_ix+b_i$ of an arrangement of 74 lines with $a_3=1764$ triangles and 3 triple points. It is obtained from the arrangement in \cref{fig:n20-collapse} by two steps of combinatorial doubling in pseudolines, analogous to \cref{sec:series}, and then stretched. This is the third member of the series $18\cdot2^t+2$. The triple points lie on line 40 together with the pairs of lines $\{32,34\}$, $\{33,37\}$, $\{35,38\}$.}
\label{tab:n74-lines}
\begin{tabular}{rrr@{\quad}rrr}
\toprule
$i$ & $m_i$ & $b_i$ & $i$ & $m_i$ & $b_i$ \\
\midrule
0 & $-6.65426468$ & $-1091.43745452$ & 37 & $0.14421861$ & $55.116043480362134$ \\
1 & $-6.64810537$ & $-1090.3864598$ & 38 & $0.14432907$ & $55.1143550503036888$ \\
2 & $-6.13520075$ & $-1003.89886196$ & 39 & $0.14480079$ & $55.20455821$ \\
3 & $-6.13374755$ & $-1003.64270684$ & 40 & $0.14497591$ & $55.21414649$ \\
4 & $-5.70667318$ & $-931.62880038$ & 41 & $0.14739226$ & $55.63755173$ \\
5 & $-5.70571304$ & $-931.45708301$ & 42 & $0.14962877$ & $56.00380269$ \\
6 & $-5.35498219$ & $-872.31688205$ & 43 & $0.15029264$ & $56.12302765$ \\
7 & $-5.35409169$ & $-872.15825542$ & 44 & $0.15480384$ & $56.85487027$ \\
8 & $-5.0732397$ & $-824.80147333$ & 45 & $0.15584749$ & $57.07155497$ \\
9 & $-5.07157073$ & $-824.51294547$ & 46 & $0.40556958$ & $98.28282147$ \\
10 & $-4.85294201$ & $-787.64847296$ & 47 & $0.4076748$ & $99.52353415$ \\
11 & $-4.85223516$ & $-787.52360093$ & 48 & $0.40962048$ & $99.85523857$ \\
12 & $-4.68721255$ & $-759.69817544$ & 49 & $0.42701422$ & $102.7788107$ \\
13 & $-4.68657182$ & $-759.58585996$ & 50 & $0.44390213$ & $105.63945342$ \\
14 & $-4.57304306$ & $-740.44325058$ & 51 & $0.48516687$ & $112.57837511$ \\
15 & $-4.57239659$ & $-740.33138239$ & 52 & $0.54937296$ & $123.42912769$ \\
16 & $-4.5057611$ & $-729.09571234$ & 53 & $0.66546018$ & $142.9647525$ \\
17 & $-4.50513696$ & $-728.98903225$ & 54 & $0.86891663$ & $177.32112487$ \\
18 & $-4.48337368$ & $-725.31945462$ & 55 & $1.47866731$ & $280.03295081$ \\
19 & $-1.5910024$ & $-237.37186659$ & 56 & $-13.67391141$ & $-2275.26702125$ \\
20 & $-0.41986925$ & $-40.02412869$ & 57 & $-13.66831944$ & $-2274.3218405$ \\
21 & $-0.41815525$ & $-39.68721726$ & 58 & $-13.29563964$ & $-2211.47504171$ \\
22 & $-0.21753219$ & $-5.89835333$ & 59 & $-13.28440237$ & $-2209.57379786$ \\
23 & $-0.20294042$ & $-3.41383947$ & 60 & $-12.60118553$ & $-2094.35988763$ \\
24 & $-0.1330755$ & $8.34706968$ & 61 & $-12.59582406$ & $-2093.44621711$ \\
25 & $-0.13299121$ & $8.37372571$ & 62 & $-11.69989116$ & $-1942.36153415$ \\
26 & $-0.10862691$ & $12.47316892$ & 63 & $-11.69402363$ & $-1941.36001817$ \\
27 & $-0.10854949$ & $12.49109486$ & 64 & $-10.71044458$ & $-1775.49620678$ \\
28 & $0.12516097$ & $51.93883116$ & 65 & $-10.70593721$ & $-1774.72250631$ \\
29 & $0.13048166$ & $52.71558026$ & 66 & $-9.72673613$ & $-1609.59816906$ \\
30 & $0.13483437$ & $53.53791167$ & 67 & $-9.72381877$ & $-1609.09187701$ \\
31 & $0.13954354$ & $54.3160826$ & 68 & $-8.80670122$ & $-1454.43810158$ \\
32 & $0.14060184$ & $54.5227052048622872$ & 69 & $-8.80451216$ & $-1454.05453181$ \\
33 & $0.14086713$ & $54.6818822830280324$ & 70 & $-7.98207725$ & $-1315.36845915$ \\
34 & $0.14291553$ & $54.8884471310340848$ & 71 & $-7.97990258$ & $-1314.98770823$ \\
35 & $0.1430402$ & $54.9155142351386022$ & 72 & $-7.26432733$ & $-1194.32251028$ \\
36 & $0.14305875$ & $54.90438146$ & 73 & $-7.26218807$ & $-1193.9484859$ \\
\bottomrule
\end{tabular}
\end{table}

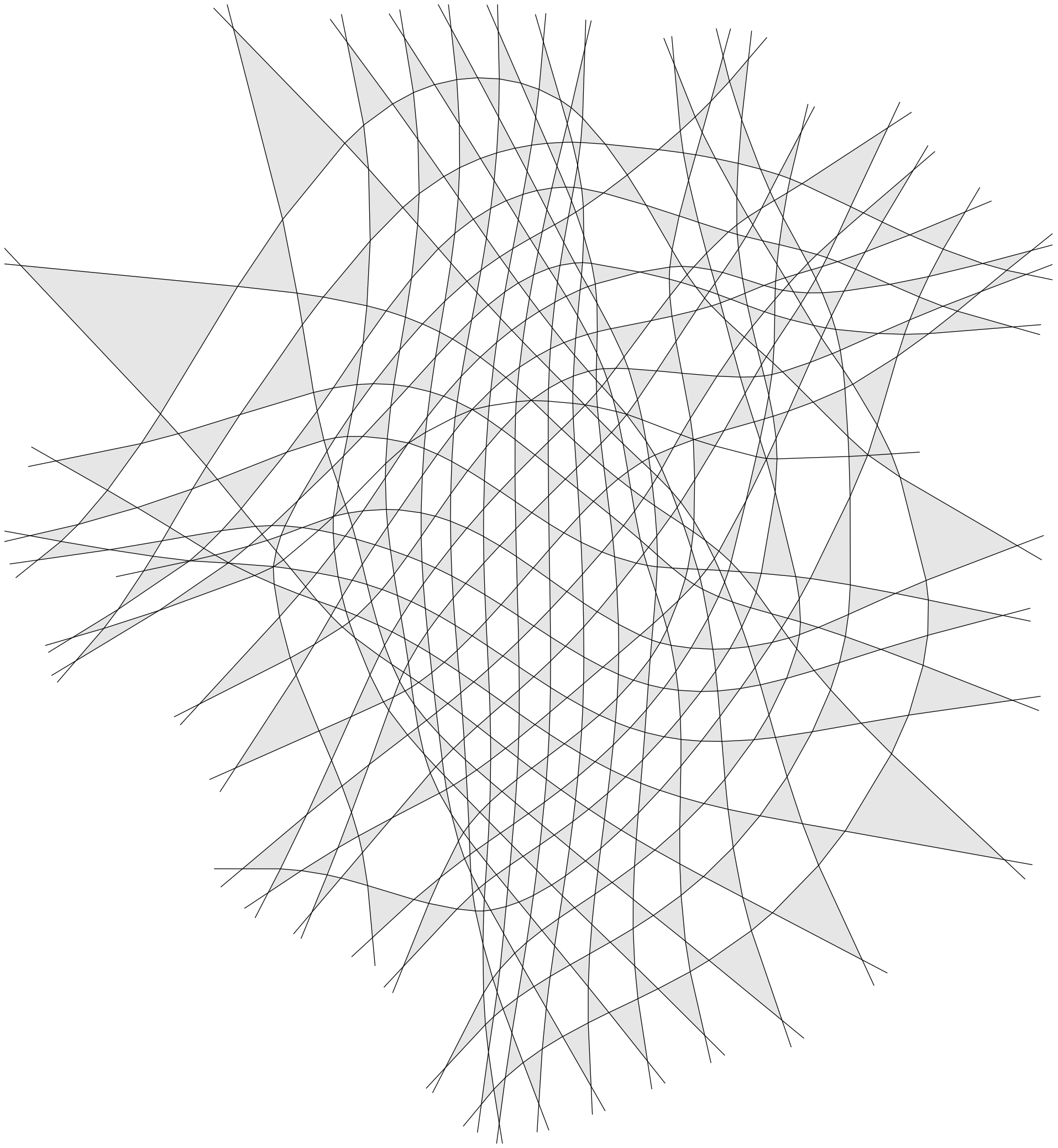
\begin{figure}[p]
\noindent\resizebox{\linewidth}{!}{\input{figures/n44/44}}

\smallskip
{\noindent\fontsize{5}{5.5}\selectfont\ttfamily 1 3 5 7 9 11 13 12 11 10 9 8 7 6 5 4 3 2 1 3 5 7 9 11 13 15 14 13 12 11 10 9 8 7 6 5 4 3 5 7 9 11 13 15 14 13 12 11 10 9 8 7 6 5 7 9 11 13 15 17 19 21 23 25 24 23 22 21 20 19 18 17 16 15 14 13 12 11 10 9 8 7 9 11 13 15 17 16 15 14 13 12 11 10 9 11 13 15 17 19 18 17 16 15 14 13 12 11 13 15 17 19 18 17 16 15 14 13 15 17 19 21 23 22 21 20 19 18 17 16 15 17 19 21 20 19 18 17 19 21 23 25 24 23 22 21 20 19 21 23 25 27 26 25 24 23 22 21 23 25 27 29 28 27 29 31 33 32 31 33 34 36 38 37 36 35 33* 36 38 40 39 38 37 36 35 34 32* 31 30 29 31* 34 36 38 40 42 41 40 39 38 37 36 35 34 33 32 31 30 29 28 27 26 25 24 23 22 21 20 19 18 17 16 15 14 13 12 11 10 9 8 7 6 5 4 3 2 1 0 1 3 5 7 9 11 13 15 17 19 21 23 25 27 29 28 27 26 25 24 23 22 21 20 19 18 17 16 15 14 13 12 11 10 9 8 7 6 5 4 3 2 1 3 5 7 9 11 13 15 17 19 21 23 25 27 29 31 30 29 28 27 26 25 24 23 22 21 20 19 18 17 16 15 14 13 12 11 10 9 8 7 6 5 4 3 5 7 9 11 13 15 17 19 21 23 25 27 29 31 34 36 35 34 33 32 31 30 29 28 27 26 25 24 23 22 21 20 19 18 17 16 15 14 13 12 11 10 9 8 7 6 5 7 9 11 13 15 17 19 21 23 25 27 29 31 33* 32 31 30 29 28 27 26 25 24 23 22 21 20 19 18 17 16 15 14 13 12 11 10 9 8 7 9 11 13 15 14 13 15 17 19 21 23 25 27 29 31 33 32 34 36 35 34 36 38 37 36 35 34 33 31* 30 29 28 27 26 25 24 23 22 21 20 19 18 17 16 15 17 19 21 23 25 27 29 31 30 32 34 36 38 40 39 38 37 36 35 34 33 32 31 29* 28 27 26 25 24 23 22 21 20 19 18 17 16 15 14 13 12 11 10 9 11 13 15 14 13 12 11 13 15 17 19 21 23 25 27 29 28 30 32 34 36 38 40 42 41 40 39 38 37 36 35 34 33 32 31 30 29 27* 26 25 24 23 22 21 20 19 18 17 16 15 14 13 12 11 10 9 8 7 6 5 4 3 2 1 3 5 7 9 8 7 6 5 4 3 5 7 9 11 10 9 8 7 6 5 7 9 11 13 15 17 16 15 17 19 18 17 19 21 20 19 18 17 19 21 23 22 21 23 25 27 26 24* 23 22 20* 19 18 16* 15 14 12* 11 10 12 14 16 15 14 13 12 14 16 18 20 22 24 23 22 24 26 28 27 26 25 24 23 22 21 20 19 18 17 16 15 14 16 18 20 22 24 26 28 30 29 28 27 26 25 24 23 22 21 20 22 24 26 28 30 32 31 30 29 28 27 26 25 24 26 28 30 32 34 33 32 31 30 29 28 27 26 25 24 23 22 24 26 28 30 32 34 36 35 34 33 32 31 30 29 28 27 26 28 30 32 34 36 38 37 36 35 34 33 32 31 30 29 28 27 26 28 30 32 34 36 38 40 39 38 37 36 35 34 33 32 31 30 29 28 30 32 34 36 38 40 42 41 40 39 38 37 36 35 34 33 32 31 30 29 28 27 26 25 24 23 22 21 20 19 18 17 16 15 14 13 12 11 9* 8 7 9 10 12 11 10 12 14 13 12 14 16 18 20 22 21 20 19 18 17 16 15 14 13 12 11 10 12 14 16 18 20 22 24 23 22 21 20 19 18 17 16 15 14 16 18 20 19 18 17 16 18 20 22 24 26 25 24 23 22 21 20 19 18 17 16 15 14 13 12 11 10 12 14 16 18 20 22 24 23 22 21 20 22 24 26 28 27 26 25 24 23 22 21 20 19 18 17 16 15 14 13 12 14 16 18 20 22 24 26 28 27 26 25 24 26 28 30 32 34 36 38 40 42\par}
\caption{An arrangement of 44 pseudolines with $a_3=609$ triangles and 12 triple points, obtained by collapse augmentation (\cref{sec:augment}). Its word is printed below the picture.}
\label{fig:n44}
\end{figure}

\clearpage

\section{Counting segments and the state table}\label{app:prunings}

This appendix justifies the two prunings announced in \cref{sec:segment-count}.
\emph{Segment pruning} (P4) discards a partial word if every completion of it has more unused segments than the balance \eqref{eq:non-simple-segm} allows for $a_3\ge K$.
\emph{State table pruning} (P5) skips a subtree if a subtree with the same state and no tighter budget constraints has already been traversed without giving words (\cref{lem:state-table}).
Along the way, \cref{app:doomed} extends the early non-triangle count (P2) to all doomed faces, which gives the doomed-face count (P6).
All three prunings are \emph{safe}: they discard only empty subtrees, so the output of the program does not change.

We keep the notation of \cref{sec:enum}: $u$ is the number of unused segments, $d$ the number of double segments (sides of two triangles at once), $t_m$ the number of vertices of multiplicity $m$, $B=\binom{n-1}2-K$ the budget of the target $K$, $S(W)$ the spending of a partial word $W$, and $a[i]$ the label of the wire in slot $i$ after $W$.
The waste $w$ \eqref{eq:waste} counts the segment units that do not go into triangles, and the target $a_3\ge K$ bounds it by the \emph{waste limit} $b$:
\begin{equation}\label{eq:waste-limit}
w(\mathcal A)\le b\overset{\mathrm{def}}{=}n(n-2)-3K .
\end{equation}

\subsection{Segment balance for arbitrary multiplicities}

In \cref{sec:enum} we already wrote down the segment balance \eqref{eq:non-simple-segm}. We now prove it in a rewritten form that uses the definition \eqref{eq:waste} of the waste $w(\mathcal{A})$.

\begin{lemma}\label{lem:segment-balance}
For every arrangement $\mathcal A$ of $n\ge3$ pseudolines,
\begin{equation}\label{eq:segm-balance}
w(\mathcal A)\overset{\mathrm{def}}{=}n(n-2)-3a_3(\mathcal A)
=u(\mathcal A)-d(\mathcal A)+\sum_v m(v)\bigl(m(v)-2\bigr).
\end{equation}
\end{lemma}

\begin{proof}
Each edge of the arrangement borders exactly two faces, which lie on opposite sides of its pseudoline, so a bounded segment is a side of at most two triangles.
The sides of a bounded face are bounded edges (a face bordering a ray is itself unbounded), and a side of a triangle is a whole edge.
Hence, if we count for each bounded segment the triangles of which it is a side, the total is $3a_3$.
Unused segments contribute 0 to this sum, double segments 2, and the others 1, so by \eqref{eq:edges}
\[
3a_3=E_b-u+d=\sum_v m(v)-n-u+d.
\]
By \eqref{eq:pair-cross}, $n(n-1)=\sum_v m(v)\bigl(m(v)-1\bigr)$, hence
\[
n(n-2)-3a_3=n(n-1)-\sum_v m(v)+u-d=u-d+\sum_v m(v)\bigl(m(v)-2\bigr).
\]
\end{proof}

We bound the right-hand side of \eqref{eq:segm-balance} by a discharging argument~\cite{cw17}.
Assign to a vertex $v$ the \emph{charge} $\mathrm{ch}(v)=m(v)\bigl(m(v)-2\bigr)$, to an unused segment the charge $+1$, to a double segment $-1$, and to the other segments~0.
By \cref{lem:segment-balance}, the charges sum to $w$.
The charge of a vertex has a direct meaning in terms of waste: each of the $m$ pseudolines through it meets $m-1$ others there at once and loses $m-2$ segments.

The search can use an upper bound on $u$ as a pruning.
For $a_3\ge K$ the waste satisfies $w\le b$ by \eqref{eq:waste-limit}, so such a bound requires only a lower bound on the total charge of vertices and double segments.
We transfer the charge of each double segment to its endpoints and bound the \emph{final charge} $\mathrm{ch}^*(v)$ of every vertex.
The transfer is possible because a double segment always has a multiple endpoint (\cref{lem:double-segment}).
A vertex of multiplicity $\ge4$ has a nonnegative final charge (\cref{lem:charge-ge4}), and a triple point has a final charge of at least $-3/2$ (\cref{cor:charge-triple}).
As a result we obtain $u\le w+\tfrac32t_3$ (\cref{cor:unused-bound}).

\subsection{Final charges of vertices}

\begin{lemma}\label{lem:double-segment}
If a bounded segment $e=PQ$ is a side of two triangles, then at least one of its endpoints has multiplicity $\ge3$.
\end{lemma}

\begin{proof}
The argument is already given in~\cite[Subsection~2.3]{pu26b} for simple arrangements, but it actually uses only the simplicity of the endpoints of $e$.
Let $e\subset\ell$, and suppose that, besides $\ell$, only the pseudolines $\ell_P$ and $\ell_Q$ pass through $P$ and $Q$, respectively.
Then a triangle with side $e$ has its sides at $P$ and $Q$ on $\ell_P$ and $\ell_Q$, and its third vertex is their common point.
Two triangles on opposite sides of $\ell$ would give two distinct common points of $\ell_P$ and $\ell_Q$.
\end{proof}

We split the double segments into two classes: $d_{\ge4}$ consists of those with at least one endpoint of multiplicity $\ge4$, and $d_3$ of all the others (by \cref{lem:double-segment}, such a segment has an endpoint of multiplicity exactly 3, and its other endpoint has multiplicity 2 or 3).
Then $d=d_3+d_{\ge4}$.

\begin{lemma}[vertices of multiplicity $\ge4$]\label{lem:charge-ge4}
Let each segment of class $d_{\ge4}$ transfer its charge $-1$ to one of its endpoints of multiplicity $\ge4$.
Then a vertex $v$ of multiplicity $m\ge4$ receives a charge of at least $-2m$, and its final charge is $\mathrm{ch}^*(v)\ge m(m-2)-2m\ge0$.
\end{lemma}
\begin{proof}
Exactly $2m$ rays and bounded edges meet at a vertex of multiplicity $m$ (two on each of the $m$ pseudolines).
The segments transferring charge to it are bounded edges, so there are at most $2m$ of them.
The inequality $2m\le m(m-2)$ holds for $m\ge4$.
\end{proof}

This does not work for a triple point: its charge is 3, while up to six double segments can be incident to it.
We therefore distribute the charge of the segments of class $d_3$, taking their other endpoints into account.

Let $P$ be a triple point.
A double segment of class $d_3$ incident to $P$ is a \emph{simple leg} of $P$ if its other endpoint is simple, and a \emph{triple leg} if its other endpoint has multiplicity 3.
Denote the number of simple legs by $\lambda_2(P)$ and the number of triple legs by $\lambda_3(P)$.
Each segment of class $d_3$ has one or two endpoints of multiplicity 3 (\cref{lem:double-segment}).
In the first case it transfers its charge $-1$ entirely to the triple endpoint, and in the second it splits the charge equally, $-\tfrac12$ to each.
The final charge of a triple point is then
\begin{equation}\label{eq:triple-charge}
\mathrm{ch}^*(P)=3-\lambda_2(P)-\tfrac12\lambda_3(P).
\end{equation}

\begin{lemma}[simple legs]\label{lem:simple-legs}
Let $h_1,\dots,h_6$ be the edges of a triple point $P$, numbered in cyclic order. No two consecutive edges are both simple legs. In particular, $\lambda_2(P)\le3$.
\end{lemma}

\begin{proof}
Indices are taken modulo 6.
Suppose, to the contrary, that $h_r$ and $h_{r+1}$ are simple legs.
Each of them is a double segment, so the three consecutive sectors (\cref{def:reduced}) adjacent to them are all triangles: the sectors between $h_{r-1}$ and $h_r$, between $h_r$ and $h_{r+1}$, and between $h_{r+1}$ and $h_{r+2}$.
The sides at $P$ of the triangular sector between $h_j$ and $h_{j+1}$ are the edges $h_j$ and $h_{j+1}$ themselves.
Let $Q_j$ be the other endpoint of $h_j$. Then the triangular sectors are $PQ_jQ_{j+1}$ ($j=r-1,r,r+1$).
Let the side $Q_rQ_{r+1}$ lie on a pseudoline $\ell$.
Besides the pseudoline containing $h_r$, only $\ell$ passes through the simple vertex $Q_r$, so the side $Q_{r-1}Q_r$ also lies on $\ell$, that is, $Q_{r-1}\in\ell$.
Similarly, $Q_{r+2}\in\ell$.
But $Q_{r-1}$ and $Q_{r+2}$ also lie on the pseudoline $\ell'$ containing the edges $h_{r-1}$ and $h_{r+2}$, which lie on opposite sides of $P$, so $\ell$ and $\ell'$ have two common points.
Moreover, $\ell\ne\ell'$: the pseudoline containing $h_r$ meets $\ell$ at $Q_r$ and $\ell'$ at $P$. This is a contradiction.
The set of simple legs contains no two consecutive elements of the 6-cycle, so it has at most three elements.
\end{proof}

\begin{corollary}\label{cor:charge-triple}
Every triple point satisfies $\mathrm{ch}^*(P)\ge-\tfrac32$.
\end{corollary}

\begin{proof}
The point $P$ is incident to six edges, so $\lambda_2(P)+\lambda_3(P)\le6$, and
\[
\mathrm{ch}^*(P)=3-\lambda_2-\tfrac12\lambda_3\ge3-\lambda_2-\tfrac12(6-\lambda_2)=-\tfrac{\lambda_2}2\ge-\tfrac32
\]
by \cref{lem:simple-legs}.
\end{proof}

The bound of this corollary is sharp. In a triangle with its three medians (six lines, four triple points), the common point of the medians has $\lambda_2=\lambda_3=3$ and a final charge of exactly $-3/2$.

\subsection{Bounding the number of unused segments}

\begin{corollary}\label{cor:unused-bound}
For every arrangement of $n\ge3$ pseudolines,
\begin{equation}\label{eq:segm-bound}
u\le w+\tfrac32 t_3 ,
\end{equation}
where $t_3$ is the number of vertices of multiplicity exactly 3.
In particular, if $a_3\ge K$, then $u\le b+\tfrac32t_3$ with the waste limit $b$ of \eqref{eq:waste-limit}.
\end{corollary}

\begin{proof}
By \cref{lem:segment-balance}, the initial charges sum to $w$.
The transfer does not change the sum, and after it the total charge of the segments equals $u$, so
\[
w=u+\sum_v\mathrm{ch}^*(v)\ \ge\ u+\sum_{P:\,m(P)=3}\mathrm{ch}^*(P)\ \ge\ u-\tfrac32t_3 ,
\]
since the final charges of vertices of multiplicity $\ge4$ are nonnegative (\cref{lem:charge-ge4}) and those of triple points are at least $-\tfrac32$ (\cref{cor:charge-triple}).
The second part follows from \eqref{eq:waste-limit}.
\end{proof}

\subsection{Doomed faces}\label{app:doomed}

\Cref{cor:unused-bound} bounds $u$ in terms of $t_3$ for a finished arrangement $\mathcal{A}$.
To use it for pruning at a partial word $W$, we need two bounds that hold for all completions of $W$: a lower bound on $u$ and an upper bound on $t_3$.
The first comes from segments each of whose two faces has closed as a non-triangle or is already doomed (\cref{lem:early-count}).
The second comes from the budget $B$ of \cref{sec:search}, strengthened by a similar argument about doomed faces (\cref{lem:triple-ceiling}).
Comparing the two bounds is the segment pruning.
Its final version also bounds the final charges of the vertices of $W$ by what is known about their edges (\cref{lem:placed-charges}).
In this subsection we describe the criteria for a face to be doomed.

In the terminology of \cref{sec:enum}, a bounded triangle is a face $i$ that was opened, received exactly one intermediate vertex and closed.

A continuation of a partial word $W$ that does not close face $i$ does not change the sets of labels in slots $0,\ldots,i$ and in slots $i+1,\ldots,n-1$ (see the proof of \cref{lem:bounded-face}).
The labels $a[i]$ and $a[i+1]$ can then change only under letters that give face $i$ intermediate vertices.
A letter with $g+m-1=i$ reverses an increasing block and puts its minimum, one of the labels in slots $0,\ldots,i$, into slot $i$. A letter with $g=i+1$ puts the maximum of its block, one of the labels in slots $i+1,\ldots,n-1$, into slot $i+1$.

\begin{lemma}\label{lem:doomed}
Suppose that in some completion of a partial word, face $i$ closes after receiving $k$ more intermediate vertices.
(a) If $a[i]>a[i+1]$, then $k\ge1$.
(b) If $a[i]>\max(a[i+1],\ldots,a[n-1])$ and $a[i+1]<\min(a[0],\ldots,a[i])$, then $k\ge2$.
\end{lemma}

\begin{proof}
At the moment of closing, the labels in the block of the closing letter increase, in particular $a[i]<a[i+1]$.
(a) Before closing, the pair $a[i],a[i+1]$ is changed only by letters that give face $i$ intermediate vertices, so there is at least one such letter.
(b) By (a), $k\ge1$, and it remains to exclude $k=1$.
If the only intermediate vertex is produced by a letter with $g=i+1$, then $a[i]$ does not change, and slot $i+1$ receives one of the labels to the right of slot $i$. By hypothesis this label is smaller than $a[i]$, and closing is impossible.
If the only intermediate vertex is produced by a letter with $g+m-1=i$, then $a[i+1]$ does not change, and slot $i$ receives one of the labels to the left of slot $i+1$. By hypothesis this label is larger than $a[i+1]$, and closing is again impossible.
\end{proof}

A face that has not closed yet is \emph{doomed} if it becomes a bounded triangle in no completion. Each of the following five criteria makes a face doomed:

\begin{enumerate}
\item face $i$ has not been opened by any letter (it is unbounded on the left);
\item the face already has $\ge2$ intermediate vertices, so it will have $\ge4$ sides if it closes, and it is unbounded if it does not;
\item the labels in slots $0,\ldots,i$ are already final, so the face will not close, that is, it is unbounded (\cref{lem:bounded-face});
\item the face has exactly one intermediate vertex and $a[i]>a[i+1]$: at least one more appears before closing (\cref{lem:doomed}(a)), so there will be $\ge2$ when it closes;
\item the face has no intermediate vertices, but $a[i]$ exceeds all labels to the right of slot $i$ and $a[i+1]$ is less than all labels to the left of slot $i+1$: by \cref{lem:doomed}(b), at least two appear before closing.
\end{enumerate}

The first three criteria are checked against the status of the face and the labels to its left, the fourth against the pair of wires in slots $i$ and $i+1$, and the fifth against this pair and all labels.
Being doomed is defined over all completions, so a face stays doomed until it closes, even if the criterion that established it stops holding (for example, criterion 5 may hold for $W$, criterion 4 for an extension of $W$, and criterion 2 for a further extension).

The early non-triangle count (P2) already counts, through the term $\delta S(W)$, the faces that have at least two intermediate vertices.
This argument uses only the fact that such a face will close and will not become a triangle, so it extends to every doomed bounded face.
Let $\delta^+S(W)\ge\delta S(W)$ be the number of bounded faces doomed by criteria 2, 4 and 5 of the list above (criteria 1 and 3 mean that the face is unbounded).
In modes II--IV (\cref{tab:modes}) the search discards a branch as soon as $S(W)+\delta^+S(W)>B$. This is the doomed-face count (P6).

\subsection{Segment pruning}

\begin{lemma}[early counting of unused segments]\label{lem:early-count}
In each of the following two cases, a segment $e$ of a partial word $W$ is, in every completion of $W$, bounded and used by no triangle.
(a) The segment $e$ is already formed (both its endpoints are crossings produced by letters of $W$), and each of the two faces bordering it is doomed or has closed as a non-triangle.
(b) The segment $e$ is the \emph{current segment} of the wire in slot $p$ (the part of this wire after its last crossing in $W$), faces $p-1$ and $p$ are doomed (for $p=0$ or $p=n-1$, one of them is the outer unbounded region, which never becomes a bounded triangle), and the wire already has a crossing in $W$ and still has an uncrossed pair, that is, some slot to the left of $p$ holds a smaller label or some slot to the right of $p$ holds a larger one.
\end{lemma}

\begin{proof}
(a) Neither of the two faces of a formed segment becomes a bounded triangle, and both its endpoints are crossings of $W$.
(b) The segment of the wire in slot $p$ borders the current faces $p-1$ and $p$, and they do not change until the next crossing on this wire: every letter that closes face $p-1$ or $p$ contains slot $p$ in its block.
Both faces are doomed, so neither becomes a bounded triangle, and $e$ is a side of no triangle.
Every pair of wires crosses exactly once, so in every completion the wire has one more crossing, and the current segment is bounded.
\end{proof}

Let $u_{\mathrm{lb}}(W)$ be the number of segments that \cref{lem:early-count} shows to be unused. Every completion of $W$ satisfies $u\ge u_{\mathrm{lb}}(W)$.

\begin{lemma}[triple-point ceiling]\label{lem:triple-ceiling}
Let $P(W)$ be the number of wire pairs already crossed in a partial word $W$.
Then every completion with $a_3\ge K$ has at most $t_3(W)+f$ triple points, where $t_3(W)$ is the number of triple points of $W$ and
\[
f=\min\left\{\left\lfloor\frac{\binom n2-P(W)}3\right\rfloor,\
B-S(W)-\delta^+S(W)\right\}.
\]
\end{lemma}

\begin{proof}
Each future triple point uses exactly three of the pairs of wires that have not crossed yet.
Vertices of higher multiplicity use more pairs.
This gives the first bound.
Further, each future triple point adds $\binom{3-1}2=1$ to $S$, and the closing of the faces counted in $\delta^+S(W)$ adds another $\delta^+S(W)$ to $S$.
A word with $a_3\ge K$ has $S\le B$, which gives the second bound.
\end{proof}

\Cref{cor:unused-bound} together with \cref{lem:early-count,lem:triple-ceiling} already allows pruning at a partial word $W$: every completion of $W$ with $a_3\ge K$ satisfies $u_{\mathrm{lb}}(W)\le u\le b+\tfrac32\bigl(t_3(W)+f\bigr)$, and every search branch that violates this inequality can be discarded.
However, this bound assumes a final charge of $-3/2$ at every triple point and $0$ at every vertex of multiplicity $\ge4$, including the vertices of $W$.
A charge of $-3/2$ occurs only with six double segments, and in a partial word many edges are already known not to be double. We therefore bound the final charges of the vertices of $W$ by what is known about their edges and leave the pessimistic case to future vertices.

A segment is \emph{decided} as non-double as soon as one of its two faces has closed as a non-triangle or is doomed by any of the five criteria, including a face doomed by the current letter itself, and as double when both faces have closed as triangles.
Both decisions are final.

Let $v$ be a vertex of multiplicity $m\ge3$ in a partial word $W$, and let $\mathrm{und}(v)$ be the number of undecided segments incident to it.
For a triple point, let $\bar\lambda_2(v)$ and $\bar\lambda_3(v)$ be the numbers of incident segments decided as double whose other endpoint is simple or has multiplicity 3, respectively. In every completion the former are simple legs and the latter triple legs of $v$.
Number the six edges of $v$ cyclically as in \cref{lem:simple-legs}. Consecutive positions bound one sector.
A set of positions is \emph{independent} if it contains no two consecutive ones. By \cref{lem:simple-legs}, the simple legs form an independent set.
Let $D$ be the positions of the segments counted in $\bar\lambda_2(v)$ as future simple legs, $U$ the positions of the undecided segments, and let
\[
\xi(v)=\max\bigl\{|I|:\ I\subseteq U,\ D\cup I\ \text{is independent}\bigr\}
\]
be the largest number of undecided segments that can still become simple legs.
For a vertex of multiplicity $m\ge4$, let $\bar\lambda(v)$ be the number of segments decided as double that transfer their charge to $v$ under the rule of \cref{lem:charge-ge4}.
A segment whose endpoints both have multiplicity $\ge4$ transfers its charge to the endpoint that occurs earlier in the word.
Let
\begin{equation}\label{eq:charge-lb}
\mathrm{ch}^*_{\mathrm{lb}}(v)=\begin{cases}
3-\bar\lambda_2-\tfrac12\bigl(\bar\lambda_3+\mathrm{und}+\xi\bigr), & m=3,\\[2pt]
m(m-2)-\bar\lambda-\mathrm{und}, & m\ge4.
\end{cases}
\end{equation}

\begin{lemma}[charges of the vertices of $W$]\label{lem:placed-charges}
In every completion of a partial word $W$, each vertex $v$ of $W$ of multiplicity $\ge3$ has $\mathrm{ch}^*(v)\ge\mathrm{ch}^*_{\mathrm{lb}}(v)$, and $\mathrm{ch}^*_{\mathrm{lb}}(v)\ge-\tfrac32$ if $m(v)=3$.
Hence every completion of $W$ with $a_3\ge K$ satisfies
\begin{equation}\label{eq:charge-cut}
2u\ \le\ 2b-\sum_{\substack{v\in W\\ m(v)\ge3}}
2\,\mathrm{ch}^*_{\mathrm{lb}}(v)+3f .
\end{equation}
\end{lemma}
\begin{proof}
Decided segments transfer to $v$ the same charge as in the completion: 0 if non-double, and their final share if double.
Let $m(v)=3$. An undecided segment decreases $2\,\mathrm{ch}^*(v)$ by 2 if it becomes a simple leg and by 1 if it becomes a triple leg, and it does not change $2\,\mathrm{ch}^*(v)$ otherwise (then it is not double, or it transfers its charge to its other endpoint of multiplicity $\ge4$).
By \cref{lem:simple-legs}, the simple legs of the completion form an independent set containing $D$, so at most $\xi$ undecided segments become simple legs.
Hence $2\,\mathrm{ch}^*(v)\ge6-\bigl(2\bar\lambda_2+\bar\lambda_3+\mathrm{und}+\xi\bigr)=2\,\mathrm{ch}^*_{\mathrm{lb}}(v)$.
Since $\bar\lambda_2+\bar\lambda_3+\mathrm{und}\le6$ and $\bar\lambda_2+\xi=|D\cup I|\le3$ for a set $I$ attaining the maximum, the subtracted sum is at most 9, that is, $\mathrm{ch}^*_{\mathrm{lb}}(v)\ge-\tfrac32$.
Let $m(v)\ge4$. Charge is transferred to $v$ by the segments counted in $\bar\lambda(v)$ and by some of the undecided segments, so $\mathrm{ch}^*(v)\ge m(m-2)-\bar\lambda-\mathrm{und}$.

As in the proof of \cref{cor:unused-bound}, $u=w-\sum_v\mathrm{ch}^*(v)$, where the sum runs over the vertices of multiplicity $\ge3$ in the completion.
The charges of the vertices of $W$ are bounded by the first part of the lemma.
There are at most $f$ future triple points (\cref{lem:triple-ceiling}), and the final charge of each is at least $-\tfrac32$ (\cref{cor:charge-triple}). The final charges of future vertices of multiplicity $\ge4$ are nonnegative (\cref{lem:charge-ge4}).
Together with $w\le b$, this proves the lemma.
\end{proof}

\emph{Segment pruning} discards a partial word $W$ if $2u_{\mathrm{lb}}(W)$ exceeds the right-hand side of \eqref{eq:charge-cut}. It is safe: a completion with $a_3\ge K$ would have $u\ge u_{\mathrm{lb}}(W)$, contrary to \eqref{eq:charge-cut}.

\subsection{The state table}

During the traversal of the tree, different partial words $W$ and $W'$ may share continuations. Which continuations can follow and how many triangles they give depends not on the whole partial word but on the state accumulated by applying its letters.
If the search found no completion of a partial word $W$ with $a_3\ge K$, then no other partial word $W'$ with the same state, or with a state that differs only by larger spending, has one either. The fact that a state has no completion is stored in a special hash table and used to skip empty subtrees.

The \emph{state} of a partial word $W$ is the tuple
\[
\sigma(W)=\bigl(\mathbf a,\ \widehat g,\ \mathbf c,\ S\bigr),
\]
where $\mathbf a=(a[0],\ldots,a[n-1])$ is the permutation of wire labels over the slots, $\widehat g$ is the initial slot of the last letter, $\mathbf c=(c[0],\ldots,c[n-2])$ is the vector of codes of the current faces, and $S=S(W)$ is the spending.

The possible continuations of a partial word, that is, those satisfying conditions (i) and (iii) of \cref{sec:search} and the forced descent (\cref{lem:forced-descent}), are determined by the pair $(\mathbf a,\widehat g)$: the increase of labels in the block of a letter is checked against $\mathbf a$, condition (iii) against $\widehat g$, and the forced descent against both.
Hence partial words with the same pair $(\mathbf a,\widehat g)$ have the same possible continuations, and it remains to compare them by face codes and spending.

The code of face $i$ records its status, and its \emph{cost} is the contribution of its closing to $S$ if the face receives exactly $k$ more intermediate vertices before closing (the indicator $\mathbf 1$ equals 1 if the condition in its subscript holds and 0 otherwise):

\begin{center}
\small
\setlength{\tabcolsep}{4pt}
\begin{tabular}{>{\raggedright\arraybackslash}p{6.5cm}cc}
\toprule
status of face $i$ & code $c[i]$ & cost \\
\midrule
face not opened (initial unbounded face) & 0 & 0 \\
face opened, no intermediate vertices & 1 & $\mathbf 1_{k\ne1}$ \\
one intermediate vertex & 2 & $\mathbf 1_{k\ne0}$ \\
$\ge2$ intermediate vertices & 3 & 1 \\
\bottomrule
\end{tabular}
\end{center}

A closed bounded face adds one to $S$ if its number of intermediate vertices differs from one, and zero otherwise. A face with code 0 is unbounded and does not contribute to $S$.
If the face will not close, its cost is zero for every code.
A letter that closes the current face $i$ opens a new face $i$.
The new face has no intermediate vertices regardless of $W$, so the contributions of the later faces $i$ to $S$ depend only on the continuation, not on $W$ itself.

We need two properties of the codes.
(a) Code 0 is determined by $\mathbf a$ alone: it means that the initial face $i$ has not closed yet, that is, slots $0,\ldots,i$ hold the labels $0,\ldots,i$.
Indeed, before face $i$ first closes, the set of labels in these slots does not change (see the proof of \cref{lem:bounded-face}), and the closing letter crosses the wires in slots $i$ and $i+1$, with labels $w_1\le i<w_2$. They cannot cross a second time, so the wire $w_2$ stays to the left of $w_1$, and the labels $0,\ldots,i$ never again occupy slots $0,\ldots,i$.
Hence partial words with the same permutation $\mathbf a$ have code 0 at the same faces, and only codes 1--3 can differ.
(b) Code 1 implies $a[i]>a[i+1]$: the letter that opened face $i$ reversed an increasing block, after which the pair $a[i],a[i+1]$ is decreasing. Only an intermediate vertex could change this, and with code 1 there was none.

\emph{Cost of a discrepancy.} Let partial words $X$ and $Y$ with the same permutation $\mathbf a$ have different codes of face $i$, $x=c_X[i]\ne y=c_Y[i]$, both at least 1.
They have the same number $k$ of new intermediate vertices, since the intermediate vertices of face $i$ before its closing are produced by letters of the common continuation.
Let $\Delta S(x,y,k)$ be the difference between the cost of face $i$ for $X$ and for $Y$.
For $x>y$ the cost table gives
\[
\begin{gathered}
\Delta S(2,1,k)=\mathbf 1_{k\ne 0}-\mathbf 1_{k\ne 1}=\mathbf 1_{k=1}-\mathbf 1_{k=0},\\
\Delta S(3,1,k)=1-\mathbf 1_{k\ne 1}=\mathbf 1_{k=1},\qquad
\Delta S(3,2,k)=1-\mathbf 1_{k\ne 0}=\mathbf 1_{k=0}.
\end{gathered}
\]
The exact number $k$ of new intermediate vertices is not known in advance, but \cref{lem:doomed} bounds it in terms of the permutation $\mathbf a$: $k\ge2$ under the hypothesis of its part (b), and $k\ge1$ if $a[i]>a[i+1]$.
For $k\ge2$ all values of $\Delta S$ vanish, and for $k\ge1$ the value $\Delta S(3,2,k)$ vanishes.
The properties of indicators also show that all three values of $\Delta S$ are at most one.

Hence for $x>y$ we have $\Delta S(x,y,k)\le\rho_i$ for every continuation, where
\begin{equation}\label{eq:face-price}
\rho_i=\begin{cases}
0,&\text{if }a[i]>\max(a[i+1],\ldots,a[n-1])\text{ and }a[i+1]<\min(a[0],\ldots,a[i]),\\
0,&\text{if }(x,y)=(3,2)\text{ and }a[i]>a[i+1],\\
1&\text{otherwise.}
\end{cases}
\end{equation}
The conditions in \eqref{eq:face-price} are determined by the common permutation $\mathbf a$ and the pair of codes $(x,y)$, so $\rho_i$ is the same quantity for $X$ and $Y$.

For $x<y$ the difference is nonpositive for every continuation: $\Delta S(1,3,k)=\mathbf 1_{k\ne1}-1$, $\Delta S(2,3,k)=\mathbf 1_{k\ne0}-1$, and $\Delta S(1,2,k)=\mathbf 1_{k\ne1}-\mathbf 1_{k\ne0}$, where $\mathbf 1_{k\ne0}=1$, since code 1 for $X$ implies $a[i]>a[i+1]$ (property (b)) and hence $k\ge1$ (\cref{lem:doomed}(a)).

\begin{lemma}[dominance]\label{lem:dominance}
Let partial words $X$ and $Y$ have the same $\mathbf a$ and $\widehat g$, and let
\[
h(X,Y)=\sum_{i:\ c_X[i]>c_Y[i]}\rho_i ,
\]
where the sum runs over the faces $i$ that will still close (by \cref{lem:bounded-face}, this is determined by $\mathbf a$).
If $S(X)+h(X,Y)\le S(Y)$, then for every completion $YV$ of $Y$, the word $XV$ is a completion of $X$, and $a_3(XV)\ge a_3(YV)$.
\end{lemma}

\begin{proof}
A sequence $V$ is a possible continuation of $X$ if and only if it is a possible continuation of $Y$, since $X$ and $Y$ share the pair $(\mathbf a,\widehat g)$.
By \cref{cor:budget}, $a_3=\binom{n-1}2-S$, so it suffices to check that $S(XV)-S(YV)\le0$.
The future contributions of multiple points to $S$ are the same, since they come from the letters of $V$.
The contributions of non-triangular faces differ only in the closing of the current faces $i$, and a face that will not close (criterion 3) contributes zero for both words.
For a face $i$ that will close, the difference of contributions is zero for equal codes and otherwise equals $\Delta S(c_X[i],c_Y[i],k)$, which is nonpositive for $c_X[i]<c_Y[i]$ and at most $\rho_i$ for $c_X[i]>c_Y[i]$.
Summing over $i$, we obtain $S(XV)-S(YV)\le S(X)-S(Y)+h(X,Y)\le0$.
\end{proof}

During the traversal of the tree, the state table is filled on backtracking: the state of a partial word $X$ is recorded when all its possible continuations have been traversed completely and none of them gave a word with $a_3\ge K$.
A node $Y$ is skipped if the table has an entry $X$ with the same $(\mathbf a,\widehat g)$ such that $X$ and $Y$ satisfy the condition of \cref{lem:dominance}; we say that $X$ \emph{dominates} $Y$.
The key $(\mathbf a,\widehat g)$ is compared in full. Eviction of entries from the hash table reduces the efficiency of the pruning but does not lose words.

\begin{lemma}[the state table loses no words]\label{lem:state-table}
The subtree of any partial word recorded in the state table contains no words with $a_3\ge K$.
Hence skipping by the state table loses no words.
\end{lemma}

\begin{proof}
We argue by induction on the order of entries. Suppose the lemma holds for all entries made before the entry of a partial word $X$; we prove it for the entry of $X$.
Suppose, to the contrary, that $X$ is recorded in the state table although its subtree contains a word $XV$ with $a_3(XV)\ge K$.
By the filling rule, $X$ is recorded only when the traversal of its subtree has found no such words, so the traversal did not reach $XV$: some partial word $Z$ on the path from $X$ to $XV$ was discarded.
Let $XV=ZV'$.
No other pruning could have discarded $Z$: all of them are safe, and $ZV'$ is a word with $a_3\ge K$.
Hence $Z$ was skipped by the entry of some $X'$ in the state table.
This entry was made before the entry of $X$, because $Z$ was processed during the traversal of the subtree of $X$, and $X$ is recorded after this traversal ends (if $X$ is the first entry, this case cannot occur at all).
By \cref{lem:dominance}, $X'V'$ is a completion of $X'$ and $a_3(X'V')\ge a_3(ZV')\ge K$, which contradicts the induction hypothesis for $X'$.
\end{proof}

\begin{remark}\label{rem:shared-table}
The proof of \cref{lem:state-table} does not use the particular order of tree traversal, only the order in which the entries appear in time. The lemma therefore also holds for a state table shared by several concurrently running processes (provided the implementation guarantees that writes are atomic), and for a state table accumulated by an interrupted run with the same $n$, $K$ and set of prunings when this run is restarted.
\end{remark}

Segment pruning and the counting of $\delta^+S$ are not part of the table key: in the proof of \cref{lem:state-table} they, like all safe prunings, discard only subtrees without words with $a_3\ge K$.

\subsection{The state table for reduced arrangements}

With the pruning of non-reduced arrangements (P7, \cref{sec:segment-count}), \cref{lem:state-table} does not apply. An entry $X$ then only means that the subtree of $X$ contains no words of reduced arrangements with $a_3\ge K$, while a continuation $V$ that gives a non-reduced arrangement after $X$ may give a reduced one after a partial word $Y$ dominated by $X$.
\Cref{cor:reduced} guarantees reducedness only for arrangements with the largest number of vertices, so skipping by the table has to compare the number of vertices as well.
Let $\ell(W)$ be the number of letters (vertices) of a partial word $W$, and let $\ell_{\max}$ be the largest number of vertices of an arrangement with $a_3\ge K$.

\begin{lemma}[the state table for reduced arrangements]\label{lem:state-table-reduced}
Suppose that $\ell(X)$ is recorded in the table together with the state of a partial word $X$, and a node $Y$ is skipped by the entry of $X$ only if both the condition of \cref{lem:dominance} and $\ell(X)\ge\ell(Y)$ hold.
Then the subtree of any recorded partial word contains no words with $a_3\ge K$ and $\ell_{\max}$ letters.
Hence the search with the pruning of non-reduced arrangements emits all words of the tree with $a_3\ge K$ and $\ell_{\max}$ letters, and if an arrangement with $a_3\ge K$ exists, the output is nonempty.
\end{lemma}

\begin{proof}
We repeat the proof of \cref{lem:state-table} for words $XV$ with $a_3(XV)\ge K$ and $\ell(XV)=\ell_{\max}$. It changes in two places.
Besides the safe prunings, the partial word $Z$ could also have been discarded by the pruning of non-reduced arrangements. But the arrangement represented by $ZV'$ has the largest number of vertices among arrangements with $a_3\ge K$, so by \cref{cor:reduced} it is reduced, and this pruning discards none of its words.
For the entry $X'$ we need $\ell(X'V')=\ell_{\max}$ in addition to $a_3(X'V')\ge K$. It follows from $\ell(X')\ge\ell(Z)$, since then $\ell(X'V')\ge\ell(ZV')=\ell_{\max}$.
The same arguments show that the search reaches every word of the tree with $a_3\ge K$ and $\ell_{\max}$ letters and emits it, and by \cref{cor:reduced} such a word with first letter $(0,2)$ exists if an arrangement with $a_3\ge K$ exists.
\end{proof}

\Cref{rem:shared-table} on concurrent processes and restarts carries over to \cref{lem:state-table-reduced} without change.

\clearpage
\bibliographystyle{unsrt}
\bibliography{references}

\end{document}